\newif\ifpictures
\picturestrue

\newif\ifcomment
\commenttrue

\documentclass[12pt]{amsart}
\usepackage{scrextend}
\usepackage{hyperref}
\usepackage{latex_base}
\usepackage{macros}
\usepackage{longtable}
\usepackage[tableposition=top]{caption}
\usepackage{tablefootnote}

\usepackage[hyphenbreaks]{breakurl}

\usepackage{pgfplots} %for the pictures
\usepackage{pgfplotstable}
\pgfplotsset{compat=1.9}
\usepackage{xcolor}
\definecolor{tubsRed}{HTML}{be1e3c} 
\definecolor{tubsGreen}{HTML}{e16d00}%{6d8300} 
\definecolor{tubsLightGray}{cmyk}{0,0,0,0.1}
\definecolor{b0}{RGB}{0, 0, 0}			%black
\definecolor{b4}{RGB}{68, 119, 170}		%blue
\definecolor{b5}{RGB}{102, 204, 238}	%cyan
\definecolor{b3}{RGB}{34, 136, 51}		%green
\definecolor{b2}{RGB}{204, 187, 68}		%yellow
\definecolor{b1}{RGB}{238, 102, 119}	%red
\definecolor{b6}{RGB}{170, 51, 119}		%purple
\definecolor{b7}{RGB}{187, 187, 187}	%grey

\usepackage{enumitem} %for referencing enumerate items
\usepackage{subcaption}
\usepackage{graphicx}
\usepackage{tikz}
\usepackage{csquotes}
\usepackage{csvsimple}
\usepackage{multirow}
\usepackage{cancel}
\algrenewcommand\algorithmicrequire{\textbf{Input:}}
\algrenewcommand\algorithmicensure{\textbf{Output:}}

\newcommand{\tilted}{tilted } %relaxed/relaxation
\newcommand{\homogeneous}{homogeneous }
\author{Sabrina C.L. Ammann}

\address{Sabrina C.L. Ammann, Technische Universit\"at Braunschweig, Institut f\"ur Mathematische Optimierung, Universit\"atsplatz 2, 38106 Braunschweig,
 Germany\medskip}
 
\email{sabrina.ammann@tu-braunschweig.de}

\author{Birte Ostermann}

\address{Birte Ostermann, Technische Universit\"at Braunschweig, Institut f\"ur Analysis und Algebra, AG Algebra, Universit\"atsplatz 2, 38106 Braunschweig,
 Germany\medskip}

\email{birte.ostermann@tu-braunschweig.de}

\author{Sebastian Stiller}

\address{Sebastian Stiller, Technische Universit\"at Braunschweig, Institut f\"ur Mathematische Optimierung, Universit\"atsplatz 2, 38106 Braunschweig,
 Germany\medskip}
 
\email{sebastian.stiller@tu-braunschweig.de}

\author{Timo de Wolff}

\address{Timo de Wolff, Technische Universit\"at Braunschweig, Institut f\"ur Analysis und Algebra, AG Algebra, Universit\"atsplatz 2, 38106 Braunschweig,
 Germany\medskip}

\email{t.de-wolff@tu-braunschweig.de}


\title[A speed-up for Helsgaun's TSP heuristic]{A speed-up for Helsgaun's TSP heuristic by relaxing the positive gain criterion}
\begin{document}

\maketitle
\begin{abstract}
The Traveling Salesman Problem (TSP) is one of the most extensively researched and widely applied combinatorial optimization problems.
It is NP-hard even in the symmetric and metric case.
Building upon elaborate research, state-of-the-art exact solvers such as CONCORDE can solve TSP instances with several ten thousand vertices.
A key ingredient for these integer programming approaches are fast heuristics to find a good initial solution, in particular the Lin-Kernighan-Helsgaun (LKH) heuristic.
For instances with few hundred vertices heuristics like LKH often find an optimal solution.
In this work we develop variations of LKH that perform significantly better on large instances.
LKH repeatedly improves an initially random tour by exchanging edges along alternating circles.
Thereby, it respects several criteria designed to quickly find alternating circles that give a feasible improvement of the tour.
Among those criteria, the positive gain criterion stayed mostly untouched in previous research.
It requires that, while constructing an alternating circle, the total gain has to be positive after each pair of edges.
We relax this criterion carefully leading to improvement steps hitherto undiscovered by LKH.
We confirm this improvement experimentally via extensive simulations on various benchmark libraries for TSP.
Our computational study shows that for large instances our method is on average 13\% faster than the latest version of LKH.
\end{abstract}

\section{Introduction}
The \struc{Traveling Salesman Problem (TSP)} is one of the most researched combinatorial optimization problems. It has found numerous and widespread applications from operations research to genome sequencing \cite[Chapter 2]{cook2011traveling}.
The TSP is strongly NP-hard even in the metric case \cite{Karp1972}.
Thus, unless P=NP, there exists no polynomial time algorithm. 
However, state-of-the-art exact solvers such as \struc{CONCORDE} \cite{concorde} find exact solutions for instances with tens of thousands of vertices. 
Metric instances of a few hundred vertices can be solved by CONCORDE on a smart phone within few seconds.

Solvers like CONCORDE are based on highly developed integer programming methods, many of which are specialized for the TSP. 
One part of these methods are heuristics to find a good starting solution.
The quality of the starting solution influences the running time and success rate of the subsequent integer programming algorithm.
For example, small
metric instances are solved so rapidly because the start heuristic used by CONCORDE --  a chained variant of the Lin-Kernighan heuristic \cite{LK_original} \struc{(LK)}  --
finds an optimal tour in almost all small metric instances. Then the subsequent integer
programming only has to prove its optimality.

Arguably, the most highly developed variant\footnote{Applegate et al. consider LKH as \enquote{milestone in the evolution of the Lin-Kernighan heuristic}, \cite[p. 466]{cook2011traveling}. Especially, LKH has repeatedly improved and holds the current record for the World TSP, see \url{https://www.math.uwaterloo.ca/tsp/world/}, last accessed 01/24, as well as records for many of the standard benchmark instances \cite{TSPLIB,vlsi,DIMACS,nationalTSP,Hougardy_Tnm}.} of the heuristic of Lin and Kernighan
\cite{LK_original} (LK) 
is that of Helsgaun \cite{LKH_original} \struc{(LKH)}.
The basic idea of these heuristics is to repeatedly improve the cost of a tour by means of edge exchange along an alternating circle.
This means that one selects a set $X$ of edges from the existing tour to be removed and a suitable set $Y$ of other edges to replace them at lower cost.
The edges in $X$ and $Y$ together form a circle that alternates between the edges from $X$ and $Y$.

Every tour $T$ can be transformed into every other tour $T’$ for the same instance by a finite set of edge exchanges along \struc{alternating circles}.
This holds because the symmetric difference for any $T$ and $T'$ can be decomposed into alternating circles.
Thus, an algorithm that tests all sets of exchanges along alternating circles finds an optimal tour. Clearly, this algorithm is neither polynomial nor practical.
The heuristic of Lin and Kernighan restricts the set of alternating circles it considers in several ways.
Thereby, LK becomes a very fast but heuristic method, i.e., it does not necessarily find an optimal solution.
Helsgaun's variants (LKH) of LK relax some of these restrictions, trading higher computational effort for stronger possibilities to improve the cost of a tour.

One of the restrictions not relaxed in LKH is the \struc{positive gain criterion}\footnote{By default, in Helgaun's code the positive gain criterion is active. It can be disabled by the user.}.
The LK and LKH heuristics build the alternating circle $(x_1, y_1, \ldots, x_k, y_k)$ edge by edge, where the edges $x_i$ will be removed from the current tour and the edges $y_i$ will be added. 
If the total cost of the removed edges is strictly larger than that of the added edges, then the exchange along the circle $(x_1, y_1,\ldots, x_i, y_i,\ldots, x_k, y_k)$  yields an improved tour. 
The positive gain criterion used by LK and LKH requires that this difference is also strictly positive for all paths of the form $P_i := (x_1, y_1, \ldots, x_i, y_i)$. Requiring positive gain also for\linebreak the paths $P_i$ with $1\leq i < k$ shall guide the search towards an alternating circle with strong cost reduction.
This guidance is heuristic as there are alternating circles that lead to cost reduction in the end, but do not fulfill the positive gain criterion at each step $i$.\linebreak
Moreover, for every alternating tour that is beneficial in the end, there is a cyclic permutation such that also every partial sum is positive, i.e., the positive gain criterion is fulfilled, cf. \cite{LK_original}. 
That means for every beneficial alternating circle there is at least one that fulfills the positive gain criterion. 
On the one hand, this observation is a strong justification\linebreak for the positive gain criterion not being relaxed.
On the other hand, the existence of a beneficial circle that fulfills the positive gain criterion does not ensure, that the heuristic finds it.
Together with other criteria, the positive gain criterion might even strictly inhibit finding certain beneficial circles.
One of these other criteria that is particularly inhibitive is that LKH only considers edges $y_i$ from candidate sets pre-computed for each vertex.

The main question to be addressed in this paper is, whether it is beneficial in terms of computation time and objective value to search for alternating circles that partially violate the positive gain criterion. 

\subsection{Our contribution}

We propose a variant of Helsgaun's heuristic with an enlarged space for searching alternating circles.
Via mildly relaxing the positive gain criterion. 
We allow pairs of edges $x_i$ and
$y_i$ in the alternating circle --- where again $x_i$ is an edge to be removed from the original tour and $y_i$ one to be added to it --- such that the total gain of exchanging along the alternating circle up to $x_i$ and $y_i$ is not positive.
Still, we are rather restrictive in our relaxation, as we do not allow for two subsequent pairs of edges $x$ and $y$ to violate the positive gain criterion.

We give an extensive computational study over established benchmark instances to show that the relaxed positive gain criterion yields a superior heuristic. Starting from this initial success we pursue further improvements by algorithm engineering. 
This results in a further variation of the relaxed positive gain criterion that yields an additional significant improvement in light of our computational experiments.

In particular, we compare our implementation to LKH-3 \cite{LKH3} on 438 benchmark instances with less than one million vertices \cref{tab:overview_instances} that originate from the TSPLIB \cite{TSPLIB}, the VLSI dataset \cite{vlsi}, the eighth DIMACS challenge \cite{DIMACS}, the national TSPs \cite{nationalTSP}, and the $T_{nm}$ instances \cite{Hougardy_Tnm}.

For large instances our heuristic with relaxed positive gain criterion reduces the computation time on average by 3.4 hours or 13.6\% \cref{fig:CR_large_instances}, and for individual instances up to 31\% \cref{fig:CR_large_instances_POPMUSIC_direct} compared to  LKH.
These improvements are attained relative to the version of LKH which is supposed to work best for large instances (i.e., using POPMUSIC \cite{LKH_popmusic} candidate sets). 
Indeed, for other versions (ALPHA\footnote{See \ref{C1:candidateset} in \Cref{subsection:LKH} for the role of candidate sets and \cite{LKH_original} and \cite{LKH_popmusic} for an explanation of ALPHA and POPMUSIC candidate sets, respectively.}) our improvements over LKH are even higher. 
The differences of our method and LKH for small and medium instances are small or even negligible.

In total, the computational results advocate the use of the relaxed positive gain criterion in practice especially for large instances and to study further variants of LKH beyond the strict, original positive gain criterion.

\subsection{Related work}
We propose an improvement for the Lin-Kernighan-Helsgaun heuristic, which is a tour improvement heuristic. These begin with a given initial tour and improve it by deleting edges of the tour and adding new ones.

\textbf{Tour improvement heuristics.} 
The \struc{2-opt} heuristic by Croes \cite{Croes2opt} and \struc{3-opt} heuristic \cite{Lin3opt} by Lin delete two and three edges respectively of a given tour to establish a better tour by adding new edges.
These heuristics set the groundwork for the tour improvement heuristic \struc{Lin-Kernighan (LK)} by Lin and Kernighan \cite{LK_original} that was already developed in 1973 and gave rise to today's state-of-the-art implementations. 
The main idea of LK is that in each step, instead of a fixed number of edges, a variable number $k$ of edges is deleted from the given tour, and $k$ new edges are added.
See \cite{rego2011traveling} for an overview of LK implementations, including the most prominent ones by Johnson and McGeoch \cite{johnson1997traveling}, Applegate, Bixby, Chv\'atal and Cook \cite{applegate1999finding}, Neto \cite{neto1999efficient}, Helsgaun \cite{LKH_original}, and Nguyen, Yoshihara, Yamamori and Yasunaga \cite{LK-NYYY}.
The Chained-LK \cite{martin1992large} by Martin, Otto and Felten iteratively restarts the LK procedure with a purtubated previous best tour. 
Its implementation by Applegate, Cook and Rohe \cite{applegate2003chained} is used as a start heuristic in the CONCORDE code \cite{concorde}.
The implementation by Neto \cite{neto1999efficient} specifically considers clustered TSP instances and suggests a modified gain criterion for these cases.

\textbf{Benchmark instances.} Helsgaun's implementation \cite{LKH_original, helsgaun2009general, LKH_popmusic} is a state-of-the-art heuristic having set many records in the benchmarking of TSP instances, such as the instances of the TSPLIB \cite{TSPLIB}, the VLSI dataset \cite{vlsi}, the national TSP instances \cite{nationalTSP}, the instances of the eighth DIMACS challenge \cite{DIMACS}, and the $T_{nm}$ instances \cite{Hougardy_Tnm}.
The current version of the code \struc{LKH-3} \cite{LKH3} has multiple extensions for constrained TSP and vehicle routing problems \cite{LKH_extension}.

\textbf{Modifications of LKH.}
J\"ager, Dong, Goldengorin, Molitor and Richter \cite{jager2014backbone}, building on \cite{richter2007improving}, improve LKH via introducing backbone search by \cite{zhang2005novel} and double bridge moves reaching the best results on some VLSI \cite{vlsi} instances.
For clustered instances, Hains, Whitley and Howe \cite{HWH2012_improving} improve LKH using generalized partition crossover \cite{whitley2009tunneling}.

In the recent years, there is active research to support and upgrade TSP heuristics via machine learning, see for example the survey \cite{bengio2021machine} as an overview.
The heuristics VSR-LKH by Zheng, He, Zhou, Jin and Li \cite{zheng2021combining}, \cite{ZHENG2023} and NeuroLKH by Xin, Song, Cao and Zhang \cite{xin2021neurolkh} combine learning methods with LKH in the process of generating the candidate sets and selecting edges.
They report quality improvements over the original LKH on instances up to 85900 vertices \cite{ZHENG2023} and 5000 vertices \cite{xin2021neurolkh}.
The recent learning approach by Wang, Zhang and Tang \cite{wang2023discovering} enhances the two previous methods.

To the best of our knowledge exists no systematic study of relaxing the positive gain criterion. 

\textbf{Non-LK based heuristics.} This work focuses on LK type heuristics. 
For an older, general survey on TSP heuristics see \cite{johnson1997traveling}. 
Besides LK-based methods the recent research considers nature inspired heuristics for example ant colony optimization \cite{skinderowicz2022improving} based on \cite{dorigo1996ant} or approaches based on swarm intelligence \cite{eberhart2001swarm}.
The nature inspired algorithm by Xie and Liu \cite{xie2008multiagent} suggests to apply the multi-agent optimization system (MAOS) to the TSP. 
They report that their results can compete with LKH on the VLSI instances.

Nagata and Kobayashi \cite{nagata2013powerful} proposed a genetic algorithm that uses edge assembly crossover (EAX) based on \cite{nagata1997edge}.
They report results that can outperform LKH on instances up to 200.000 vertices \cite{kotthoff2015improving}.
A recent preprint \cite{fu2023hierarchical} proposes a destroy-and-repair approach for large TSP instances that can compete with LKH having set a new record for two instances with three and ten million vertices.

\section*{Acknowledgments}
All authors were supported by the German Federal Ministry for Economic Affairs and Climate Action (BMWK), project ProvideQ.

\section{Preliminaries}
\label{section:preliminaries}
\subsection{Notation and basic definitions}
\label{subsection:notation}
\begin{figure}[t] %Running Example
	\begin{center}
		\begin{subfigure}[b]{0.29\textwidth}
			\centering
			\resizebox{\linewidth}{!}{
				\begin{tikzpicture}
					%nodes
					\node (a) at (0,0) {};
					\node (b) at (4,0) {};
					\node (c) at (6, -3.46) {};
					\node (d) at (4,-6.93) {};
					\node (e) at (0,-6.93) {};
					\node (f) at (-2,-3.46) {};
					
					%Auxiliary nodes for edge labels
					%length 6
					\node (h1)  [thick, font=\fontsize{20}{0}\selectfont] at (1.2, -4.2) {6};
					\node (h2)  [thick, font=\fontsize{20}{0}\selectfont] at (2.4, -4.7) {6};
					\node (h2)  [thick, font=\fontsize{20}{0}\selectfont] at (3, -3.15) {6};
					%length 5
					\node (h3)  [thick, font=\fontsize{20}{0}\selectfont] at (1.0, -0.9) {5};
					\node (h4)  [thick, font=\fontsize{20}{0}\selectfont] at (-0.8, -2.4) {5};
					\node (h5)  [thick, font=\fontsize{20}{0}\selectfont] at (-0.2, -5.3) {5};
					\node (h6)  [thick, font=\fontsize{20}{0}\selectfont] at (2.5, -6.5) {5};
					\node (h7)  [thick, font=\fontsize{20}{0}\selectfont] at (4.6, -4.7) {5};
					\node (h9)  [thick, font=\fontsize{20}{0}\selectfont] at (3.8, -1.4) {5};
	
					%length 4 edges and nodes
					\draw[thick, lightgray, font=\fontsize{20}{0}\selectfont]  (a) edge node[above, black]{4} (b);
					\draw[thick, lightgray, font=\fontsize{20}{0}\selectfont]  (d) edge node[below, black]{4} (e);
					%length 3 edges and nodes
					\draw[thick, lightgray, font=\fontsize{20}{0}\selectfont]  (b) edge node[above right, black]{3} (c);
					\draw[thick, lightgray, thick, font=\fontsize{20}{0}\selectfont]  (c) edge node[below right, black]{3} (d);
					\draw[thick, lightgray, font=\fontsize{20}{0}\selectfont]  (e) edge node[below left, black]{3} (f);
					\draw[thick, lightgray, font=\fontsize{20}{0}\selectfont, thick]  (f) edge node[above left, black]{3} (a);
					
					%length 6 edges
					\draw[thick, lightgray]  (a) edge (d);
					\draw[thick, lightgray]  (b) edge (e);
					\draw[thick, lightgray]  (c) edge (f);
					
					%length 5 edges
					\draw[thick, lightgray]  (a) edge (c);
					\draw[thick, lightgray]  (b) edge (d);
					\draw[thick, lightgray]  (c) edge (e);
					\draw[thick, lightgray]  (d) edge (f);
					\draw[thick, lightgray]  (e) edge (a);
					\draw[thick, lightgray]  (f) edge (b);
					
					%Nodes above edges
					\draw[fill] (0,0) circle (2pt);
					\draw[fill] (4,0) circle (2pt);
					\draw[fill] (6, -3.46) circle (2pt);
					\draw[fill] (4,-6.93) circle (2pt);
					\draw[fill] (0,-6.93) circle (2pt);
					\draw[fill] (-2,-3.46) circle (2pt);
					
			\end{tikzpicture}}
			\caption{}
			\label{LKedgeweights:figure}
		\end{subfigure}
		\hspace{1.5em}
		\begin{subfigure}[b]{0.29\textwidth} %Starttour
			\centering
			\resizebox{\linewidth}{!}{
				\begin{tikzpicture}
					\node (a) at (0,0) {};
					\node (b) at (4,0) {};
					\node (c) at (6, -3.46) {};
					\node (d) at (4,-6.93) {};
					\node (e) at (0,-6.93) {};
					\node (f) at (-2,-3.46) {};
					\node (d2) at (4,-7.6) {}; %to shift image
					
					%length 6
					\draw[thick]  (a) edge (d);
					\draw[thick]  (b) edge (e);
					\draw[thick, lightgray]  (c) edge (f);
					
					%length 5
					\draw[thick, lightgray]  (a) edge (c);
					\draw[thick, lightgray]  (b) edge (d);
					\draw[thick, lightgray]  (c) edge (e);
					\draw[thick, lightgray]  (d) edge (f);
					\draw[thick, lightgray]  (e) edge (a);
					\draw[thick, lightgray]  (f) edge (b);
										
					%length 4
					\draw[thick, lightgray]  (a) edge (b);
					\draw[thick, lightgray]  (d) edge  (e);
					
					%length 3
					\draw[thick]  (b) edge  (c);
					\draw[thick]  (c) edge  (d);
					\draw[thick]  (e) edge (f);
					\draw[thick]  (f) edge (a);
				
					%Nodes above edges
					\draw[fill] (0,0) circle (2pt);
					\draw[fill] (4,0) circle (2pt);
					\draw[fill] (6, -3.46) circle (2pt);
					\draw[fill] (4,-6.93) circle (2pt);
					\draw[fill] (0,-6.93) circle (2pt);
					\draw[fill] (-2,-3.46) circle (2pt);
				\end{tikzpicture}
			}
			\caption{} 
			\label{starttour}
		\end{subfigure}
		\hspace{1.5em}
		\begin{subfigure}[b]{0.29\textwidth}
			\centering
			\resizebox{\linewidth}{!}{
				\begin{tikzpicture}	
					\node (a) at (0,0) {};
					\node (b) at (4,0) {};
					\node (c) at (6, -3.46) {};
					\node (d) at (4,-6.93) {};
					\node (e) at (0,-6.93) {};
					\node (f) at (-2,-3.46) {};
					\node (s1) at (4,-7.6) {}; %to shift image
					
					%length 6
					\draw[dashed, red, line width=0.05cm, dash pattern=on 8pt off 5pt]  (a) edge (d);
					\draw[dashed, red, line width=0.05cm, dash pattern=on 8pt off 5pt]  (b) edge (e);
					\draw[thick, lightgray]  (c) edge (f);
					
					%length 5
					\draw[thick, lightgray]  (a) edge (c);
					\draw[thick, lightgray]  (b) edge (d);
					\draw[thick, lightgray]  (c) edge (e);
					\draw[thick, lightgray]  (d) edge (f);
					\draw[thick, lightgray]  (e) edge (a);
					\draw[thick, lightgray]  (f) edge (b);
					
					%length 4
					\draw[thick, blue]  (a) edge (b);
					\draw[thick, blue]  (d) edge  (e);
					
					%length 3
					\draw[thick]  (b) edge  (c);
					\draw[thick]  (c) edge  (d);
					\draw[thick]  (e) edge (f);
					\draw[thick]  (f) edge (a);
					
					%Nodes above edges
					\draw[fill] (0,0) circle (2pt);
					\draw[fill] (4,0) circle (2pt);
					\draw[fill] (6, -3.46) circle (2pt);
					\draw[fill] (4,-6.93) circle (2pt);
					\draw[fill] (0,-6.93) circle (2pt);
					\draw[fill] (-2,-3.46) circle (2pt);
			\end{tikzpicture}}
			\caption{}
			\label{2opt}
		\end{subfigure}
		\caption{Running example adjusted from \cite{thesis_BO}:
			\cref{LKedgeweights:figure} depicts edge costs.
			A possible starting tour is in \cref{starttour}.
			Example of a $2$-opt-move in \cref{2opt}:  
			Edges of starting tour are black, we delete dashed red edges and insert blue edges.} 
		\label{lkhimprovedexample}
	\end{center}
\end{figure}
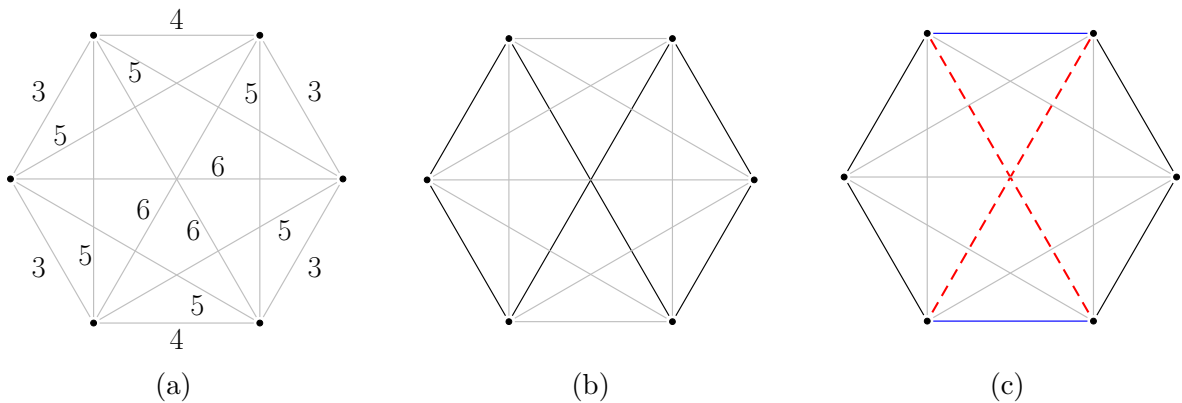

Throughout this article we use the following \linebreak notation. 
If not stated otherwise, we consider a \struc{complete}, \struc{weighted} and \struc{undirected} \linebreak \struc{graph} $\struc{G = (V,E)}$, consisting of $n$ \struc{vertices} $\struc{V} := \{1,\hdots,n\}$ and \struc{edges} \linebreak ${\struc{E}:=\{(i,j) : (i,j) \in\ V\times V, i\neq j\}}$. 
The edge weights are called \struc{costs} and assigned by a function $\struc{c}: E \to \mathbb{R}_{> 0}$.
A \struc{tour} $\struc{T}$ is a \struc{Hamiltonian cycle} of $G$, i.e. a cycle that includes every vertex exactly once.
The \struc{Traveling Salesman Problem (TSP)} asks for a tour $T$ in which the sum of costs $\struc{c(T)}$ over the edges traversed is minimal. 
We consider the \struc{metric} and \struc{symmetric} TSP, i.e. the edge costs fulfill the triangle equality, and for all edges $(i,j)$ it holds that $c(i,j) = c(j,i)$, respectively.
See \cref{LKedgeweights:figure} for the graph that we use as a running example.
The unique optimal tour $T^*$ includes only the outer edges and has total cost $c(T^*)=20$.

\subsection{Lin-Kernighan-Helsgaun (LKH)}\label{subsection:LKH}

LKH is a heuristic that iteratively tries to improve a given tour $T$.
That means, given $T$ as a starting tour, we try to establish a new, improved tour $T'$ with $c(T') < c(T)$.
We pursue this via deleting edges from $T$ and adding new ones.
When constructing the sets of edges to be deleted respectively added we obey certain heuristic limiting criteria that we describe below.
This approach has significantly less running time than an exhaustive search over all sets of edges that can be exchanged.
However, since those criteria are heuristic there is no guarantee that the optimal tour is found.

We briefly recall the concept of LKH and those details that are relevant to our improvement.
We follow the notation as given in the original publication of Helsgaun \cite{LKH_original}:
A \struc{$k$-opt-move} is the action of deleting $k$ edges from the current tour $T$ and adding $k$ edges such that we obtain a new feasible tour $T'$ with $c(T') < c(T)$. 
For an example of a $2$-opt-move derived from a starting tour, see \cref{starttour} and \cref{2opt}.
During a $k$-opt-move, the set \struc{$X=\{x_1, \hdots, x_k\}$} contains the edges that we delete from the tour while we add the set of edges \struc{$Y=\{y_1, \hdots, y_k\}$}.
In what follows we restrict $X$ and $Y$ to sets that form an \struc{alternating circle} \struc{$P = [x_1, y_1,\dots, x_k,y_k]$}.
We label the vertices that are incident with the edges in $X$ and $Y$ ascendingly as \struc{$t_j$} for $j=1,\dots,2k$. 
In particular, starting with vertex $t_1$, we delete the edge $x_1 = (t_1, t_2)$. 
Then we add edge $y_1 = (t_2, t_3)$.
Accordingly, we can think of a $k$-opt-move as a \struc{sequential edge exchange} replacing $x_i = (t_{2i-1}, t_{2i})$ by $y_i = (t_{2i}, t_{2i+1})$ for all $i=1,\dots, k$.
Therefore, we can also write the alternating path $P$ with respect to the vertices as $P=[t_1, t_2,\dots, t_{2k}, t_{2k+1}]$.
As $T'$ has to be a tour, we have $t_{2k+1} = t_1$.
Consider the following example that applies the notation and points out the basic procedure and idea of LKH:	
\begin{figure}[t] %Example of LKH without adjustment
	\begin{center}
		\begin{subfigure}[b]{0.29\textwidth}
			\centering
			\resizebox{\linewidth}{!}{
				\begin{tikzpicture}
					\node (a) at (0,0) {};
					\node (b) at (4,0) {};
					\node (c) at (6, -3.46) {};
					\node (d) at (4,-6.93) {};
					\node (e) at (0,-6.93) {};
					\node (f) at (-2,-3.46) {};
					\node (s1) at (4,-7.5) {}; %to shift image
					\node (s2) at (6.5,-3.46) {}; %to shift image
					
					%length 6
					\draw[dotted, red, font=\fontsize{20}{0}\selectfont, line width=0.05cm]  (a) edge (d);
					\node (h2)  [thick, red, font=\fontsize{20}{0}\selectfont] at (2.05, -4.85) {$x_1$?};
					\draw[thick]  (b) edge (e);
					\draw[thick, lightgray]  (c) edge (f);
					
					%length 5
					\draw[thick, lightgray]  (a) edge (c);
					\draw[thick, lightgray]  (b) edge (d);
					\draw[thick, lightgray]  (c) edge (e);
					\draw[thick, lightgray]  (d) edge (f);
					\draw[thick, lightgray]  (e) edge (a);
					\draw[thick, lightgray]  (f) edge (b);
					
					%length 4
					\draw[thick, lightgray]  (a) edge (b);
					\draw[thick, lightgray]  (d) edge (e);
					
					%length 3
					\draw[thick]  (b) edge  (c);
					\draw[thick]  (c) edge  (d);
					\draw[thick]  (e) edge (f);
					\draw[dotted, red, font=\fontsize{20}{0}\selectfont, line width=0.05cm]  (f)  edge node[above left, red]{$x_1?$} (a);
					
					%Nodes above edges
					\draw[fill] (0,0) circle (2pt) node[left, thick, font=\fontsize{20}{0}\selectfont] {$t_1$};
					\draw[fill] (4,0) circle (2pt);
					\draw[fill] (6, -3.46) circle (2pt);
					\draw[fill] (4,-6.93) circle (2pt) node[right, thick, font=\fontsize{20}{0}\selectfont] {$t_2$?};
					\draw[fill] (0,-6.93) circle (2pt);
					\draw[fill] (-2,-3.46) circle (2pt) node[left, thick, font=\fontsize{20}{0}\selectfont] {$t_2$?};	
				\end{tikzpicture}
			}
			\caption{Selection of $x_1$ to delete.} 
			\label{lkh2}
		\end{subfigure}
		\hspace{1.em}
		\begin{subfigure}[b]{0.29\textwidth}
			\centering
			\resizebox{\linewidth}{!}{
				\begin{tikzpicture}
					\node (a) at (0,0) {};
					\node (b) at (4,0) {};
					\node (c) at (6, -3.46) {};
					\node (d) at (4,-6.93) {};
					\node (e) at (0,-6.93) {};
					\node (f) at (-2,-3.46) {};
					\node (a2) at (0.1,0) {};
					\node (d2) at (4,-6.757) {};
					
					%length 6
					\draw[dashed, red, font=\fontsize{20}{0}\selectfont, line width=0.05cm, dash pattern=on 8pt off 5pt]  (a) edge (d);
					\node (h2)  [thick, red, font=\fontsize{20}{0}\selectfont] at (2.2, -4.85) {$x_1$};
					\draw[thick]  (b) edge (e);
					\draw[thick, lightgray]  (c) edge (f);
					\draw[dotted, blue, line width=0.05cm]  (a2) edge (d2);
					\node (h10)  [thick, blue, font=\fontsize{20}{0}\selectfont] at (3.2, -4) {$y_1$?};
					
					%length 5
					\draw[thick, lightgray]  (a) edge (c);
					\draw[dotted, blue, line width=0.05cm]  (b) edge (d);
					\draw[thick, lightgray]  (c) edge (e);
					\draw[dotted, blue, line width=0.05cm]  (d) edge (f);
					\draw[thick, lightgray]  (e) edge (a);
					\draw[thick, lightgray]  (f) edge (b);
					\node (h5)  [thick, blue, font=\fontsize{20}{0}\selectfont] at (0.45, -4.2) {$y_1$?};
					\node (h9)  [thick, blue, font=\fontsize{20}{0}\selectfont] at (3.5, -2.1) {$y_1$?};
					
					%length 4
					\draw[thick, lightgray]  (a) edge (b);
					\draw[dotted, blue, line width=0.05cm]  (d) edge (e);
					\node (y2) [thick, blue, font=\fontsize{20}{0}] at (1.7, -6.65) {$y_1$?};
					
					%length 3
					\draw[thick]  (b) edge  (c);
					\draw[thick]  (c) edge  (d);
					\draw[thick]  (e) edge (f);
					\draw[thick]  (f)  edge (a);
					
					%Nodes above edges
					\draw[fill] (0,0) circle (2pt) node[left, thick, black, font=\fontsize{20}{0}\selectfont, align=left] {$t_1$\\$t_3$?};
					\draw[fill] (4,0) circle (2pt) node[right, thick, black, font=\fontsize{20}{0}\selectfont] {$t_3?$};
					\draw[fill] (6, -3.46) circle (2pt) {};
					\draw[fill] (4,-6.93) circle (2pt)node[right, thick, black, font=\fontsize{20}{0}\selectfont] {$t_2$};
					\draw[fill] (0,-6.93) circle (2pt)  node[left, black, font=\fontsize{20}{0}\selectfont] {$t_3$?};
					\draw[fill] (-2,-3.46) circle (2pt) node[left, thick, font=\fontsize{20}{0}\selectfont, align=left] {$t_3$?};
						
				\end{tikzpicture}
			}
			\caption{Selection of $y_1$ to add.} 
			\label{lkh3}
		\end{subfigure}
		\hspace{1em}
		\begin{subfigure}[b]{0.29\textwidth}
			\centering
			\resizebox{\linewidth}{!}{
				\begin{tikzpicture}
					\node (a) at (0,0) {};
					\node (b) at (4,0) {};
					\node (c) at (6, -3.46) {};
					\node (d) at (4,-6.93) {};
					\node (e) at (0,-6.93) {};
					\node (f) at (-2,-3.46) {};
					\node (a2) at (0.1,0) {};
					\node (d2) at (4,-6.757) {};
					\node (s2) at (7,-3.46) {}; %to shift image
					
					%length 6
					\draw[dashed, red, font=\fontsize{20}{0}\selectfont, line width=0.05cm, dash pattern=on 8pt off 5pt]  (a) edge (d);
					\node (h2)  [thick, red, font=\fontsize{20}{0}\selectfont] at (2.2, -4.85) {$x_1$};
					\draw[dashed, red, line width=0.05cm, dash pattern=on 8pt off 5pt]  (b) edge (e);
					\node (h1)  [thick, red, font=\fontsize{20}{0}\selectfont] at (1.05, -4.2) {$x_2$};
					\draw[thick, lightgray]  (c) edge (f);

					%length 5
					\draw[thick, lightgray]  (a) edge (c);
					\draw[thick, lightgray]  (b) edge (d);
					\draw[thick, lightgray]  (c) edge (e);
					\draw[thick, lightgray]  (d) edge (f);
					\draw[thick, lightgray]  (e) edge (a);
					\draw[thick, lightgray]  (f) edge (b);
					
					%length 4
					\draw[thick, blue, font=\fontsize{20}{0}\selectfont]  (a) edge node[above, blue]{$y_2$} (b);
					\draw[thick, blue, line width=0.05cm]  (d) edge (e);
					\node (y2) [thick, blue, font=\fontsize{20}{0}] at (1.7, -6.65) {$y_1$};
					
					%length 3
					\draw[thick]  (b) edge  (c);
					\draw[thick]  (c) edge  (d);
					\draw[thick]  (e) edge (f);
					\draw[thick]  (f)  edge (a);
					
					%Nodes above edges
					\draw[fill] (0,0) circle (2pt) node[above, thick, black, font=\fontsize{20}{0}\selectfont, align=left] {$t_1$};
					\draw[fill] (4,0) circle (2pt) node[right, thick, black, font=\fontsize{20}{0}\selectfont] {$t_4$};
					\draw[fill] (6, -3.46) circle (2pt) {};
					\draw[fill] (4,-6.93) circle (2pt)node[right, thick, black, font=\fontsize{20}{0}\selectfont] {$t_2$};
					\draw[fill] (0,-6.93) circle (2pt)  node[left, black, font=\fontsize{20}{0}\selectfont] {$t_3$};
					\draw[fill] (-2,-3.46) circle (2pt) node[left, thick, font=\fontsize{20}{0}\selectfont, align=left] {};
						
				\end{tikzpicture}
			}
			\caption{Feasible edge exchanges.} 
			\label{lkh4}
		\end{subfigure}
		
		\caption{Display of \cref{exampleLKHnormal}. Black edges belong to starting tour. Red edges indicate removal. Blue edges indicate insertion. In \cref{lkh2} we select one of the red dotted edges to be removed, in \cref{lkh3} we choose one of the blue dotted edges to be added.} 
		\label{lkhexample}
	\end{center}
\end{figure}
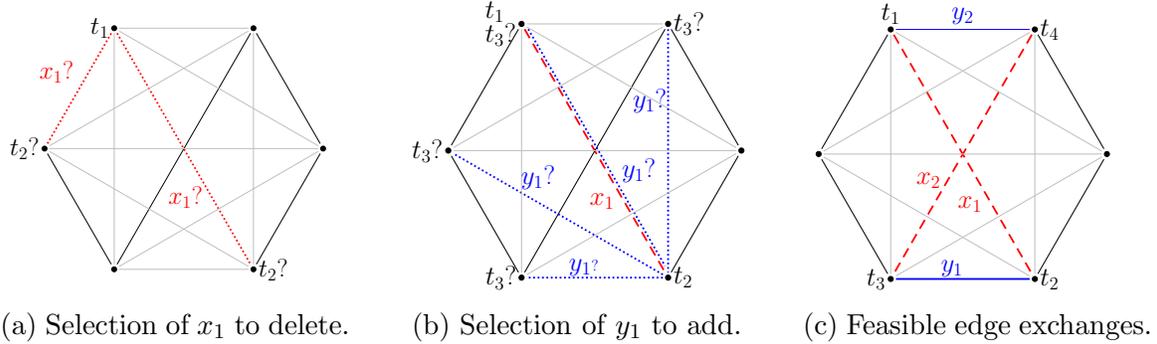

\begin{examplex}
	\label{exampleLKHnormal}
	The tour $T$ in \cref{starttour} has total cost $c(T) = 4 \cdot 3 + 2 \cdot 6 = 24$. 
	Consider \cref{lkhexample}: 
	We assume $T$ to be the starting tour. 
	We randomly choose $t_1$ as the starting vertex of the alternating path. 
	The red dotted edges in \cref{lkh2} are the only two edges of the tour $T$ that are incident with $t_1$.
	Thus, we need to select one of them to be the first edge $x_1 = (t_1, t_2)$ to be deleted.
	The idea for choosing a pair of $(x_i, y_i)$ is to delete an edge with preferably high cost and then add an edge with preferably lower cost.
	In this example, we select $\color{red}x_1$ \color{black} as indicated by the red dashed edge in \cref{lkh3}.
	
	We now have to choose an edge that is incident with $t_2$ to be added as $y_1 = (t_2, t_3)$.
	In this example, there are four such edges indicated by the blue dotted lines in \cref{lkh3}.
	Following the goal of adding an edge with low cost, we choose to insert the blue solid edge \color{blue}$y_1$ \color{black} in \cref{lkh4}. 
	This exchange leads to a cost difference of $c(\color{red}x_1\color{black}) -  c(\color{blue}y_1\color{black})  = \color{red}6\color{black} - \color{blue}4 \color{black}= 2 > 0$.
	Continuing the sequential edge exchange, we decrease the total cost and arrive at the tour $T'$ in \cref{lkh4} with $c(T') = 4\cdot 3+ 2\cdot 6 - 2\cdot\color{red}{ 6}\color{black} +2\cdot \color{blue}{4}\color{black}=20 < c(T)$. 
	In this case, $T'$ is the only optimal tour.
\end{examplex}
In the previous example we demonstrated the intuition for selecting edges for exchanges. 
Now we describe the heuristic criteria \ref{C1:candidateset}, \ref{C2:gain}, \ref{C3:feasibility}, \ref{C4:sequential} and \ref{C5:disjunctivity} used in LKH\footnote{See \cite{LKH_original} and \cite{LKH_koptTR} for a more detailed description of the criteria used in LKH, and its extension LKH-3 \cite{LKH_extension} that modify the rules originally developed in LK \cite{LK_original}.}.
Choosing such a set of heuristic criteria, one has to balance the computational effort for finding the sets to be exchanged and the improvement achieved by this exchange.
 
Suppose that we have removed the edge $x_i = (t_{2i-1}, t_{2i})$ and that we need to add an edge $y_i$ instead. 
The following two criteria concern the selection of $y_i$:
\begin{enumerate}[label=(C\arabic*)] 
	\item To limit the amount of potential edges $y_i$, for each vertex $t_{2i}$ of the graph, we determine a \struc{candidate set} of certain edges incident to $t_{2i}$, and an ranking of these edges during a preprocessing phase.
	We select the edge $y_i=(t_{2i}, t_{2i+1})$ only among the candidate set of the vertex $t_{2i}$.
	We choose the first edge in the given ranking that fulfills the subsequent criteria.
	In LKH, the construction of the candidate sets and the ranking of its included edges is based on \struc{$\alpha$-values} \cite{LKH_original} or the meta-heuristic \struc{POPMUSIC} \cite{LKH_popmusic}. \label{C1:candidateset}
	\item 
	We define
	\begin{align*}
		\struc{G_i} := \sum_{\ell = 1}^{i} \underbrace{c(x_{\ell}) - c(y_{\ell})}_{\struc{g_{\ell}}}
	\end{align*} 
	as the \struc{total gain} in costs that is achieved after deleting and adding $i$ many edges during a $k$-opt move.
	By \struc{$g_{\ell}$} we denote the cost difference of the exchanged edges $x_{\ell}$ and $y_{\ell}$ in step $\ell$.
	The \struc{positive gain criterion} only allows to add a potential edge $y_i$ after deleting $x_i$ if the total gain $G_i$ is positive.
	
	This heuristically pursues the goal of only performing exchanges that decrease the total cost.
	Exchanging edge $x_i$ by $y_i$ with non-positive cost difference \linebreak${g_{i}= c(x_i)-c(y_i) \leq 0}$ is possible, but only if the total gain $G_i$ after the first $i$ exchanges is positive. \label{C2:gain}
\end{enumerate}

We require the following criterion to delete an edge $x_i$:
\begin{enumerate}[label=(C\arabic*)]
	\addtocounter{enumi}{2}
	\item
	For a fixed parameter $r$  
	the \struc{feasibility criterion} demands for every multiple of $r$, i.e. $ir$ with $i \in \mathbb{N}\setminus\{0\}$, that we can \struc{close up} the tour if we delete $x_{ir}=(t_{2ir-1}, t_{2ir})$.
	Closing up the tour means that deleting the edges $\{x_1,\dots,x_{ir}\}$ from the tour and adding the edges $\{y_1,\dots, y_{ir-1}, (t_{2ir}, t_1)\}$ instead would result in a tour. 
	
	Note that the criterion does not require the tour resulting from the alternating circle $P$ to have lower cost than the original tour.
	From the extension LKH-2 \cite{LKH_koptTR} on, the user can choose $r$ as a parameter.
	By default, $r$ is set to 5.
	 \label{C3:feasibility}
\end{enumerate}

The following two criteria restrict the selection of $x_i$ as well as the selection of $y_i$:
\begin{enumerate}[label=(C\arabic*)]
	\addtocounter{enumi}{3}
	\item
	The \struc{sequential exchange criterion} demands that the removed and added edges form an alternating path.
	
	If LKH does not find an improvement of the tour via sequential edge exchanges, LKH additionally tries exchanging edges that do not form an alternating path.
	We do not further describe these \struc{non-sequential edge exchanges}, see \cite{helsgaun2009general}. \label{C4:sequential}
	\item Using the same fixed parameter $r$ as in \ref{C3:feasibility} the \struc{disjunctivity criterion} in LKH demands that every edge of the form $x_{ir}$ is not an element of $\{y_1, y_2,\dots, y_{ir-1}\}$.
	Except for the edges of the form $x_{ir}$ the sets $X$ and $Y$ may overlap.
	The name of the criterion originates from LK \cite{LK_original} where $X$ and $Y$ still have to be disjoint. \label{C5:disjunctivity}
\end{enumerate}
The algorithm for finding an alternating path starts by choosing a random vertex $t_1$.
Then it builds an alternating path according to the criteria \ref{C1:candidateset}, \ref{C2:gain}, \ref{C3:feasibility}, \ref{C4:sequential} and \ref{C5:disjunctivity}.
As termination criterion we check for each $i\geq 2$ after deleting $x_i$ whether we can close up the tour with positive gain.
The final edge $(t_{2i}, t_1)$, which closes the alternating path into an alternating circle, is not required to be contained in a candidate set.

Now we pay further attention to the use of the positive gain criterion \ref{C2:gain}. 
We provide a pseudocode in \cref{alg:LKH_Gain_Criterion}  representing the selection of the edge $y_i = (t_{2i},t_{2i+1})$ to be added.
We assume that we have already exchanged $i-1$ edges and selected the next edge  $x_i = (t_{2i-1}, t_{2i})$ to be deleted.
Thus, we input an alternating path $P = [t_1,..., t_{2i}]$ that ends with an edge $x_i$ to the function \texttt{FindSequentialEdgeToBeAdded}.
To find $y_i = (t_{2i}, t_{2i+1})$ this function searches for the first edge in the candidate set of $t_{2i}$ that satisfies the positive gain criterion \ref{C2:gain}\footnote{Note that the shown pseudocode is a simplification. 
		In the original LKH code there is a parameter \enquote{GainCriterionUsed}, which is \enquote{true} or \enquote{false} depending on whether the user wants to include \ref{C2:gain} as a criterion
		or not. 
		In this paper we always set GainCriterionUsed to the default value \enquote{true}. 
		Therefore, we did not insert it in the pseudocode.
		In the equivalent of \cref{line:alg_if} the original code checks whether \enquote{GainCriterionUsed = false \textbf{or} $G_i > 0$} instead of \enquote{$G_i > 0$}.},
as well as \ref{C4:sequential} and \ref{C5:disjunctivity}.
For simplification we denote the latter two criteria via the boolean function \texttt{OtherCriteriaSatisfied}.
If we find a candidate satisfying all criteria we return \enquote{Found} and append the corresponding edge to the alternating path. 
If no element of the candidate set satisfies the criteria, then the algorithm returns \enquote{Not found}.
\begin{figure}
	\begin{algorithm}[\texttt{FindSequentialEdgeToBeAdded}]~
	\label{alg:LKH_Gain_Criterion}
	
		\noindent \INPUT{Alternating path $P = [t_1, ..., t_{2i}]$, gain $G_{i-1}$ of $[t_1,...t_{2i-1}]$}\\
		\OUTPUT{\enquote{Found} / \enquote{Not found}, alternating path $P = [t_1, ..., t_{2i}, t_{2i+1}]$}
		\begin{algorithmic}[1] 
			\For{$(t_{2i},t_{2i+1})$ in \texttt{CandidateSet}($t_{2i}$)} 
			\State $G_i \gets G_{i-1} + c(t_{2i-1}, t_{2i}) - c(t_{2i}, t_{2i+1})$
			\If{ $G_i > 0$  \textbf{and}  $\texttt{OtherCriteriaSatisfied}(P, (t_{2i}, t_{2i+1}))$} \label{line:alg_if}
			\State \Return \enquote{Found}, $P = [t_1, ..., t_{2i}, t_{2i+1}]$
			\EndIf	
			\EndFor
			\State \Return \enquote{Not found}
		\end{algorithmic}	
		
	\end{algorithm}
\caption*{Simplified example of the usage of the positive gain criterion \ref{C2:gain} in the original LKH code \cite{LKH3}.
Based on this example code, we provide our variation in \cref{alg:LKH_0.1_Gain_Criterion}.}
\end{figure}

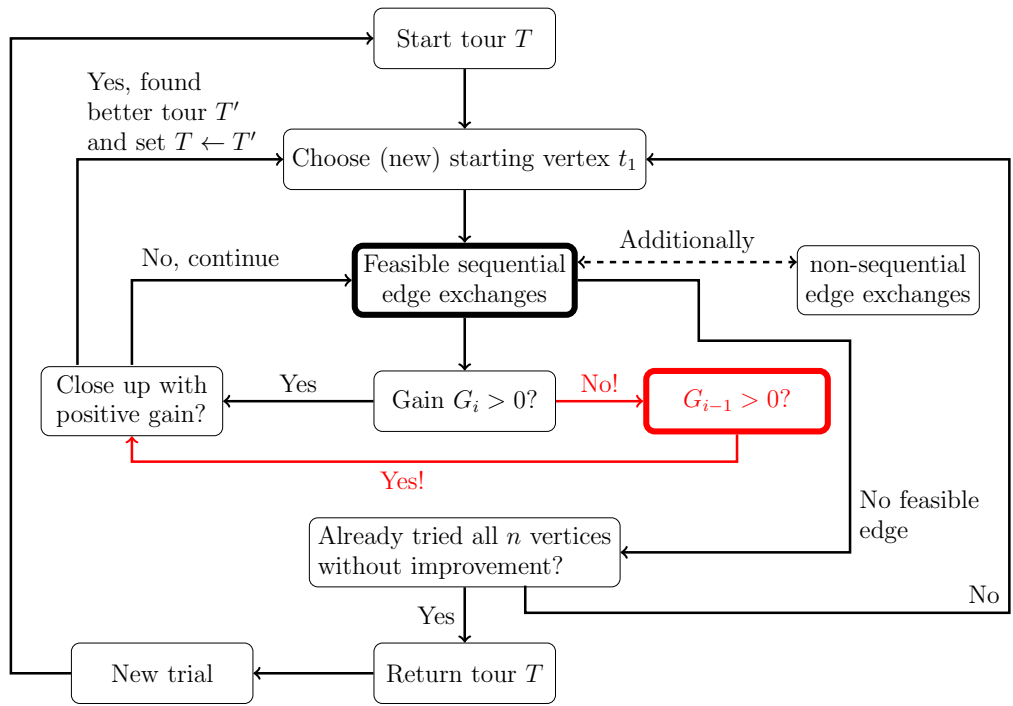
\begin{figure}[t] %Flowchart of LKH with adjusted gain criterion
	\begin{center}
		\centering
		\resizebox{0.85\linewidth}{!}{
			\begin{tikzpicture}[node distance=2cm]
				%Nodes
				\node (starttour) [rectangle, rounded corners, minimum width=3cm, minimum height=1cm,text centered, draw=black, fill=white] {Start tour $T$};
				\node (startvertex) [rectangle, rounded corners, minimum width=3cm, minimum height=1cm,text centered, draw=black, fill=white, below of=starttour] {Choose (new) starting vertex $t_1$};
				
				\node (exchanges) [rectangle, rounded corners, minimum width=3cm, minimum height=1cm,text centered, draw=black, fill=white, below of=startvertex, line width=1mm, align=center] {Feasible sequential \\edge exchanges};
				
				\node (nonseq) [rectangle, rounded corners, minimum width=3cm, minimum height=1cm,text centered, draw=black, fill=white, xshift=5cm, right of=exchanges, align = center] {non-sequential\\edge exchanges};
				
				\node (criteria) [rectangle, rounded corners, minimum width=3cm, minimum height=1cm,text centered, draw=black, fill=white, below of=exchanges] {Gain $G_i>0$?};
				
				\node (improvement) [rectangle, rounded corners, minimum width=3cm, minimum height=1cm,text centered, draw=black, fill=white, left of=criteria, xshift=-3.5cm, align=left] {Close up with \\positive gain?};
				
				\node (allvertices) [rectangle, rounded corners, minimum width=3cm, minimum height=1cm,text centered, draw=black, fill=white, yshift = -0.5cm, below of=criteria, align = left] {Already tried all $n$ vertices\\without improvement?};
				
				\node (finish) [rectangle, rounded corners, minimum width=3cm, minimum height=1cm,text centered, draw=black, fill=white, below of=allvertices] {Return tour $T$};
				
				\node (newrun) [rectangle, rounded corners, minimum width=3cm, minimum height=1cm,text centered, draw=black, fill=white, left of=finish, align=left, xshift = -3cm] {New trial};
				
				\node (negGain) [line width= 1mm, red, rectangle, rounded corners, minimum width=3cm, minimum height=1cm,text centered, draw=red, fill=white, right of=criteria, xshift=2.5cm, align=left] {$G_{i-1} > 0$?};
				
				%Arrows
				\draw [->, very thick, very thick] (starttour) -- (startvertex);
				
				\draw [->, very thick] (startvertex) -- node[anchor=west] {} (exchanges);
				
				\draw [->, very thick, align=left] (exchanges) -- node[anchor=west] {} (criteria);
				
				\draw [->, very thick, align=left] (criteria) -- node[anchor=south] {Yes} (improvement);
				
				\draw [red, ->, very thick, align=left] (criteria) -- node[anchor=south] {No!} (negGain);
				
				\draw [red, ->, very thick, align=left] (negGain) |-  + (0,-1) -- ++(-5,-1) node[anchor=north east] {Yes!} -| (improvement);
				
				\draw [->, very thick, align=left] (improvement) + (-0.9,0.6) |- node[anchor=south west] {Yes, found\\ better tour $T'$\\ and set $T\leftarrow T'$} (startvertex);
				
				\draw [->, very thick, align=left] (improvement) + (0.0,0.6)  |- node[anchor= south west] {No, continue} (exchanges);
				
				\draw [dashed, <->, very thick, align=left] ([yshift=0.3cm]nonseq.west) --  node[anchor=south] {Additionally} ([yshift=0.3cm]exchanges.east);
				
				\draw [->, very thick, align=left] (exchanges.east) -|  + (2,0) |- + (4.5,-1) |- node[ anchor=south west]{No feasible\\ edge} (allvertices);
				
				\draw [->, very thick, align=left] (allvertices) -- node[anchor=east] {Yes} (finish);
				
				\draw [->, very thick, align=left] (allvertices)  + (1,-0.55) |-  + (3,-1) -- ++(9cm,-1) node[anchor=south east] {No} |- (startvertex);
				
				\draw [->, very thick] (finish) --  node[anchor=south] {} (newrun);
				
				\draw [->, very thick] (newrun) --
				++(-2.5cm,0) |-  node[anchor=south] {} (starttour);
			\end{tikzpicture}
		}
		\caption{Simplified procedure of improved LKH procedure with relaxed positive gain criterion displayed in red. Figure adjusted from \cite{thesis_BO}.
		}
		\label{flowchartLKimproved}
	\end{center}
\end{figure}

The criteria mentioned above regulate the edge exchanges that are the core of the heuristic. We describe an outer framework determining how and how often to repeat the search.
For the original LKH heuristic, consider the black features in \cref{flowchartLKimproved}:
We provide a simplified visualization of the procedure of one \struc{run} which is a restart of the entire heuristic.
Each run is independent from previous runs and consists of $n = |V|$ \struc{trials}.
A trial of LKH consists of improving a tour $T$ by repeatedly searching for $k$-opt-moves via sequential edge exchanges as explained above or additional non-sequential edge exchanges.
Once we find a tour $T'$ that has lower cost than the original tour $T$, we set $T\leftarrow T'$.
Then we repeat the procedure.
If we can not find a better tour, then we choose a new starting vertex $t_1$.
A trial stops when we have tried all vertices in $V$ as the starting vertex $t_1$ without finding a better tour.
The next trial begins with a new starting tour that partly depends on the edges from the current best tour of the run, see \cite{LKH_originalTR} for a description.

\section{Relaxing the positive gain criterion}
\label{refinements}
We aim to change the rules of the LKH procedure in a way such that it can find more alternating circles and thereby reach better TSP solutions. 
To this end we relax the hitherto untouched positive gain criterion \ref{C2:gain}. 
Our idea is that in particular for large instances alternating circles which have non-positive gain in the beginning can still be beneficial in the end. 

Pursuing this idea leads to the following trade-off. 
Relaxing the positive gain criterion enlarges the possibilities to find edge exchanges when constructing an alternating circle. 
Therewith we enlarge the search space of alternating circles that the procedure can find. 
But, eventually, we want the complete alternating circle to have total positive gain. 
Allowing for long series of edge exchanges with non-positive gain, while constructing the circle, might lead into situations where total positive gain in the end is hardly attainable. 
There is a trade-off between enlarging the search space of findable paths by relaxing the positive gain criterion and the chance to eventually close the circle with total positive gain.

For this reason we relax the positive gain criterion in a very careful way. Namely, we still require that if the gain is non-positive after the $i$th exchange it has to be positive after the $(i+1)$st. 
We call this the \struc{\homogeneous relaxation \ref{C2*:newgain}} of the positive gain criterion.
The computational results in \cref{section:experiments} show clearly that \ref{C2*:newgain} is a better trade-off than the original positive gain criterion \ref{C2:gain}.
On average the running times for large instances go down by 2.6 hours or 10.4\% and the quality goes up by 0.007\% compared to the best LKH version\footnote{For large instances Helsgaun recommends a POPMUSIC candidate set. See \cref{fig:CR_large_instances,%
		fig:CR_large_instances_fully_relaxed}.\label{popmusic_fully_relaxed}}.
 
Still, the \homogeneous relaxation is a stipulated choice. Therefore, we applied an algorithm engineering type circle to search for better solutions to the trade-offs. 
In our algorithm engineering we translate the insights from the computational study into a new variant for the positive gain criterion, and then test it again computationally. 
This leads us to a second variant of the positive gain criterion, namely, the \struc{\tilted relaxation \ref{C2**:tiltedgain}} of the positive gain criterion.
Here we restrict the \homogeneous relaxation by not allowing a violation of the positive gain criterion \ref{C2:gain} for certain exchanges in a $k$-opt-move.
It turns out in the experiments for large instances that this slight adjustment causes additional 0.8 hours or 3.2 percentage points decrease in running time on average, leading to a total decrease of \textbf{13.6\%}, and the same quality of solution\footref{popmusic_fully_relaxed}.

We first formally describe the \homogeneous relaxation \ref{C2*:newgain} of the positive gain criterion in \cref{subsection:relaxation} and then the \tilted relaxation \ref{C2**:tiltedgain} as an extension of \ref{C2*:newgain} in \cref{subsection:algorithm engineering}. 
\subsection{The \homogeneous relaxation of the positive gain criterion}
\label{subsection:relaxation}
In \cite{LKH_original} and \cite{LK_original} the authors argue that if $G_i = \sum_{\ell=1}^i g_{\ell} > 0$, then there exists a permutation $\pi$ such that every partial sum of $G_i = \sum_{\ell = 1}^ig_{\pi(\ell)}$ is positive, justifying the positive gain criterion \ref{C2:gain}.
However, being restricted to the candidate sets might not always permit finding this permutation.
Therefore, our idea is to relax the positive gain criterion by allowing intermediate steps with non-positive gain.
This increases the number of edges within the candidate set that satisfy all criteria.
We define the \struc{\homogeneous relaxation} of the positive gain criterion:
\begin{enumerate}[label=(C\arabic**)]
	\addtocounter{enumi}{1}
	\item Recall the gain $G_i := \sum_{\ell = 1}^{i} c(x_{\ell}) - c(y_{\ell})$ is the sum of cost differences after exchanging the first $i$ edges in a $k$-opt-move.
	By \ref{C2:gain} we may add an edge $y_i$ after deleting $x_i$ if $G_i >0$.
	\textbf{For sequential edge exchanges, we additionally permit adding $\mathbf{y_i}$ after deleting $\mathbf{x_i}$ with $\mathbf{G_i \leq 0}$ if $\mathbf{G_{i-1}>0}$.} 
	In particular, we assume $G_0$ to be positive. 
	%From now on, we refer to the relaxed criterion simply as \struc{gain criterion}. 
	\label{C2*:newgain}
\end{enumerate}
%We make no further adjustments. 
Restricting the non-positive gain $G_i \leq 0$ to the case of the previous gain $G_{i-1}$ being positive enlarges the number of feasible edges in the candidate set.
Note that we do not allow consecutive edge exchanges with non-positive gain.
Moreover, the termination criterion still ensures positive gain in the end.

Observe the red features in \cref{flowchartLKimproved} as visualization: 
Now, we additionally consider a previously prohibited sequential edge exchange feasible if $G_i \leq 0$ and $G_{i-1} > 0$.
To obtain an improved tour, the overall gain of course needs to be positive.

The following example displayed in \cref{lkhimprovedexample.2} shows that the \homogeneous relaxation \ref{C2*:newgain} gives rise to more possibilities for feasible tours:
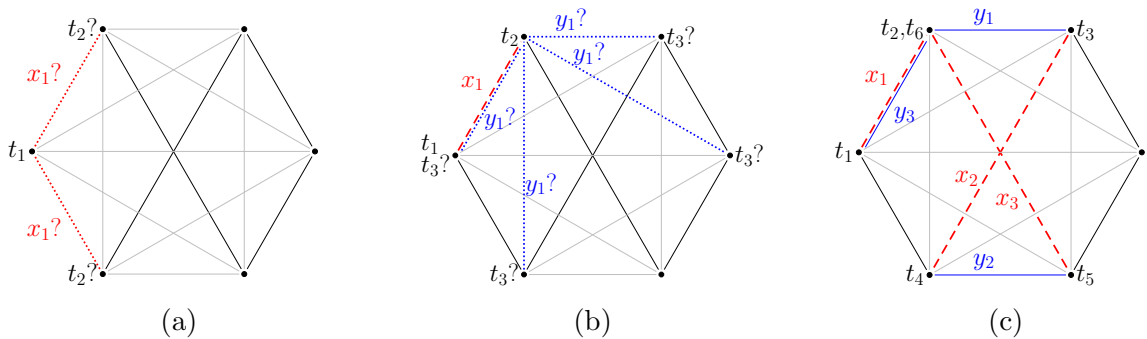
\begin{figure}[t] %Example of LKH with our gain criterion
	\begin{center}
		\begin{subfigure}[b]{0.3\textwidth}
			\centering
			\resizebox{\linewidth}{!}{
				\begin{tikzpicture}
					\node (a) at (0,0) {};
					\node (b) at (4,0) {};
					\node (c) at (6, -3.46) {};
					\node (d) at (4,-6.93) {};
					\node (e) at (0,-6.93) {};
					\node (f) at (-2,-3.46) {};
					\node (s1) at (4,-7) {}; %to shift image
					\node (s2) at (7,-3.46) {}; %to shift image
					
					%length 6
					\draw[thick]  (a) edge (d);
					\draw[thick]  (b) edge (e);
					\draw[thick, lightgray]  (c) edge (f);
					
					%length 5
					\draw[thick, lightgray]  (a) edge (c);
					\draw[thick, lightgray]  (b) edge (d);
					\draw[thick, lightgray]  (c) edge (e);
					\draw[thick, lightgray]  (d) edge (f);
					\draw[thick, lightgray]  (e) edge (a);
					\draw[thick, lightgray]  (f) edge (b);
					
					%length 4
					\draw[thick, lightgray]  (a) edge (b);
					\draw[thick, lightgray]  (d) edge (e);
					
					%length 3
					\draw[thick]  (b) edge  (c);
					\draw[thick]  (c) edge  (d);
					\draw[dotted, red, font=\fontsize{20}{0}\selectfont, line width=0.05cm]  (e) edge node[below left, red]{$x_1?$} (f);
					\draw[dotted, red, font=\fontsize{20}{0}\selectfont, line width=0.05cm]  (f)  edge node[above left, red]{$x_1?$} (a);
					
					%Nodes above edges
					\draw[fill] (0,0) circle (2pt) node[left, thick, font=\fontsize{20}{0}\selectfont] {$t_2$?};
					\draw[fill] (4,0) circle (2pt);
					\draw[fill] (6, -3.46) circle (2pt);
					\draw[fill] (4,-6.93) circle (2pt);
					\draw[fill] (0,-6.93) circle (2pt)  node[left, thick, font=\fontsize{20}{0}\selectfont] {$t_2$?};
					\draw[fill] (-2,-3.46) circle (2pt) node[left, thick, font=\fontsize{20}{0}\selectfont] {$t_1$};	
				\end{tikzpicture}
			}
			\caption{} 
			\label{lkhimproved2.2}
		\end{subfigure}
		\hspace{1em}
		\begin{subfigure}[b]{0.3\textwidth}
			\centering
			\resizebox{\linewidth}{!}{
				\begin{tikzpicture}
					\node (a) at (0,0) {};
					\node (b) at (4,0) {};
					\node (c) at (6, -3.46) {};
					\node (d) at (4,-6.93) {};
					\node (e) at (0,-6.93) {};
					\node (f) at (-2,-3.46) {};
					\node (s1) at (4,-7) {}; %to shift image
					\node (s2) at (6.5,-3.46) {}; %to shift image
					\node (a2) at (0, -0.173) {};
					\node (f2) at (-1.9, -3.46) {};
					
					%length 6
					\draw[thick]  (a) edge (d);
					\draw[thick]  (b) edge (e);
					\draw[thick, lightgray]  (c) edge (f);
					
					%length 5
					\draw[dotted, blue, font=\fontsize{20}{0}\selectfont, line width=0.05cm]  (a) edge (c);
					\node (h3)  [thick, blue, font=\fontsize{20}{0}\selectfont] at (1.9, -0.5) {$y_1$?};
					\draw[thick, lightgray]  (b) edge (d);
					\draw[thick, lightgray]  (c) edge (e);
					\draw[thick, lightgray]  (d) edge (f);
					\draw[dotted, blue, font=\fontsize{20}{0}\selectfont, line width=0.05cm]  (e) edge (a);
					\node (h5)  [thick, blue, font=\fontsize{20}{0}\selectfont] at (0.5, -4.3) {$y_1$?};
					\draw[thick, lightgray]  (f) edge (b);
					
					%length 4
					\draw[dotted, blue, font=\fontsize{20}{0}\selectfont, line width=0.05cm]  (a)  edge node[above left, blue]{$y_1$?} (b);
					\draw[thick, lightgray]  (d) edge (e);
					
					%length 3
					\draw[thick]  (b) edge  (c);
					\draw[thick]  (c) edge  (d);
					\draw[thick]  (e) edge (f);
					\draw[dashed, red, thick, font=\fontsize{20}{0}\selectfont, line width=0.05cm, dash pattern=on 8pt off 5pt]  (f)  edge node[above left, red]{$x_1$} (a);
					\draw[dotted, blue, font=\fontsize{20}{0}\selectfont, line width=0.05cm]  (f2)  edge (a2);
					\node (h6) [thick, blue, font=\fontsize{20}{0}\selectfont] at (-0.7, -2.4) {$y_1$?};
					
					%Nodes above edges
					\draw[fill] (0,0) circle (2pt) node[left, thick, font=\fontsize{20}{0}\selectfont, align = left] {$t_2$};
					\draw[fill] (4,0) circle (2pt)  node[right, thick, font=\fontsize{20}{0}\selectfont] {$t_3$?};
					\draw[fill] (6, -3.46) circle (2pt)  node[right, thick, font=\fontsize{20}{0}\selectfont] {$t_3$?};
					\draw[fill] (4,-6.93) circle (2pt);
					\draw[fill] (0,-6.93) circle (2pt)  node[left, thick, font=\fontsize{20}{0}\selectfont] {$t_3$?};
					\draw[fill] (-2,-3.46) circle (2pt) node[left, thick, font=\fontsize{20}{0}\selectfont, align=left] {$t_1$\\$t_3$?};	
				\end{tikzpicture}
			}
			\caption{} 
			\label{lkhimproved2.3}
		\end{subfigure}
		\hspace{1em}
		\begin{subfigure}[b]{0.3\textwidth}
			\centering
			\resizebox{\linewidth}{!}{
				\begin{tikzpicture}
				\node (a) at (0,0) {};
				\node (b) at (4,0) {};
				\node (c) at (6, -3.46) {};
				\node (d) at (4,-6.93) {};
				\node (e) at (0,-6.93) {};
				\node (f) at (-2,-3.46) {};
				\node (s1) at (4,-7) {}; %to shift image
				\node (s2) at (7,-3.46) {}; %to shift image
				\node (a2) at (0, -0.173) {};
				\node (f2) at (-1.9, -3.46) {};

				%length 6
				\draw[dashed, red, thick, font=\fontsize{20}{0}\selectfont, line width=0.05cm, dash pattern=on 8pt off 5pt]  (a) edge (d);
				\node (h2)  [thick, red, font=\fontsize{20}{0}\selectfont] at (2.2, -4.85) {$x_3$};
				\draw[dashed, red, thick, font=\fontsize{20}{0}\selectfont, line width=0.05cm, dash pattern=on 8pt off 5pt]  (b) edge (e);
				\node (h1)  [thick, red, font=\fontsize{20}{0}\selectfont] at (1.05, -4.2) {$x_2$};
				\draw[thick, lightgray]  (c) edge (f);
				
				%length 5
				\draw[thick, lightgray]  (a) edge (c);
				\draw[thick, lightgray]  (b) edge (d);
				\draw[thick, lightgray]  (c) edge (e);
				\draw[thick, lightgray]  (d) edge (f);
				\draw[thick, lightgray]  (e) edge (a);
				\draw[thick, lightgray]  (f) edge (b);
				
				%length 4
				\draw[thick, blue, font=\fontsize{20}{0}]  (a)  edge node[above left, blue]{$y_1$} (b);
				\draw[thick, blue, font=\fontsize{20}{0}] (d) edge node[above left, blue]{$y_2$} (e);
				
				%length 3
				\draw[thick]  (b) edge  (c);
				\draw[thick]  (c) edge  (d);
				\draw[thick]  (e) edge (f);
				\draw[dashed, red, thick, font=\fontsize{20}{0}\selectfont, line width=0.05cm, dash pattern=on 8pt off 5pt]  (f)  edge node[above left, red]{$x_1$} (a);
				\draw[thick, blue, font=\fontsize{20}{0}]  (f2)  edge (a2);
				\node (h6) [thick, blue, font=\fontsize{20}{0}\selectfont] at (-0.7, -2.4) {$y_3$};
			
				%Nodes above edges
				\draw[fill] (0,0) circle (2pt) node[left, thick, font=\fontsize{20}{0}\selectfont, align = left] {$t_2$,$t_6$};
				\draw[fill] (4,0) circle (2pt)  node[right, thick, font=\fontsize{20}{0}\selectfont] {$t_3$};
				\draw[fill] (6, -3.46) circle (2pt);
				\draw[fill] (4,-6.93) circle (2pt) node[right, thick, font=\fontsize{20}{0}\selectfont] {$t_5$};
				\draw[fill] (0,-6.93) circle (2pt)   node[left, thick, font=\fontsize{20}{0}\selectfont] {$t_4$};
				\draw[fill] (-2,-3.46) circle (2pt) node[left, thick, font=\fontsize{20}{0}\selectfont] {$t_1$};	
				\end{tikzpicture}
			}
			\caption{} 
			\label{lkhimproved2.4}
		\end{subfigure}
		
		\caption{Display of \cref{exampleLKHimproved}: Relaxed positive gain criterion. 
			Black edges belong to starting tour. Red edges indicate removal. Blue edges indicate insertion. In \cref{lkhimproved2.2} we choose one of the red dotted edges to be removed, in \cref{lkhimproved2.3} we select one of the blue dotted edges to be added.} 
		\label{lkhimprovedexample.2}
	\end{center}
\end{figure}
\begin{examplex}
	\label{exampleLKHimproved}
	We use the starting tour $T$ in \cref{starttour} but select a starting vertex $t_1$ that is different from the one in \cref{exampleLKHnormal}.
	In \cref{lkhimproved2.2} there are two edges incident to $t_1$ that belong to $T$. Thus, we have two choices for edges that we could delete as $x_1$. 
	After deleting the red dashed \color{red}$x_1$ \color{black} in \cref{lkhimproved2.3}, none of the blue dotted edges that are incident to $t_2$ could serve as $y_1$ by the positive gain criterion \ref{C2:gain}, as every gain would be non-positive.
	The same holds for the other choice of $x_1$ in \cref{lkhimproved2.2}.
	
	Now we consider the \homogeneous relaxation \ref{C2*:newgain}:
	Assume that all blue dotted edges in \cref{lkhimproved2.3} belong to the candidate set of $t_2$.
	Searching through the candidate set in a determined ranking and checking the criteria \ref{C4:sequential}, \ref{C5:disjunctivity} as well as the \homogeneous relaxation \ref{C2*:newgain} we may select the blue solid edge in \cref{lkhimproved2.4} as $\color{blue}y_1$\color{black}.
	Exchanging the edges \color{red}$x_1$ \color{black} and \color{blue}$y_1$ \color{black} leads to negative gain $G_1 = c(\color{red}x_1\color{black}) - c(\color{blue}y_1\color{black}) =\color{red}3\color{black}-\color{blue}4\color{black} = -1$.
	Given that we deleted \color{red}$x_2$ \color{black} in \cref{lkhimproved2.4} we now have to find $y_2$. 
	Only \color{blue} $y_2$ \color{black} as in \cref{lkhimproved2.4} satisfies all criteria \ref{C4:sequential}, \ref{C5:disjunctivity}, \ref{C2*:newgain} because we now demand a positive gain $G_2 = -1 + c(\color{red}x_2\color{black}) - c(\color{blue}y_2\color{black}) = -1 + \color{red}6 \color{black}- \color{blue}4 \color{black}= 1$.
	Finally, we close the tour with positive gain $G_3 = 1 + c(\color{blue}y_3\color{black}) - c(\color{red}x_3\color{black}) = 1 +\color{red} 6 \color{black} - \color{blue}3 \color{black} = 4$ in \cref{lkhimproved2.4}.
	We thus obtain the optimal tour as in \cref{exampleLKHnormal}.
	
	In this example, using the \homogeneous relaxation \ref{C2*:newgain} of the positive gain criterion leads to the optimal tour for \emph{every} starting vertex.
	In contrast, if using the original positive gain criterion \ref{C2:gain}, only four of the six vertices lead to the optimal tour.
	Thus, this example demonstrates that modifying the positive gain criterion can lead to more possibilities for finding tours.
\end{examplex}
\begin{figure}
\begin{algorithm}[\texttt{AdjustedFindSequentialEdgeToBeAdded}]~
	\label{alg:LKH_0.1_Gain_Criterion}
	
	\noindent \INPUT{Alternating path $P = [t_1, ..., t_{2i}]$, gain $G_{i-1}$ of $[t_1,...t_{2i-1}]$}\\
	\OUTPUT{\enquote{Found} / \enquote{Not found}, alternating path $P = [t_1, ..., t_{2i}, t_{2i+1}]$}
	\begin{algorithmic}[1] 
		\For{$(t_{2i},t_{2i+1})$ in \texttt{CandidateSet}($t_{2i}$)} 
		\State $G_i \gets G_{i-1} + c(t_{2i-1}, t_{2i}) - c(t_{2i}, t_{2i+1})$
		\If{ \textcolor{red}{(}$G_i > 0$ \textcolor{red}{\textbf{or} $G_{i-1} > 0$)} \textbf{and} $\texttt{OtherCriteriaSatisfied}(P, (t_{2i}, t_{2i+1}))$} \label{line:alg_0.1_if}
		\State \Return \enquote{Found}, $P = [t_1, ..., t_{2i}, t_{2i+1}]$ 
		\EndIf	
		\EndFor
		\State \Return \enquote{Not found}
	\end{algorithmic}	
\end{algorithm}
\caption*{Simplified example of the usage of the \homogeneous relaxation of the positive gain criterion in our modified code based on the example code in \cref{alg:LKH_Gain_Criterion}. We highlight the modified part in red.}
\end{figure}
\cref{alg:LKH_0.1_Gain_Criterion} is our adjusted version of \cref{alg:LKH_Gain_Criterion} for selecting an edge $y_i$ to be added using \ref{C2*:newgain} instead of \ref{C2:gain}. 
We highlight our changes in \color{red}red\color{black}.
As in the original LKH we iterate over all edges in the candidate set of vertex $t_{2i}$. 
Using the if-clause in \cref{line:alg_0.1_if} of the function \texttt{AdjustedFindSequentialEdgeToBeAdded} we allow $G_i$ to be negative if $G_{i-1}$ is positive.
Following the ranking of the candidate set, we want to select the first edge $(t_{2i}, t_{2i+1})$ as $y_i$ that satisfies the \homogeneous relaxation \ref{C2*:newgain} as well as criteria
\ref{C3:feasibility}, \ref{C4:sequential} and \ref{C5:disjunctivity}, summarized as function \texttt{OtherCriteriaSatisfied}. 
If we find such a feasible $y_i$, then we add it to the alternating path $P$.
Otherwise we return ``Not found''.
In the latter case, $P$ does not lead to a feasible $k$-opt-move, and we have to discard it.
Such failure to enlarge the alternating path until it can be closed as a beneficial alternating circle happens less often with a relaxed positive gain criterion \ref{C2*:newgain} compared to \ref{C2:gain}.

While the pseudo code above provides the idea of our \homogeneous relaxation \ref{C2*:newgain}, we present the core of our variant as C-code in \cref{lst:CodeBest5opt} in the Appendix.

\subsection{Further improvement by algorithm engineering}
\label{subsection:algorithm engineering}
The previously described \homogeneous relaxation enlarges the search space of alternating circles. By algorithm engineering we developed the \tilted relaxation \ref{C2**:tiltedgain} of the positive gain criterion. Recall the gain $G_i := \sum_{\ell = 1}^{i} c(x_{\ell}) - c(y_{\ell})$ is the sum of cost differences after exchanging the first $i$ edges in a $k$-opt-move.
	By \ref{C2:gain} we may add an edge $y_i$ after deleting $x_i$ if $G_i >0$.
	For sequential edge exchanges, we additionally permit adding $y_i$ after deleting $x_i$ with $G_i \leq 0$ if $G_{i-1}>0$.
	Technically, for \ref{C2*:newgain} we assume $G_0$ to be positive. 
	The \tilted relaxation \ref{C2**:tiltedgain} restricts the enlarged search space of the \homogeneous relaxation \ref{C2*:newgain} by the following two ideas inspired by our algorithm engineering:

\begin{enumerate}
\item In the \homogeneous relaxation, we allow $G_{i} \leq 0$ for every $i$ if $G_{i-1} > 0$.
In the \tilted relaxation we restrict these $i$ to those that satisfy one of the following conditions: 
\begin{enumerate}
	\item $i-1$ is not a multiple of $k$.
	\item The positive gain criterion \ref{C2:gain} was satisfied for $i-2$,  i.e. $G_{i-2} > 0$.
\end{enumerate}
\item Our implementation of \ref{C2*:newgain} always assumes $G_0$ to be positive.
Thus, we allow $G_1 \leq 0$. 
In the \tilted relaxation we only assume $G_0$ to be positive for the first exchange in the first trial of the first run.
For later $k$-opt-moves, trials and runs, we assume $G_0$ to have the same sign as the last $G_i$ or, in case $i$ is a multiple of $k$, as the last $G_{i-1}$.
\end{enumerate}
By condition (2) the positivity of an exchange may depend on the positivity of previous $k$-opt-moves, trials and runs. 
Summarizing these ideas we define the \struc{\tilted relaxation} of the positive gain criterion as follows: 

\begin{enumerate}[label=(C\arabic***)]
	\addtocounter{enumi}{1}
	\item For sequential edge exchanges, we permit adding $y_i$ after deleting $x_i$ with $G_i \leq 0$ if $G_{i-1}>0$ and either $i-1$ is not a multiple of $k$ or  $G_{i-2}>0$. 
	Additionally, we assume $G_0$ to have the same sign as the last $G_i$ or,  in case $i$ is a multiple of $k$, as the last $G_{i-1}$.
	If no previous $G_i$ exists, then we assume $G_0$ to be positive. 
	\label{C2**:tiltedgain}
\end{enumerate}

See \cref{subsection:algo_engineering_results} for an empirical comparison of the \homogeneous relaxation and the \tilted relaxation.

\section{Experiments}
\label{section:experiments}
In this section we evaluate our modification of the LKH heuristic by comparing the results with LKH on standard benchmark sets. 
First, we describe the graphs of the benchmark instances and then our experimental setting.
Finally, we present and assess our computational results.

\subsection{Instances}
\label{subsec:Ins}
We use different datasets of graphs which are commonly used for empirical analysis of TSP-heuristics or algorithms to compare our modification to LKH. 
We provide an overview of these datasets in \cref{tab:overview_instances}. 
We distinguish between \struc{small}, \struc{medium} and \struc{large} instances depending on the number of nodes in the graph. 
Small instances have 1000 or less nodes, medium instances between 1001 and 30000 nodes, and large instances more than 30000 nodes. 
We chose the limit for the medium instances because we were able to solve almost every such instance in under two days 
\footnote{The following instance was not solved in two days with either an ALPHA or a POPMUSIC generated candidate set: National\_TSPs\_medium\_ho14473.par. We excluded this instance from further considerations for the medium instances.}.
\begin{table}
\caption{Overview of the graphs with less than one million nodes used in our benchmark computations.
}
\begin{tabular}{p{5cm}rrrll}
instance set & \multicolumn{3}{c}{\# instances} & source\\
&small&medium&large&\\
\hline
Symmetric TSP-Library & 78&31&2 & \cite{TSPLIB}\\
Asymmetric TSP-Library & 19&0&0 & \cite{TSPLIB}\tablefootnote{Asymmetric instances are transformed into symmetric instances during preprocessing \cite{LKH_original}. A POPMUSIC candidate set is not applicable for asymmetric problems.}\\
VLSI (Dataset) & 15&72&15 &\cite{vlsi} \\
Randomized Instances (8th DIMACS challenge) & 24&19&10 &Used first in \cite{DIMACS}\\
National TSPs & 6&19&2 & \cite{nationalTSP}\\
$T_{nm}$-Instances & 123 & 4 & 0 & Constructed in \cite{Hougardy_Tnm}\\
\end{tabular}
\label{tab:overview_instances}
\end{table}

\subsection{Experimental design}
\label{subsec:Exp_Des}
The entire experiment is based on the software by Keld Helsgaun \cite{LKH3} which was written in C. 
We use the version LKH 3.0.8 and refer to this code as \struc{(original) LKH}.
We introduce our new variables in files \enquote{LKH.h} and \enquote{ReadParameters.c}. 
Our modification is implemented in \enquote{BestKOptMove.c} and \enquote{Best5KOptMove.c}.  
%Our adjusted code for \ref{C2**:tiltedgain} can be found in \url{...}\footnote{
%Our code will be available soon.}.
We execute our computations on a cluster with one node and 256 AMD EPYC 7742 64-Core processors. 
Every instance runs on one cpu with a frequency of 3315.439 Mhz. 
To track the runtime we measure the execution time. 
Regarding the LKH implementation we use the parameters presented in \cref{tab:LKH_parameters}. 
We apply the default paramters described in \cite{LKH_original} and \cite{LKH_popmusic} if not stated differently: 
\begin{itemize}
\item We perform 100 \struc{RUNS} for the small instances and ten runs for the medium and large instances. 
For the {$T_{nm}\text{-instances}$} we used deviating parameters according to the published results by \cite{Helsgaun_Tnm}. 
Therefore, the number of runs for $T_{nm}$-instances is set to 100 for all instances. 
\item In real-world applications we do not know the cost of the optimal tour in advance. 
If one knows this optimal value in advance and passes it to the algorithm, one can also set the parameter \struc{STOP\_AT\_OPTIMUM} to YES. 
Thus, a run stops when the cost of the current tour equals the given optimum. 
We set the paramenter STOP\_AT\_OPTIMUM to NO to be able to compare LKH and our variant in a near real-world environment. 
Thus, we also consider the \struc{checkout time}, i.e. the time spent searching for an improvement although one optimal tour has already been found.
This differs from Helsgaun's approach, who set this parameter to YES. 
\item Our candidate sets contain a maximum of 100 candidates for $T_{nm}$-instances and five candidates for the other instances. 
This parameter is called \struc{MAX\_CANDIDATES}. 
\item For most instances we provide results for two different \struc{CANDIDATE\_SET\_TYPES}, i.e. types of generating the candidate sets: POPMUSIC and ALPHA. 
If the candidate set is generated according to POPMUSIC, then we set \struc{POPMUSIC\_TRIALS} to 0 and \struc{POPMUSIC\_SAMPLE\_SIZE}  to 100 for the $T_{nm}$-instances. 
For the other instances we use the default values as given in \cref{tab:LKH_parameters}. 
See \cite{LKH_popmusic} for further explanation of the candidate set generation by POPMUSIC.
\end{itemize}
\begin{table}
	\caption{LKH parameters for the $T_{nm}$ instances according to \cite{Helsgaun_Tnm} and all other instances mostly according to \cite{LKH_original, LKH_popmusic}.}
	\begin{center}
		\begin{tabular}{lrr}
			Parameter & value for $T_{nm}$ & value for rest \\
			\hline
			RUNS & 100 & 10/100\\
			STOP\_AT\_OPTIMUM & NO & NO\\		
			MAX\_ CANDIDATES & 100 & 5\\
			CANDIDATE\_ SET\_ TYPE & POPMUSIC & POPMUSIC/ \linebreak ALPHA\\
			POPMUSIC\_TRIALS & 0 & 1/-\\
			POPMUSIC\_SAMPLE\_SIZE & 100 & 10/-
		\end{tabular}
	\end{center}
	\label{tab:LKH_parameters}
\end{table}

Since nine instances did not finish in one month of computation time with the original code or our variant or both, we compare the original LKH and our variant for these instances by giving both the same time limit for each run. 
To speed up the computation process we split the usual ten runs into two starts with five runs each.
We provide all parameters for this evaluation in \cref{tab:LKH_large_parameters}. 
\begin{itemize}
\item Setting the initial \struc{SEED} for random number generation ensures that for each variant the same initial tour is built for the different runs. 
One start is processed with seed 1 and the other start with seed 42.
\item We set the \struc{TIME\_LIMIT} to three days per run excluding the preprocessing time to ensure a reasonable overall computation time. 
\end{itemize} 
\begin{table}
	\caption{Second LKH parameters for the large instances that do not finish in one month.}
	\begin{center}
		\begin{tabular}{lr}
			Parameter & value \\
			\hline
			RUNS & 5\\			
			STOP\_AT\_OPTIMUM & NO\\		
			MAX\_ CANDIDATES & 5\\
			CANDIDATE\_ SET\_ TYPE & POPMUSIC/ \linebreak ALPHA\\
			SEED & 1, 42\\			
			TIME\_LIMIT & 259200s
		\end{tabular}
	\end{center}
	\label{tab:LKH_large_parameters}
\end{table}
\subsection{Computational results}
In this subsection we compare the results of the instances described in \cref{subsec:Ins} with the setting given in \cref{subsec:Exp_Des} for our code and the original code of Helsgaun.
For every instance and run we track the computation time of our variant as well as the original LKH. 
\begin{itemize}
\item We refer to the average time per run of one instance as \struc{TimeAvg}.
\item The original LKH and our variant produce a feasible tour in every run by construction.
We refer to the minimal cost after all runs of one instance as \struc{CostMin}. 
The average cost over all runs of one instance is denoted as \struc{CostAvg}.
\item If the optimum for the given instance is known, then we also track the gap between the optimum and the cost of the resulting tour in this run. Analogous to the cost, we call the minimal gap after all runs \struc{GapMin} and the average gap \struc{GapAvg}.
\end{itemize}
The following computational results show that our modification substantially decreases the computation time in comparison to the original LKH with a growing number of nodes. 
The significant time savings occur for large instances.
For completeness we also include the results for the small and medium instances. 
We display the average values of the parameters mentioned above for small, medium and large instances in \cref{fig:CR_small_instances,%
fig:CR_medium_instances,%
fig:CR_large_instances,%
fig:CR_large_running_instances_POPMUSIC,%
fig:CR_large_running_instances_alpha,%
fig:CR_large_stop_instances_POPMUSIC,%
fig:CR_large_stop_instances_alpha}.
We provide analogous tables in the appendix in \cref{fig:CR_median_large_instances,%
fig:CR_median_medium_instances,%
fig:CR_median_small_instances}
and \cref{fig:CR_quantil_large_instances,%
fig:CR_quantil_medium_instances,%
fig:CR_quantil_small_instances} calculating the median and the 95-th quantil respectively.
See \cref{fig:CR_large_instances_POPMUSIC_direct,%
	fig:CR_large_instances_ALPHA_direct,%
fig:CR_medium_instances_ALPHA_direct,%
fig:CR_medium_instances_POPMUSIC_direct,%
fig:CR_small_instances_ALPHA_direct,%
fig:CR_small_instances_POPMUSIC_direct,%
fig:CR_medium_instances_POPMUSIC_direct_100,%
fig:CR_small_instances_POPMUSIC_direct_100} in the appendix for a direct comparison of the instances.
Green and bold indicate a better performance of our variant compared to Helsgaun's original version, red and italic indicate a worse performance of our variant.

\textbf{Small and medium instances.}
In \cref{fig:CR_small_instances} we display the average over all small instances. 
The difference of quality between the results of the original code and our variant is not significant. 
Our variant is on average half a second slower per run than the original LKH.
Depending on the candidate set the average gap is slightly in favor of the original LKH or our variant. 
The minimal gap is in average the same for all instances. 

\begin{table}
\caption{Comparison of our variant using the \tilted relaxation \ref{C2**:tiltedgain} and the original LKH in minimal gap, average gap and average time per run on \textbf{small instances} with different candidate set types. Green and bold indicate a better performance of our variant, red and italic a worse performance of our variant. 
}
\begin{tabular}{lrrrrr}%
\hline
     Version & Candidates & \# Ex. & GapMin & GapAvg & TimeAvg
     \\ 
     &&& in \% & in \% & in s\\ \hline
    \begin{filecontents*}{Latex_Summary_Small.csv}
Version;CandidateModel;Candidates;Examples;Runs;Successes;TrialsMin;TrialsAvg;TrialsMax;Optimum;CostMin;CostAvg;CostMax;GapMin;GapAvg;GapMax;TimeMin;TimeAvg;TimeMax
LKH 3.0.8;ALPHA;5;142;95.07;84.34;405.5;405.5;405.5;2626438;2626438;2626675;2627681;0.0;6.64e-05;4.60e-04;1.512;3.234;6.399
LKH with \ref{C2**:tiltedgain};ALPHA;5;142;95.07;81.61;405.5;405.5;405.5;2626438;2626438;\textcolor{red}{\textit{2626846}};\textcolor{red}{\textit{2628161}};0.0;\textcolor{red}{\textit{7.97e-05}};\textcolor{green}{\textbf{4.17e-04}};\textcolor{red}{\textit{1.598}};\textcolor{red}{\textit{3.696}};\textcolor{red}{\textit{7.866}}
LKH 3.0.8;POPMUSIC;5;123;95.12;87.28;428.8;428.8;428.8;3031128;3031129;3031275;3031656;2.32e-03;7.70e-05;2.82e-04;1.897;3.502;6.860
LKH with \ref{C2**:tiltedgain};POPMUSIC;5;123;95.12;87.21;428.8;428.8;428.8;3031128;3031129;\textcolor{green}{\textbf{3031264}};\textcolor{red}{\textit{3032323}};2.32e-03;\textcolor{green}{\textbf{7.21e-05}};\textcolor{red}{\textit{2.95e-04}};\textcolor{red}{\textit{1.988}};\textcolor{red}{\textit{4.006}};\textcolor{red}{\textit{8.602}}
LKH 3.0.8;POPMUSIC;100;123;100.0;20.90;251.6;251.6;251.6;Error;4156602;4156666;4156689;0.481;4.82e-03;4.82e-03;0.311;0.509;0.748
LKH with \ref{C2**:tiltedgain};POPMUSIC;100;123;100.0;19.80;251.6;251.6;251.6;Error;\textcolor{green}{\textbf{4156587}};\textcolor{green}{\textbf{4156651}};\textcolor{green}{\textbf{4156670}};0.481;4.82e-03;4.82e-03;\textcolor{red}{\textit{0.561}};\textcolor{red}{\textit{0.934}};\textcolor{red}{\textit{1.454}}
\end{filecontents*}
\csvreader[head to column names, separator=semicolon]{Latex_Summary_Small.csv}{}
    { \Version & \Candidates\ \CandidateModel & \Examples & \GapMin & \GapAvg  & \TimeAvg\\}
    \\[-\normalbaselineskip]\hline
\end{tabular}
\label{fig:CR_small_instances}
\end{table}
\begin{table}
\caption{Comparison of our variant using the \tilted relaxation \ref{C2**:tiltedgain} and the original LKH in minimal gap, average gap and average time per run on \textbf{medium instances} with different candidate set types.}
\begin{tabular}{lrrrrr}%
\hline
     Version & Candidates & \# Ex. & GapMin & GapAvg & TimeAvg% specify table head
     \\ 
     &&&  in \% &  in \% & in s\\ \hline
    \begin{filecontents*}{Latex_Summary_Medium.csv}
Version;CandidateModel;Candidates;Examples;Runs;Successes;TrialsMin;TrialsAvg;TrialsMax;Optimum;CostMin;CostAvg;CostMax;GapMin;GapAvg;GapMax;TimeMin;TimeAvg;TimeMax
LKH 3.0.8;ALPHA;5;140;10.00;4.143;6957;6957;6957;Error;4863897;4864286;4865108;3.99e-03;1.75e-04;4.05e-04;1785;2129;2591
LKH with \ref{C2**:tiltedgain};ALPHA;5;140;10.00;3.857;6957;6957;6957;Error;4863897;\textcolor{red}{\textit{4864538}};\textcolor{red}{\textit{4865433}};\textcolor{red}{\textit{5.88e-03}};\textcolor{red}{\textit{2.60e-04}};\textcolor{red}{\textit{5.23e-04}};\textcolor{green}{\textbf{1693}};\textcolor{green}{\textbf{2043}};\textcolor{green}{\textbf{2499}}
LKH 3.0.8;POPMUSIC;5;140;10.00;4.743;6957;6957;6957;Error;4865171;4865546;4866138;0.068;8.27e-04;1.04e-03;2138;2550;3090
LKH with \ref{C2**:tiltedgain};POPMUSIC;5;140;10.00;4.557;6957;6957;6957;Error;\textcolor{green}{\textbf{4864873}};\textcolor{green}{\textbf{4865371}};\textcolor{green}{\textbf{4866057}};\textcolor{green}{\textbf{0.053}};\textcolor{green}{\textbf{7.56e-04}};\textcolor{green}{\textbf{1.02e-03}};\textcolor{green}{\textbf{2053}};\textcolor{green}{\textbf{2454}};\textcolor{green}{\textbf{2948}}
LKH 3.0.8;POPMUSIC;100;4;100.0;0.0;3501;3501;3501;Error;65404885;65406344;65407232;0.0;2.40e-05;4.10e-05;5.880;9.738;15.79
LKH with \ref{C2**:tiltedgain};POPMUSIC;100;4;100.0;0.0;3501;3501;3501;Error;\textcolor{green}{\textbf{65402841}};\textcolor{green}{\textbf{65405417}};\textcolor{green}{\textbf{65405788}};0.0;\textcolor{red}{\textit{3.90e-05}};\textcolor{red}{\textit{4.45e-05}};\textcolor{red}{\textit{6.543}};\textcolor{red}{\textit{13.80}};\textcolor{red}{\textit{25.31}}
\end{filecontents*}
\csvreader[head to column names, separator=semicolon]{Latex_Summary_Medium.csv}{}
    { \Version & \Candidates\ \CandidateModel & \Examples & \GapMin & \GapAvg & \TimeAvg\\}
    \\[-\normalbaselineskip]\hline
\end{tabular}
\label{fig:CR_medium_instances}
\end{table}
In \cref{fig:CR_medium_instances} we display the average over all medium instances.
On the one hand, the computation times of our variant averaged over all medium instances is roughly 90 seconds or four percent lower than the original LKH.
On the other hand, the average of the ratios between our and LKH's running times is 112\% where the average is taken over all medium instances with five POPMUSIC candidates.
Averaging over all medium instances with five ALPHA candidates this ratio is 115\%.
This discrepancy originates from our algorithm being slightly slower on various instances with shorter running times but substantially faster on few instances with longer running times,  see \cref{fig:CR_medium_instances_POPMUSIC_direct} and \cref{fig:CR_medium_instances_ALPHA_direct}. 
The solution quality of our variant is slightly worse or better than the one of the original LKH depending on the chosen candidate set type.
For the $T_{nm}$-instances with 100 candidates our variant achieves the same minimal and average gap as the original LKH. 
However, our variant spends on average four seconds or 42\% more on one run than the original code.
 
Overall for small and medium instances the solution quality and the running times are almost on the same level for both variants.\\
\begin{table}
\caption{Comparison of our variant using the \tilted relaxation \ref{C2**:tiltedgain} and the original LKH in minimal cost, average cost and average time per run on \textbf{large instances} with different candidate set types.}
\begin{tabular}{lrrrrr}%
\hline
     Version & Candidates & \# Ex. & CostMin & CostAvg & TimeAvg % specify table head
     \\ \hline
    \begin{filecontents*}{Latex_Summary_Large.csv}
Version;CandidateModel;Candidates;Examples;Runs;Successes;TrialsMin;TrialsAvg;TrialsMax;Optimum;CostMin;CostAvg;CostMax;GapMin;GapAvg;GapMax;TimeMin;TimeAvg;TimeMax
LKH 3.0.8;ALPHA;5;20;10.00;0.050;46688;46688;46688;Error;32842020;32856000;32882369;8.90e-03;3.47e-04;7.80e-04;65935;77171;111442
LKH with \ref{C2**:tiltedgain};ALPHA;5;20;10.00;0.0;46688;46688;46688;Error;\textcolor{red}{\textit{32856187}};\textcolor{red}{\textit{32880367}};\textcolor{red}{\textit{32930336}};\textcolor{green}{\textbf{8.02e-03}};\textcolor{red}{\textit{4.50e-04}};\textcolor{red}{\textit{1.18e-03}};\textcolor{green}{\textbf{51723}};\textcolor{green}{\textbf{61126}};\textcolor{green}{\textbf{78546}}
LKH 3.0.8;POPMUSIC;5;20;10.00;0.0;46688;46688;46688;Error;32826656;32832477;32837146;7.47e-03;2.05e-04;3.46e-04;80793;91276;103658
LKH with \ref{C2**:tiltedgain};POPMUSIC;5;20;10.00;0.0;46688;46688;46688;Error;\textcolor{green}{\textbf{32824384}};\textcolor{green}{\textbf{32832052}};\textcolor{red}{\textit{32841860}};\textcolor{green}{\textbf{6.38e-03}};\textcolor{red}{\textit{2.24e-04}};\textcolor{red}{\textit{4.35e-04}};\textcolor{green}{\textbf{69651}};\textcolor{green}{\textbf{78864}};\textcolor{green}{\textbf{88445}}
\end{filecontents*}
\csvreader[head to column names, separator=semicolon]{Latex_Summary_Large.csv}{}
    { \Version & \Candidates\ \CandidateModel & \Examples & \CostMin & \CostAvg & \TimeAvg \;s\\}
    \\[-\normalbaselineskip]\hline
\end{tabular}
\label{fig:CR_large_instances}
\end{table}
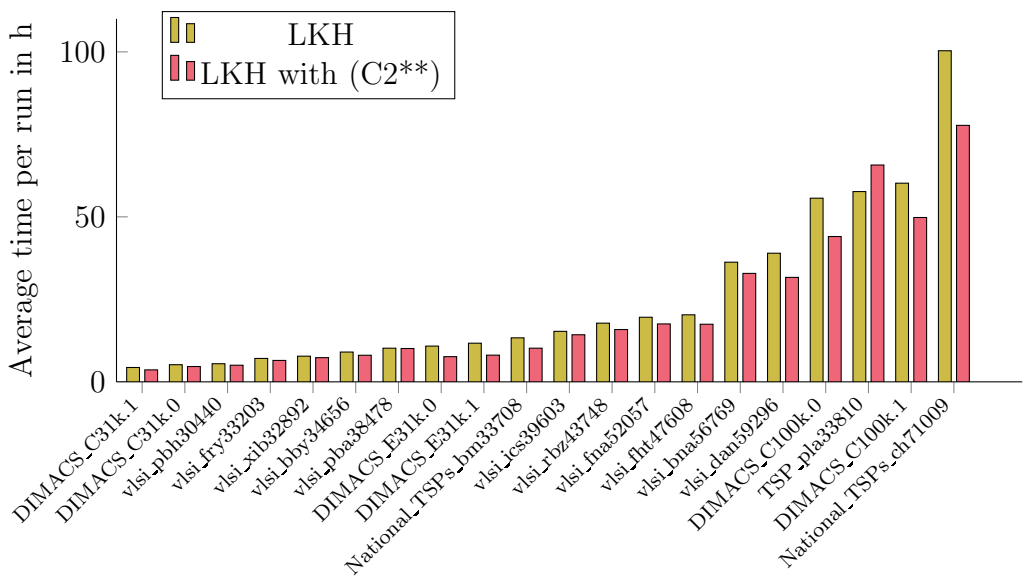
\begin{figure}
\centering 
\pgfplotstableread[col sep=semicolon, 
/pgf/number format/read comma as period%
]{test_large.csv}\datatable
\begin{tikzpicture} 
\begin{axis}
[
width=0.85\linewidth, % Scale the plot to \linewidth
height=0.4\linewidth,
xmin=0,
xmax=20,
ylabel={Average time per run in h},
ymin=0,
axis lines*=left,
enlarge x limits=0.03,
ymax=110,
x tick label style ={font=%\fontsize{4}{6}\selectfont,
\tiny, rotate=45,anchor=east},
ybar = 2pt,
bar width=0.3,
xtick=data,
xticklabels from table={\datatable}{Problem},
legend style={at={(0.05,0.9)},anchor=west}
] 
\addplot[fill=b2] table [y=TimeAvgLKHh, /pgf/number format/read comma as period,%
 col sep=semicolon, x expr=\coordindex]{\datatable};
\addplot[fill=b1] table [y=TimeAvgOurh, /pgf/number format/read comma as period,%
 col sep=semicolon, x expr=\coordindex]{\datatable};
\legend{LKH, LKH with \ref{C2**:tiltedgain}} 
\end{axis} 
\end{tikzpicture} 
\caption{Comparison of our variant using the \tilted relaxation \ref{C2**:tiltedgain} and the original LKH on large instances ($>30000$ vertices) in the average time per run with 5 POPMUSIC candidates. Instances that did not finish in one month of computation time are not displayed.}
\label{fig:AvgTime_large_instances}
\end{figure}
\textbf{Large instances.} For the large instances we do not compute the gap to the optimum because the optimum is not known for all instances.
Instead, we compare the minimal and the average cost. 
We distinguish between the large instances that finished their runs in one month and the rest. 
The averaged results of the 20 former instances in \cref{fig:CR_large_instances} show that with an ALPHA candidate set our variant using \tilted relaxation is about 4.5 hours or 20.8\% faster per run than the original LKH. 
In terms of average of ratios of running times this amounts to 17.8\%.
With POPMUSIC candidates the average time reduction are 3.4 hours or 13.6\% see \cref{fig:AvgTime_large_instances}.
The average of ratios is 13.1\% see \cref{fig:ratio_large_instances}.
With ALPHA candidate sets, the minimal cost achieved by our variant us compared to LKH increases on average 0.07\% over all large instances. 
With the POPMUSIC candidate set the minimal cost decreases on average 0.007\% over all large instances.
Note that Helsgaun recommends POPMUSIC for large-scale instances and all tested variants achieve significantly better tours with POPMUSIC candidate sets.
In \cref{subsection:Add_Tables_large_instances} one can see a detailed comparison of the minimal cost and the average time per run for each large instance and both candidate sets in \cref{fig:CR_large_instances_POPMUSIC_direct,%
fig:CR_large_instances_ALPHA_direct}.

Nine large instances were not solvable in one month.
Among these nine instances, six instances finished at least one run with a POPMUSIC candidate set and five instances with an ALPHA candidate set. 
We present these results in \cref{fig:CR_large_running_instances_POPMUSIC,%
fig:CR_large_running_instances_alpha}. 
The tables show the minimal cost and the average time per run over all finished runs. \\
We now give an example that exhibits the most significant time saving we could achieve:
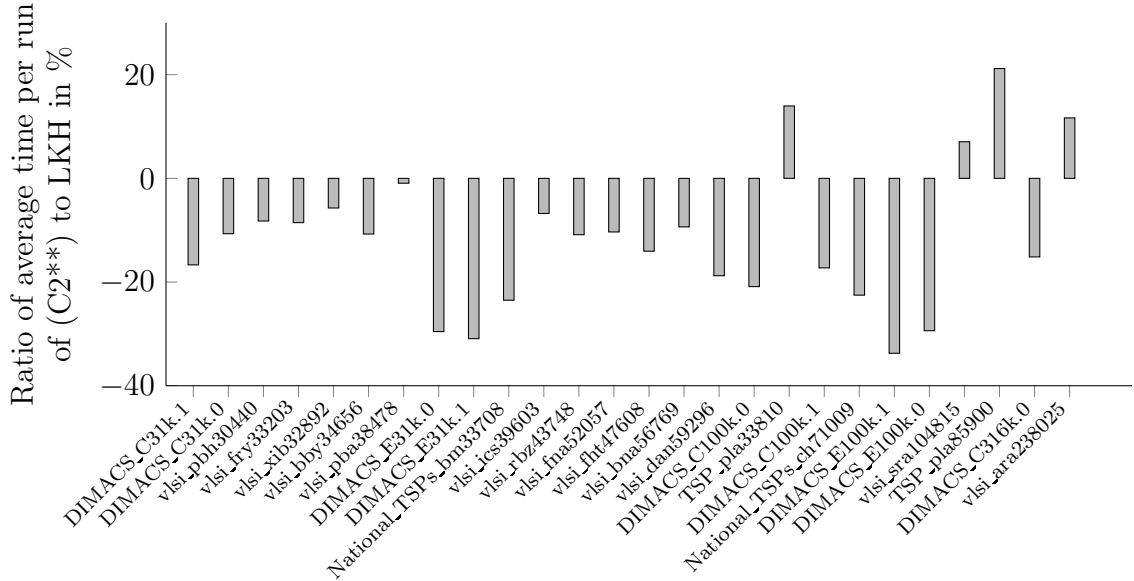
\begin{figure}
\centering 
\pgfplotstableset{col sep=semicolon}
\pgfplotstablevertcat{\unsortedtable}{test_large.csv} % loads `dataA.csv' -> `\output'
\pgfplotstablevertcat{\unsortedtable}{test_large_running.csv} % appends rows of dataB.csv
\pgfplotstablesort[sort key=TimeAvgLKH]{\datatable}{\unsortedtable}
\begin{tikzpicture} 
\begin{axis}
[
width=0.9\linewidth, % Scale the plot to \linewidth
height=0.4\linewidth,
xmin=0,
xmax=26,
ylabel style={align=center},
ylabel={Ratio of average time per run\\ of \ref{C2**:tiltedgain} to LKH in \%},
ymin=-40,
axis lines*=left,
enlarge x limits=0.03,
ymax=30,
x tick label style ={font=%\fontsize{4}{6}\selectfont,
\tiny, rotate=45,anchor=east}, 
ybar = 2pt,%Säulendiagramm
bar width=0.3,%Breite der Säulen
xtick=data,
xticklabels from table={\datatable}{Problem},
legend style={at={(0.05,0.95)},anchor=west}
] 
\addplot[fill=b7] table [y=TimeAvgChange, /pgf/number format/read comma as period,%
 col sep=semicolon, x expr=\coordindex]{\datatable};
\end{axis} 
\end{tikzpicture} 
\caption{Ratios of the average time per run of our variant using the  \tilted relaxation \ref{C2**:tiltedgain} to the average time per run of the original LKH on large instances ($>30000$ vertices) in the average time per run with 5 POPMUSIC candidates. Instances that did not finish in one month of computation time are also displayed.}
\label{fig:ratio_large_instances}
\end{figure}
\begin{table}
\caption{Direct comparison of our variant using the \tilted relaxation \ref{C2**:tiltedgain} and the original LKH with five POPMUSIC candidates in minimal cost and average time per run for all \textbf{large instances that did not finish in one month of computation time}.}
\begin{tabular}{lrrrrrr}%
\hline
 \multicolumn{1}{c}{Problem} & \multicolumn{2}{c}{CostMin} & %\multicolumn{2}{c}{CostAvg} & 
\multicolumn{2}{c}{TimeAvg in s} & \multicolumn{1}{c}{TimeAvg} \\ 
& LKH 3.0.8 & with \ref{C2**:tiltedgain} & LKH 3.0.8 & with \ref{C2**:tiltedgain} & \multicolumn{1}{c}{ratio in \%}\\ \hline% specify table head
    \begin{filecontents*}{Latex_Running_Comparison_Large_POPMUSIC.csv}
Problem;CandidateModel;Candidates;Optimum;RunsLKH;RunsOur;CostMinLKH;CostMinOur;CostAvgLKH;CostAvgOur;CostMaxLKH;CostMaxOur;GapMinLKH;GapMinOur;GapAvgLKH;GapAvgOur;GapMaxLKH;GapMaxOur;TimeMinLKH;TimeMinOur;TimeAvgLKH;TimeAvgOur;TimeMaxLKH;TimeMaxOur;CostMinChange;CostAvgChange;CostMaxChange;GapMinChange;GapAvgChange;GapMaxChange;TimeMinChange;TimeAvgChange;TimeMaxChange
vlsi\_ara238025;POPMUSIC;5;578761;2.000;2.000;579046;\textcolor{green}{\textbf{578931}};579134;\textcolor{green}{\textbf{579028}};579223;\textcolor{green}{\textbf{579126}};4.92e-04;\textcolor{green}{\textbf{2.94e-04}};6.45e-04;\textcolor{green}{\textbf{4.63e-04}};7.98e-04;\textcolor{green}{\textbf{6.31e-04}};4715563;\textcolor{red}{\textit{5350400}};5160046;\textcolor{red}{\textit{5761612}};5604529;\textcolor{red}{\textit{6172824}};-1.99e-02;-1.83e-02;-1.67e-02;-40.24;-28.22;-20.93;+13.46;+11.66;+10.14
TSP\_pla85900;POPMUSIC;5;142382641;10.00;8.000;142390914;\textcolor{red}{\textit{142392721}};142413210;\textcolor{green}{\textbf{142407413}};142440190;\textcolor{green}{\textbf{142433446}};5.80e-05;\textcolor{red}{\textit{7.10e-05}};2.15e-04;\textcolor{green}{\textbf{1.74e-04}};4.04e-04;\textcolor{green}{\textbf{3.57e-04}};1112645;\textcolor{red}{\textit{1337066}};1149094;\textcolor{red}{\textit{1392517}};1230683;\textcolor{red}{\textit{1464697}};+1.27e-03;-4.07e-03;-4.73e-03;+22.41;-19.07;-11.63;+20.17;+21.18;+19.01
DIMACS\_E100k.1;POPMUSIC;5;Error;10.00;10.00;225692504;\textcolor{green}{\textbf{225671959}};225734644;\textcolor{green}{\textbf{225680495}};225759877;\textcolor{green}{\textbf{225688956}};nan;0.0;nan;3.80e-05;nan;7.50e-05;748667;\textcolor{green}{\textbf{472123}};803160;\textcolor{green}{\textbf{532180}};965099;\textcolor{green}{\textbf{614401}};-9.10e-03;-2.40e-02;-3.14e-02;$\pm$ 0;$\pm$ 0;$\pm$ 0;-36.94;-33.74;-36.34
DIMACS\_E100k.0;POPMUSIC;5;Error;10.00;10.00;225852535;\textcolor{green}{\textbf{225806177}};225883814;\textcolor{green}{\textbf{225812934}};225920212;\textcolor{green}{\textbf{225820383}};nan;0.0;nan;3e-05;nan;6.30e-05;771576;\textcolor{green}{\textbf{526575}};815162;\textcolor{green}{\textbf{575585}};899004;\textcolor{green}{\textbf{741841}};-2.05e-02;-3.14e-02;-4.42e-02;$\pm$ 0;$\pm$ 0;$\pm$ 0;-31.75;-29.39;-17.48
DIMACS\_C316k.0;POPMUSIC;5;nan;5.000;6.000;186899030;\textcolor{red}{\textit{186935218}};186910144;\textcolor{red}{\textit{186953074}};186941216;\textcolor{red}{\textit{186987561}};nan;nan;nan;nan;nan;nan;1944189;\textcolor{green}{\textbf{1657905}};2112022;\textcolor{green}{\textbf{1791947}};2281218;\textcolor{green}{\textbf{2023697}};+1.94e-02;+2.30e-02;+2.48e-02;$\pm$ 0;$\pm$ 0;$\pm$ 0;-14.73;-15.15;-11.29
vlsi\_sra104815;POPMUSIC;5;251342;10.00;9.000;251437;\textcolor{green}{\textbf{251383}};251462;\textcolor{green}{\textbf{251406}};251496;\textcolor{green}{\textbf{251438}};3.78e-04;\textcolor{green}{\textbf{1.63e-04}};4.77e-04;\textcolor{green}{\textbf{2.56e-04}};6.13e-04;\textcolor{green}{\textbf{3.82e-04}};760748;\textcolor{red}{\textit{811546}};834368;\textcolor{red}{\textit{893250}};943314;\textcolor{red}{\textit{986642}};-2.15e-02;-2.23e-02;-2.31e-02;-56.88;-46.33;-37.68;+6.68;+7.06;+4.59
\end{filecontents*}
\csvreader[head to column names, separator=semicolon]{Latex_Running_Comparison_Large_POPMUSIC.csv}{}
    {\Problem & \CostMinLKH & \CostMinOur & \TimeAvgLKH & \TimeAvgOur & \TimeAvgChange \\}
    \\[-\normalbaselineskip]\hline
\end{tabular}
\label{fig:CR_large_running_instances_POPMUSIC}
\end{table}
\begin{examplex}
Consider the DIMACS instance E100k.1 in \cref{fig:CR_large_running_instances_POPMUSIC,%
fig:CR_large_running_instances_alpha}. 
With a POPMUSIC generated candidate set the average time per run of our variant is approximately 6.2 days. 
Since the original LKH needs approximately 9.3 days we reduce the average time per run by one third. 
With an ALPHA candidate set the computation time is reduced even more. 
Our variant needs approximately 5.3 days per run, the original LKH needs 9.8 days. 
This represents a time saving of 46.07\% using our variant instead of the original one. 
Additionally, we obtain a better minimal cost with our variant for both candidate set types. 
\end{examplex}

As shown in \cref{fig:ratio_large_instances,%
fig:CR_large_running_instances_POPMUSIC}
%our variant using 
the \tilted relaxation decreases the average time per run with POPMUSIC candidate sets by 6.4\% on average compared to the original LKH. 
For a few instances our variants perform worse.
For the TSP-library instance pla85900 the average time per run is increased by 21.18\%.  
Also, the minimal cost for this instance is increased by 0.001\%.

\begin{table}
\caption{Direct comparison of our variant using the \tilted relaxation \ref{C2**:tiltedgain} and the original LKH with five ALPHA candidates in minimal cost and average time per run for all \textbf{large instances that did not finish in one month of computation time}.
}
\begin{tabular}{lrrrrr}%
\hline
\multicolumn{1}{c}{Problem}  & \multicolumn{2}{c}{CostMin} & \multicolumn{2}{c}{TimeAvg in s} & \multicolumn{1}{c}{TimeAvg} \\ 
& LKH 3.0.8 & with \ref{C2**:tiltedgain} & LKH 3.0.8 & with \ref{C2**:tiltedgain} & \multicolumn{1}{c}{ratio in \%}\\ \hline% specify table head
    \begin{filecontents*}{Latex_Running_Comparison_Large_ALPHA.csv}
Problem;CandidateModel;Candidates;Optimum;RunsLKH;RunsOur;CostMinLKH;CostMinOur;CostAvgLKH;CostAvgOur;CostMaxLKH;CostMaxOur;GapMinLKH;GapMinOur;GapAvgLKH;GapAvgOur;GapMaxLKH;GapMaxOur;TimeMinLKH;TimeMinOur;TimeAvgLKH;TimeAvgOur;TimeMaxLKH;TimeMaxOur;CostMinChange;CostAvgChange;CostMaxChange;GapMinChange;GapAvgChange;GapMaxChange;TimeMinChange;TimeAvgChange;TimeMaxChange
vlsi\_sra104815;ALPHA;5;251342;10.00;10.00;251389;\textcolor{green}{\textbf{251373}};251444;\textcolor{green}{\textbf{251402}};251508;\textcolor{green}{\textbf{251465}};1.87e-04;\textcolor{green}{\textbf{1.23e-04}};4.04e-04;\textcolor{green}{\textbf{2.38e-04}};6.60e-04;\textcolor{green}{\textbf{4.89e-04}};622045;\textcolor{red}{\textit{624171}};741206;\textcolor{red}{\textit{741771}};904440;\textcolor{red}{\textit{904865}};-6.36e-03;-1.67e-02;-1.71e-02;-34.22;-41.09;-25.91;+0.34;+0.08;+0.05
vlsi\_ara238025;ALPHA;5;578761;2.000;3.000;578970;\textcolor{green}{\textbf{578914}};579091;\textcolor{green}{\textbf{578978}};579212;\textcolor{green}{\textbf{579050}};3.61e-04;\textcolor{green}{\textbf{2.64e-04}};5.70e-04;\textcolor{green}{\textbf{3.74e-04}};7.79e-04;\textcolor{green}{\textbf{4.99e-04}};3633175;\textcolor{green}{\textbf{3286557}};4260252;\textcolor{green}{\textbf{3749923}};4887328;\textcolor{green}{\textbf{4380766}};-9.67e-03;-1.95e-02;-2.80e-02;-26.87;-34.39;-35.94;-9.54;-11.98;-10.36
TSP\_pla85900;ALPHA;5;142382641;10.00;10.00;142420045;\textcolor{green}{\textbf{142405913}};142437076;\textcolor{green}{\textbf{142416622}};142448988;\textcolor{green}{\textbf{142433435}};2.63e-04;\textcolor{green}{\textbf{1.63e-04}};3.82e-04;\textcolor{green}{\textbf{2.38e-04}};4.66e-04;\textcolor{green}{\textbf{3.57e-04}};358090;\textcolor{red}{\textit{409078}};385645;\textcolor{red}{\textit{468678}};462521;\textcolor{red}{\textit{571563}};-9.92e-03;-1.44e-02;-1.09e-02;-38.02;-37.70;-23.39;+14.24;+21.53;+23.58
DIMACS\_E100k.0;ALPHA;5;Error;10.00;10.00;225819355;\textcolor{green}{\textbf{225804698}};225855098;\textcolor{green}{\textbf{225812262}};225882144;\textcolor{green}{\textbf{225823927}};nan;0.0;nan;3.30e-05;nan;8.50e-05;741552;\textcolor{green}{\textbf{406378}};803876;\textcolor{green}{\textbf{442230}};900291;\textcolor{green}{\textbf{492377}};-6.49e-03;-1.90e-02;-2.58e-02;$\pm$ 0;$\pm$ 0;$\pm$ 0;-45.20;-44.99;-45.31
DIMACS\_E100k.1;ALPHA;5;Error;10.00;10.00;225725610;\textcolor{green}{\textbf{225670787}};225746192;\textcolor{green}{\textbf{225680825}};225774528;\textcolor{green}{\textbf{225690780}};nan;0.0;nan;4.40e-05;nan;8.90e-05;797771;\textcolor{green}{\textbf{418123}};843327;\textcolor{green}{\textbf{454843}};954691;\textcolor{green}{\textbf{517754}};-2.43e-02;-2.90e-02;-3.71e-02;$\pm$ 0;$\pm$ 0;$\pm$ 0;-47.59;-46.07;-45.77
\end{filecontents*}
\csvreader[head to column names, separator=semicolon]{Latex_Running_Comparison_Large_ALPHA.csv}{}
    {\Problem  & \CostMinLKH & \CostMinOur &\TimeAvgLKH & \TimeAvgOur &\TimeAvgChange \\}
    \\[-\normalbaselineskip]\hline
\end{tabular}
\label{fig:CR_large_running_instances_alpha}
\end{table}

With ALPHA candidate sets 
averaging over the five instances shown in \cref{fig:CR_large_running_instances_alpha} the \tilted relaxation reduces the time per run by 16.3\%. 
For the DIMACS instances, i.e. E100k.0 and E100k.1, this reduction is approximately 45\%.
For the TSP\_pla85900 instance the \tilted relaxation increases the average time per run also with ALPHA candidate sets. 
The time save per run with the \tilted relaxation can be used to find slightly better solutions on average. 
Note that at least for those instances where the optimum is known the solutions are already very close to optimality. 
Thus, more than slight improvements are not possible. 

\begin{table}[h]
\caption{Direct comparison of our variant using the \tilted relaxation \ref{C2**:tiltedgain} and the original LKH with five POPMUSIC candidates in minimal cost and average cost per run for all \textbf{large instances that did not finish in one month of computation time} with a \textbf{time limit} of three days.}
\begin{tabular}{lrrrrrr}%
\hline
Problem  & \multicolumn{2}{c}{CostMin} & %\multicolumn{2}{c}{CostAvg} & 
\multicolumn{2}{c}{CostAvg} & \multicolumn{1}{c}{CostMin} \\ 
& LKH 3.0.8 & with \ref{C2**:tiltedgain} & LKH 3.0.8 & with \ref{C2**:tiltedgain} & \multicolumn{1}{c}{ratio in \%}\\ \hline% specify table head
    \begin{filecontents*}{Latex_Running_Comparison_Large_stop_POPMUSIC.csv}
Problem;CandidateModel;Candidates;Optimum;RunsLKH;RunsOur;CostMinLKH;CostMinOur;CostAvgLKH;CostAvgOur;CostMaxLKH;CostMaxOur;GapMinLKH;GapMinOur;GapAvgLKH;GapAvgOur;GapMaxLKH;GapMaxOur;TimeMinLKH;TimeMinOur;TimeAvgLKH;TimeAvgOur;TimeMaxLKH;TimeMaxOur;CostMinChange;CostAvgChange;CostMaxChange;GapMinChange;GapAvgChange;GapMaxChange;TimeMinChange;TimeAvgChange;TimeMaxChange
TSP\_pla85900;POPMUSIC;5;15343;10.00;10.00;142395857;\textcolor{red}{\textit{142401344}};142418416;\textcolor{green}{\textbf{142413718}};142469731;\textcolor{green}{\textbf{142430510}};9.30e-05;\textcolor{red}{\textit{1.31e-04}};2.51e-04;\textcolor{green}{\textbf{2.18e-04}};6.12e-04;\textcolor{green}{\textbf{3.36e-04}};259201;\textcolor{red}{\textit{259203}};259207;\textcolor{green}{\textbf{259206}};259210;\textcolor{red}{\textit{259214}};+3.85e-03;-3.30e-03;-2.75e-02;+40.86;-13.15;-45.10;+0.00;-0.00;+0.00
DIMACS\_E100k.0;POPMUSIC;5;43167;10.00;10.00;225871989;\textcolor{green}{\textbf{225812844}};225897966;\textcolor{green}{\textbf{225818738}};225917533;\textcolor{green}{\textbf{225827680}};0.0;0.0;1.15e-04;\textcolor{green}{\textbf{2.60e-05}};2.02e-04;\textcolor{green}{\textbf{6.55e-05}};259201;259201;259204;\textcolor{green}{\textbf{259203}};259209;\textcolor{green}{\textbf{259206}};-2.62e-02;-3.51e-02;-3.98e-02;$\pm$ 0;-77.39;-67.57;$\pm$ 0;-0.00;-0.00
DIMACS\_E100k.1;POPMUSIC;5;43149;10.00;10.00;225743656;\textcolor{green}{\textbf{225681751}};225774364;\textcolor{green}{\textbf{225685683}};225809192;\textcolor{green}{\textbf{225693252}};0.0;0.0;1.36e-04;\textcolor{green}{\textbf{1.75e-05}};2.90e-04;\textcolor{green}{\textbf{5.10e-05}};259202;\textcolor{green}{\textbf{259200}};259206;\textcolor{green}{\textbf{259203}};259208;\textcolor{green}{\textbf{259206}};-2.74e-02;-3.93e-02;-5.13e-02;$\pm$ 0;-87.13;-82.41;-0.00;-0.00;-0.00
vlsi\_sra104815;POPMUSIC;5;26296;10.00;10.00;251466;\textcolor{green}{\textbf{251413}};251502;\textcolor{green}{\textbf{251454}};251555;\textcolor{green}{\textbf{251506}};4.91e-04;\textcolor{green}{\textbf{2.82e-04}};6.34e-04;\textcolor{green}{\textbf{4.46e-04}};8.47e-04;\textcolor{green}{\textbf{6.51e-04}};259201;\textcolor{red}{\textit{259202}};259203;\textcolor{red}{\textit{259206}};259206;\textcolor{red}{\textit{259212}};-2.11e-02;-1.91e-02;-1.95e-02;-42.57;-29.65;-23.14;+0.00;+0.00;+0.00
vlsi\_ara238025;POPMUSIC;5;8887;10.00;10.00;579233;\textcolor{green}{\textbf{579146}};579368;\textcolor{green}{\textbf{579237}};579459;\textcolor{green}{\textbf{579322}};8.16e-04;\textcolor{green}{\textbf{6.66e-04}};1.05e-03;\textcolor{green}{\textbf{8.22e-04}};1.21e-03;\textcolor{green}{\textbf{9.70e-04}};259203;\textcolor{red}{\textit{259206}};259215;\textcolor{green}{\textbf{259214}};259229;259229;-1.50e-02;-2.26e-02;-2.36e-02;-18.38;-21.71;-19.83;+0.00;-0.00;$\pm$ 0
DIMACS\_E316k.0;POPMUSIC;5;6944;10.00;10.00;401650710;\textcolor{green}{\textbf{401518962}};401673366;\textcolor{green}{\textbf{401535147}};401716850;\textcolor{green}{\textbf{401553991}};0.0;0.0;5.60e-05;\textcolor{green}{\textbf{4.05e-05}};1.64e-04;\textcolor{green}{\textbf{8.75e-05}};259207;\textcolor{red}{\textit{259210}};259222;\textcolor{green}{\textbf{259217}};259236;\textcolor{green}{\textbf{259227}};-3.28e-02;-3.44e-02;-4.05e-02;$\pm$ 0;-27.68;-46.65;+0.00;-0.00;-0.00
DIMACS\_C316k.0;POPMUSIC;5;41374;10.00;10.00;186956459;\textcolor{green}{\textbf{186931858}};186977382;\textcolor{green}{\textbf{186965247}};187015850;\textcolor{green}{\textbf{187013758}};0.0;0.0;1.12e-04;\textcolor{red}{\textit{1.79e-04}};3.18e-04;\textcolor{red}{\textit{4.38e-04}};259201;\textcolor{red}{\textit{259202}};259204;\textcolor{green}{\textbf{259203}};259208;\textcolor{green}{\textbf{259205}};-1.32e-02;-6.49e-03;-1.12e-03;$\pm$ 0;+59.82;+37.74;+0.00;-0.00;-0.00
vlsi\_lra498378;POPMUSIC;5;7348;10.00;10.00;2170006;\textcolor{green}{\textbf{2169464}};2170345;\textcolor{green}{\textbf{2169824}};2170938;\textcolor{green}{\textbf{2170226}};9.07e-04;\textcolor{green}{\textbf{6.57e-04}};1.06e-03;\textcolor{green}{\textbf{8.24e-04}};1.34e-03;\textcolor{green}{\textbf{1.01e-03}};259211;259211;259217;\textcolor{red}{\textit{259223}};259225;\textcolor{red}{\textit{259238}};-2.50e-02;-2.40e-02;-3.28e-02;-27.56;-22.26;-24.63;$\pm$ 0;+0.00;+0.01
vlsi\_lrb744710;POPMUSIC;5;2032;10.00;10.00;1613738;\textcolor{green}{\textbf{1613266}};1613841;\textcolor{green}{\textbf{1613467}};1614070;\textcolor{green}{\textbf{1613736}};1.56e-03;\textcolor{green}{\textbf{1.26e-03}};1.62e-03;\textcolor{green}{\textbf{1.39e-03}};1.76e-03;\textcolor{green}{\textbf{1.55e-03}};259219;259219;259256;\textcolor{red}{\textit{259278}};259286;\textcolor{red}{\textit{259322}};-2.92e-02;-2.32e-02;-2.07e-02;-19.23;-14.20;-11.93;$\pm$ 0;+0.01;+0.01
\hline Average;POPMUSIC;5.000;21616;10.00;10.00;131914790;\textcolor{green}{\textbf{131884450}};131928505;\textcolor{green}{\textbf{131892502}};131949464;\textcolor{green}{\textbf{131903776}};4.29e-04;\textcolor{green}{\textbf{3.33e-04}};5.60e-04;\textcolor{green}{\textbf{4.40e-04}};7.49e-04;\textcolor{green}{\textbf{5.73e-04}};259205;\textcolor{red}{\textit{259206}};259215;\textcolor{red}{\textit{259217}};259224;\textcolor{red}{\textit{259229}};-2.10e-02;-2.30e-02;-2.90e-02;nan;-26.01;-31.45;3.50e-04;8.97e-04;1.74e-03
Median;POPMUSIC;5.000;15343;10.00;10.00;142395857;\textcolor{red}{\textit{142401344}};142418416;\textcolor{green}{\textbf{142413718}};142469731;\textcolor{green}{\textbf{142430510}};9.30e-05;\textcolor{red}{\textit{1.31e-04}};2.51e-04;\textcolor{green}{\textbf{2.18e-04}};6.12e-04;\textcolor{green}{\textbf{4.38e-04}};259202;\textcolor{red}{\textit{259203}};259207;\textcolor{green}{\textbf{259206}};259210;\textcolor{red}{\textit{259214}};-2.50e-02;-2.30e-02;-2.80e-02;nan;-22.60;-24.57;2.57e-04;-2.76e-04;-7.91e-05
95th Quantil;POPMUSIC;5.000;43160;10.00;10.00;331339221;\textcolor{green}{\textbf{331236515}};331363206;\textcolor{green}{\textbf{331248583}};331397123;\textcolor{green}{\textbf{331263466}};1.30e-03;\textcolor{green}{\textbf{1.02e-03}};1.40e-03;\textcolor{green}{\textbf{1.16e-03}};1.59e-03;\textcolor{green}{\textbf{1.34e-03}};259216;259216;259242;\textcolor{red}{\textit{259256}};259266;\textcolor{red}{\textit{259288}};-2.95e-03;-4.58e-03;-8.54e-03;nan;30.45;18.08;1.18e-03;6.03e-03;1.00e-02
\end{filecontents*}
\csvreader[head to column names, separator=semicolon,
    late after line=\ifnumequal{\thecsvrow}{10}{\\\hline}{\\},
    late after last line=\\\hline,]{Latex_Running_Comparison_Large_stop_POPMUSIC.csv}{}
    {\Problem & \CostMinLKH & \CostMinOur & \CostAvgLKH & \CostAvgOur & \CostMinChange}
\end{tabular}
\label{fig:CR_large_stop_instances_POPMUSIC}
\end{table}
\begin{table}[h]
\caption{Direct comparison of our variant using the \tilted relaxation \ref{C2**:tiltedgain} and the original LKH with five ALPHA candidates in minimal cost and average cost per run for all \textbf{large instances that did not finish in one month of computation time} with a \textbf{time limit} of three days.}
\begin{tabular}{lrrrrr}%
\hline
Problem  & \multicolumn{2}{c}{CostMin} & \multicolumn{2}{c}{CostAvg} & \multicolumn{1}{c}{CostMin} \\ 
& LKH 3.0.8 & with \ref{C2**:tiltedgain} & LKH 3.0.8 & with \ref{C2**:tiltedgain} & \multicolumn{1}{c}{ratio in \%}\\ \hline% specify table head
    \begin{filecontents*}{Latex_Running_Comparison_Large_stop_ALPHA.csv}
Problem;CandidateModel;Candidates;Optimum;RunsLKH;RunsOur;CostMinLKH;CostMinOur;CostAvgLKH;CostAvgOur;CostMaxLKH;CostMaxOur;GapMinLKH;GapMinOur;GapAvgLKH;GapAvgOur;GapMaxLKH;GapMaxOur;TimeMinLKH;TimeMinOur;TimeAvgLKH;TimeAvgOur;TimeMaxLKH;TimeMaxOur;CostMinChange;CostAvgChange;CostMaxChange;GapMinChange;GapAvgChange;GapMaxChange;TimeMinChange;TimeAvgChange;TimeMaxChange
TSP\_pla85900;ALPHA;5;46976;10.00;10.00;142418874;\textcolor{red}{\textit{142434986}};142432185;\textcolor{red}{\textit{142452978}};142446040;\textcolor{red}{\textit{142493352}};2.55e-04;\textcolor{red}{\textit{3.67e-04}};3.48e-04;\textcolor{red}{\textit{4.94e-04}};4.45e-04;\textcolor{red}{\textit{7.77e-04}};259201;259201;259202;259202;259204;259204;+1.13e-02;+1.46e-02;+3.32e-02;+43.92;+41.95;+74.61;$\pm$ 0;$\pm$ 0;$\pm$ 0
DIMACS\_E100k.0;ALPHA;5;53237;10.00;10.00;225854660;\textcolor{green}{\textbf{225803922}};225877921;\textcolor{green}{\textbf{225811661}};225909128;\textcolor{green}{\textbf{225821364}};0.0;0.0;1.03e-04;\textcolor{green}{\textbf{3.45e-05}};2.41e-04;\textcolor{green}{\textbf{7.70e-05}};259202;\textcolor{red}{\textit{259203}};259203;\textcolor{red}{\textit{259205}};259205;\textcolor{red}{\textit{259208}};-2.25e-02;-2.93e-02;-3.88e-02;$\pm$ 0;-66.50;-68.05;+0.00;+0.00;+0.00
DIMACS\_E100k.1;ALPHA;5;52678;10.00;10.00;225717093;\textcolor{green}{\textbf{225680993}};225741713;\textcolor{green}{\textbf{225685384}};225761658;\textcolor{green}{\textbf{225690892}};0.0;0.0;1.09e-04;\textcolor{green}{\textbf{1.95e-05}};1.98e-04;\textcolor{green}{\textbf{4.35e-05}};259201;259201;259203;259203;259207;\textcolor{green}{\textbf{259206}};-1.60e-02;-2.50e-02;-3.13e-02;$\pm$ 0;-82.11;-78.03;$\pm$ 0;$\pm$ 0;-0.00
vlsi\_sra104815;ALPHA;5;32162;10.00;10.00;251433;\textcolor{green}{\textbf{251397}};251466;\textcolor{green}{\textbf{251417}};251504;\textcolor{green}{\textbf{251438}};3.62e-04;\textcolor{green}{\textbf{2.19e-04}};4.94e-04;\textcolor{green}{\textbf{2.98e-04}};6.47e-04;\textcolor{green}{\textbf{3.84e-04}};259200;\textcolor{red}{\textit{259202}};259202;\textcolor{red}{\textit{259206}};259206;\textcolor{red}{\textit{259213}};-1.43e-02;-1.95e-02;-2.62e-02;-39.50;-39.68;-40.65;+0.00;+0.00;+0.00
vlsi\_ara238025;ALPHA;5;12640;10.00;10.00;579186;\textcolor{green}{\textbf{579073}};579433;\textcolor{green}{\textbf{579172}};579634;\textcolor{green}{\textbf{579327}};7.33e-04;\textcolor{green}{\textbf{5.39e-04}};1.16e-03;\textcolor{green}{\textbf{7.11e-04}};1.51e-03;\textcolor{green}{\textbf{9.78e-04}};259203;\textcolor{green}{\textbf{259202}};259209;\textcolor{green}{\textbf{259208}};259217;\textcolor{red}{\textit{259220}};-1.95e-02;-4.50e-02;-5.30e-02;-26.47;-38.71;-35.23;-0.00;-0.00;+0.00
DIMACS\_C316k.0;ALPHA;5;26906;10.00;10.00;4146826698;\textcolor{red}{\textit{5467724892}};5744491051;\textcolor{red}{\textit{6097665242}};7827580207;\textcolor{green}{\textbf{7193228842}};0.0;0.0;0.390;\textcolor{green}{\textbf{0.115}};0.896;\textcolor{green}{\textbf{0.315}};259201;259201;259215;\textcolor{green}{\textbf{259205}};259246;\textcolor{green}{\textbf{259211}};+31.85;+6.15;-8.10;$\pm$ 0;-70.51;-64.84;$\pm$ 0;-0.00;-0.01
DIMACS\_E316k.0;ALPHA;5;7502;10.00;10.00;401608912;\textcolor{green}{\textbf{401509132}};401638980;\textcolor{green}{\textbf{401518121}};401671244;\textcolor{green}{\textbf{401534600}};0.0;0.0;7.45e-05;\textcolor{green}{\textbf{2.20e-05}};1.55e-04;\textcolor{green}{\textbf{6.35e-05}};259204;\textcolor{red}{\textit{259209}};259223;\textcolor{green}{\textbf{259218}};259239;\textcolor{green}{\textbf{259227}};-2.48e-02;-3.01e-02;-3.40e-02;$\pm$ 0;-70.47;-59.03;+0.00;-0.00;-0.00
vlsi\_lra498378;ALPHA;5;3954;10.00;10.00;2694856;\textcolor{red}{\textit{4067073}};3460386;\textcolor{red}{\textit{5288952}};5321168;\textcolor{red}{\textit{7085603}};0.243;\textcolor{red}{\textit{0.876}};0.596;\textcolor{red}{\textit{1.440}};1.454;\textcolor{red}{\textit{2.268}};259214;\textcolor{green}{\textbf{259210}};259254;\textcolor{green}{\textbf{259239}};259322;\textcolor{green}{\textbf{259276}};+50.92;+52.84;+33.16;+260.49;+141.61;+55.98;-0.00;-0.01;-0.02
vlsi\_lrb744710;ALPHA;5;2434;10.00;10.00;1613576;\textcolor{green}{\textbf{1613436}};1613753;\textcolor{green}{\textbf{1613495}};1613983;\textcolor{green}{\textbf{1613552}};1.45e-03;\textcolor{green}{\textbf{1.37e-03}};1.56e-03;\textcolor{green}{\textbf{1.40e-03}};1.71e-03;\textcolor{green}{\textbf{1.44e-03}};259215;\textcolor{red}{\textit{259227}};259248;\textcolor{red}{\textit{259257}};259282;\textcolor{red}{\textit{259285}};-8.68e-03;-1.60e-02;-2.67e-02;-5.52;-10.26;-15.79;+0.00;+0.00;+0.00
\hline Average;ALPHA;5.000;26499;10.00;10.00;571951698;\textcolor{red}{\textit{718851656}};749565210;\textcolor{red}{\textit{788985158}};981237174;\textcolor{green}{\textbf{910922108}};0.027;\textcolor{red}{\textit{0.098}};0.110;\textcolor{red}{\textit{0.173}};0.262;\textcolor{red}{\textit{0.287}};259205;\textcolor{red}{\textit{259206}};259218;\textcolor{green}{\textbf{259216}};259236;\textcolor{green}{\textbf{259228}};9.19;6.54;2.76;nan;-21.65;-25.65;6.22e-04;-6.78e-04;-3.32e-03
Median;ALPHA;5.000;26906;10.00;10.00;142418874;\textcolor{red}{\textit{142434986}};142432185;\textcolor{red}{\textit{142452978}};142446040;\textcolor{red}{\textit{142493352}};2.55e-04;\textcolor{green}{\textbf{2.19e-04}};4.94e-04;4.94e-04;6.47e-04;\textcolor{red}{\textit{7.77e-04}};259202;259202;259209;\textcolor{green}{\textbf{259206}};259217;\textcolor{green}{\textbf{259213}};-1.40e-02;-2.00e-02;-3.10e-02;nan;-39.68;-40.60;1.77e-04;-2.53e-04;5.79e-06
95th Quantil;ALPHA;5.000;53013;10.00;10.00;2648739583;\textcolor{red}{\textit{3441238587}};3607350222;\textcolor{red}{\textit{3819206393}};4857216622;\textcolor{green}{\textbf{4476551145}};0.146;\textcolor{red}{\textit{0.526}};0.514;\textcolor{red}{\textit{0.910}};1.231;\textcolor{red}{\textit{1.487}};259215;\textcolor{red}{\textit{259220}};259251;\textcolor{green}{\textbf{259250}};259306;\textcolor{green}{\textbf{259281}};43.29;34.16;19.91;nan;101.70;67.10;3.51e-03;2.76e-03;2.22e-03
\end{filecontents*}
\csvreader[head to column names, 
    separator=semicolon, 
    late after line=\ifnumequal{\thecsvrow}{10}{\\\hline}{\\},
    late after last line=\\\hline,]{Latex_Running_Comparison_Large_stop_ALPHA.csv}{}
    {\Problem & \CostMinLKH & \CostMinOur &\CostAvgLKH & \CostAvgOur & \CostMinChange}
\end{tabular}
\label{fig:CR_large_stop_instances_alpha}
\end{table}
\cref{fig:CR_large_stop_instances_POPMUSIC,%
	fig:CR_large_stop_instances_alpha} show the results when we restrict the \tilted relaxation to the same number of runs and the same ten initial tours as the original LKH. 
Concerning the changes in the minimal cost with the POPMUSIC candidate sets, confer \cref{fig:CR_large_stop_instances_POPMUSIC}, we achieved slightly lower minimal cost for eight of nine instances and worse for one instance.
With ALPHA candidate sets confer \cref{fig:CR_large_stop_instances_alpha} the results are similar except for two instances, where neither our variants nor the original LKH finished one run. 
As already recommended by Helsgaun our results confirm using POPMUSIC candidate sets instead of ALPHA candidate sets. 

To summarize, with a growing number of nodes, we substantially decrease the computation time using \ref{C2**:tiltedgain} instead of \ref{C2:gain}. 
For the small instances our variant is half a second slower per run. 
However, it achieves an average time reduction of about 90 seconds or four percent for the medium instances. 
For the large instances we reduce the average time per run by 13.1\% and 17.8\%.
This holds in particular for the DIMACS instances, i.e. DIMACS\_E100k.0 and DIMACS\_E100k.1, for which we achieve a time reduction of about 45 \% with ALPHA and about 30\% with POPMUSIC candidate sets. 
\subsection{Comparing \homogeneous and \tilted relaxation}
\label{subsection:algo_engineering_results}
In \cref{fig:CR_large_instances_fully_relaxed,% 
fig:CR_medium_instances_fully_relaxed,% 
fig:CR_small_instances_fully_relaxed} we compare our two variants, i.e. the \homogeneous relaxation \ref{C2*:newgain} and the \tilted relaxation \ref{C2**:tiltedgain} of the positive gain criterion.
For all large instances the minimal cost and the average cost are the same for both variants, see \cref{fig:CR_large_instances_fully_relaxed}.
Note that we used the same random seed for both variants.
However, the variant using \ref{C2**:tiltedgain} is 0.8 hours faster per run than the variant using \ref{C2*:newgain} considering a POPMUSIC candidate set. 
%\birte{Saetze anpassen! Richtige \% nehmen}
This means that with \ref{C2**:tiltedgain} the running time decreases by 13.6\% compared to LKH, whereas with \ref{C2*:newgain} the running time decreases by 10.4\%.
With an ALPHA candidate set, we likewise report that our variant using \ref{C2**:tiltedgain} is 0.8 hours faster per run than our variant using \ref{C2*:newgain}. 
This means that with \ref{C2**:tiltedgain} the running time decreases by 20.8\% compared to LKH, whereas with \ref{C2*:newgain} the running time decreases by 17\%.

In \cref{fig:CR_small_instances_fully_relaxed,%
	fig:CR_medium_instances_fully_relaxed} we depict the averaged results over all small and medium instances, respectively.
The solution quality stays exactly the same for both of our relaxations while the running time is slightly worse for the variant using the \tilted relaxation \ref{C2**:tiltedgain} compared to the variant using the \homogeneous relaxation \ref{C2*:newgain}.
\begin{table}[h]
\caption{Comparison of our \homogeneous relaxation \ref{C2*:newgain} variant and our \tilted \ref{C2**:tiltedgain} relaxation variant in minimal cost, average cost and average time per run on \textbf{large instances} with different candidate set types. Green and bold indicate a better performance of our \tilted relaxation, red and italic a worse performance.}
\begin{tabular}{lrrrrr}%
\hline
     Version & Candidates & \# Ex. & CostMin & CostAvg & TimeAvg
     \\ \hline
    \begin{filecontents*}{Homogeneous_Latex_Summary_Large.csv}
Version;CandidateModel;Candidates;Examples;Runs;Successes;TrialsMin;TrialsAvg;TrialsMax;Optimum;CostMin;CostAvg;CostMax;GapMin;GapAvg;GapMax;TimeMin;TimeAvg;TimeMax
\Homogeneous relaxation;ALPHA;5;20;10.00;0.0;46688;46688;46688;Error;32856187;32880367;32930336;8.02e-03;4.50e-04;1.18e-03;55546;64059;85135
\Tilted relaxation;ALPHA;5;20;10.00;0.0;46688;46688;46688;Error;32856187;32880367;32930336;8.02e-03;4.50e-04;1.18e-03;\textcolor{green}{\textbf{51723}};\textcolor{green}{\textbf{61126}};\textcolor{green}{\textbf{78546}}
\Homogeneous relaxation;POPMUSIC;5;20;10.00;0.0;46688;46688;46688;Error;32824384;32832052;32841860;6.38e-03;2.24e-04;4.35e-04;71868;81789;91839
\Tilted relaxation;POPMUSIC;5;20;10.00;0.0;46688;46688;46688;Error;32824384;32832052;32841860;6.38e-03;2.24e-04;4.35e-04;\textcolor{green}{\textbf{69651}};\textcolor{green}{\textbf{78864}};\textcolor{green}{\textbf{88445}}
\end{filecontents*}
\csvreader[head to column names, separator=semicolon]{Homogeneous_Latex_Summary_Large.csv}{}
    { \Version & \Candidates\ \CandidateModel & \Examples & \CostMin & \CostAvg & \TimeAvg \;s\\}
    \\[-\normalbaselineskip]\hline
\end{tabular}
\label{fig:CR_large_instances_fully_relaxed}
\end{table}
\begin{table}
\caption{Comparison of our \homogeneous relaxation \ref{C2*:newgain} variant and our \tilted  relaxation \ref{C2**:tiltedgain} variant in minimal gap, average gap and average time per run on \textbf{medium instances} with different candidate set types.}
\begin{tabular}{lrrrrr}%
\hline
     Version & Candidates & \# Ex. & GapMin & GapAvg & TimeAvg % specify table head
     \\ 
     &&&  in \% &  in \% & in s\\ \hline
    \begin{filecontents*}{Homogeneous_Latex_Summary_Medium.csv}
Version;CandidateModel;Candidates;Examples;Runs;Successes;TrialsMin;TrialsAvg;TrialsMax;Optimum;CostMin;CostAvg;CostMax;GapMin;GapAvg;GapMax;TimeMin;TimeAvg;TimeMax
\Homogeneous relaxation;ALPHA;5;140;10.00;3.857;6957;6957;6957;Error;4863897;4864538;4865433;5.88e-03;2.60e-04;5.23e-04;1622;1962;2403
\Tilted relaxation;ALPHA;5;140;10.00;3.857;6957;6957;6957;Error;4863897;4864538;4865433;5.88e-03;2.60e-04;5.23e-04;\textcolor{red}{\textit{1693}};\textcolor{red}{\textit{2043}};\textcolor{red}{\textit{2499}}
\Homogeneous relaxation;POPMUSIC;5;140;10.00;4.557;6957;6957;6957;Error;4864873;4865371;4866057;0.053;7.56e-04;1.02e-03;2006;2386;2859
\Tilted relaxation;POPMUSIC;5;140;10.00;4.557;6957;6957;6957;Error;4864873;4865371;4866057;0.053;7.56e-04;1.02e-03;\textcolor{red}{\textit{2053}};\textcolor{red}{\textit{2454}};\textcolor{red}{\textit{2948}}
\Homogeneous relaxation;POPMUSIC;100;4;100.0;0.0;3501;3501;3501;Error;65402934;65404116;65404771;0.0;2.85e-05;4.15e-05;6.308;11.28;22.65
\Tilted relaxation;POPMUSIC;100;4;100.0;0.0;3501;3501;3501;Error;\textcolor{green}{\textbf{65402841}};\textcolor{red}{\textit{65405417}};\textcolor{red}{\textit{65405788}};0.0;\textcolor{red}{\textit{3.90e-05}};\textcolor{red}{\textit{4.45e-05}};\textcolor{red}{\textit{6.543}};\textcolor{red}{\textit{13.80}};\textcolor{red}{\textit{25.31}}
\end{filecontents*}
\csvreader[head to column names, separator=semicolon]{Homogeneous_Latex_Summary_Medium.csv}{}
    { \Version & \Candidates\ \CandidateModel & \Examples & \GapMin & \GapAvg & \TimeAvg\\}
    \\[-\normalbaselineskip]\hline
\end{tabular}
\label{fig:CR_medium_instances_fully_relaxed}
\end{table}
\begin{table}
\caption{Comparison of our \homogeneous relaxation \ref{C2*:newgain} variant and our \tilted relaxation \ref{C2**:tiltedgain} variant in minimal gap, average gap and average time per run on \textbf{small instances} with different candidate set types.}
\begin{tabular}{lrrrrr}%
\hline
     Version & Candidates & \# Ex. & GapMin & GapAvg & TimeAvg % specify table head
     \\ 
     &&&  in \% &  in \% & in s\\ \hline
    \begin{filecontents*}{Homogeneous_Latex_Summary_Small.csv}
Version;CandidateModel;Candidates;Examples;Runs;Successes;TrialsMin;TrialsAvg;TrialsMax;Optimum;CostMin;CostAvg;CostMax;GapMin;GapAvg;GapMax;TimeMin;TimeAvg;TimeMax
\Homogeneous relaxation;ALPHA;5;142;95.07;81.61;405.5;405.5;405.5;2626438;2626438;2626846;2628161;0.0;7.97e-05;4.17e-04;1.521;3.514;7.463
\Tilted relaxation;ALPHA;5;142;95.07;81.61;405.5;405.5;405.5;2626438;2626438;2626846;2628161;0.0;7.97e-05;4.17e-04;\textcolor{red}{\textit{1.598}};\textcolor{red}{\textit{3.696}};\textcolor{red}{\textit{7.866}}
\Homogeneous relaxation;POPMUSIC;5;123;95.12;87.21;428.8;428.8;428.8;3031128;3031129;3031264;3032323;2.32e-03;7.21e-05;2.95e-04;1.898;3.830;8.205
\Tilted relaxation;POPMUSIC;5;123;95.12;87.21;428.8;428.8;428.8;3031128;3031129;3031264;3032323;2.32e-03;7.21e-05;2.95e-04;\textcolor{red}{\textit{1.988}};\textcolor{red}{\textit{4.006}};\textcolor{red}{\textit{8.602}}
\Homogeneous relaxation;POPMUSIC;100;123;100.0;20.39;251.6;251.6;251.6;Error;4156597;4156676;4156722;0.481;4.82e-03;4.82e-03;0.522;0.879;1.280
\Tilted relaxation;POPMUSIC;100;123;100.0;19.80;251.6;251.6;251.6;Error;\textcolor{green}{\textbf{4156587}};\textcolor{green}{\textbf{4156651}};\textcolor{green}{\textbf{4156670}};0.481;4.82e-03;4.82e-03;\textcolor{red}{\textit{0.561}};\textcolor{red}{\textit{0.934}};\textcolor{red}{\textit{1.454}}
\end{filecontents*}
\csvreader[head to column names, separator=semicolon]{Homogeneous_Latex_Summary_Small.csv}{}
    { \Version & \Candidates\ \CandidateModel & \Examples & \GapMin & \GapAvg & \TimeAvg\\}
    \\[-\normalbaselineskip]\hline
\end{tabular}
\label{fig:CR_small_instances_fully_relaxed}
\end{table}

\section{Conclusion}
The Lin-Kernighan-Helsgaun heuristic (LKH) is a highly developed and arguably the most powerful heuristic for solving the Traveling Salesman Problem.
We propose new and stronger variants of LKH by relaxing the positive gain criterion, namely the \homogeneous and the \tilted relaxation.
The \homogeneous relaxation contains the basic idea of relaxing the positive gain criterion, while the \tilted relaxation is a variant developed by algorithm engineering.
An extensive computational study shows that both variants have substantial improvements for large instances, while the \tilted relaxation is prevailing overall. For large instances, the \tilted variant is on average 13\% faster, for some instances about 30\%, compared to the best current version of LKH.
The quality of solutions stays virtually the same.

To the best of our knowledge this paper is the first systematic experimental study of LKH with a relaxed positive gain criterion.
This study advocates to use the relaxed positive gain criterion for solving large TSP instances and motivates to further investigate relaxations of the positive gain criterion.

\bibliographystyle{amsalpha}

\bibliography{bibliography}
%\clearpage
\appendix
\section{Additional tables}
\label{appendix}
\subsection{Median Summary}
\label{mediansummaries}
\begin{filecontents*}{Latex_Median_Summary_Large.csv}
Version;CandidateModel;Candidates;Examples;Runs;Successes;TrialsMin;TrialsAvg;TrialsMax;Optimum;CostMin;CostAvg;CostMax;GapMin;GapAvg;GapMax;TimeMin;TimeAvg;TimeMax
LKH 3.0.8;ALPHA;5;20;10.00;0.0;36567;36567;36567;102989;161763;161776;161786;6.20e-03;2.40e-04;3.96e-04;36776;41924;46393
LKH with \ref{C2**:tiltedgain};ALPHA;5;20;10.00;0.0;36567;36567;36567;102989;\textcolor{green}{\textbf{161751}};\textcolor{green}{\textbf{161771}};\textcolor{red}{\textit{161800}};\textcolor{red}{\textit{7.15e-03}};\textcolor{green}{\textbf{2.16e-04}};\textcolor{green}{\textbf{3.30e-04}};\textcolor{green}{\textbf{29734}};\textcolor{green}{\textbf{36522}};\textcolor{green}{\textbf{42755}}
LKH 3.0.8;POPMUSIC;5;20;10.00;0.0;36567;36567;36567;102989;161751;161771;161792;5.90e-03;1.98e-04;3.58e-04;41687;51549;60512
LKH with \ref{C2**:tiltedgain};POPMUSIC;5;20;10.00;0.0;36567;36567;36567;102989;\textcolor{green}{\textbf{161748}};\textcolor{green}{\textbf{161764}};\textcolor{green}{\textbf{161784}};\textcolor{green}{\textbf{4.85e-03}};\textcolor{green}{\textbf{1.90e-04}};\textcolor{green}{\textbf{3.09e-04}};\textcolor{green}{\textbf{37580}};\textcolor{green}{\textbf{44042}};\textcolor{green}{\textbf{52305}}
\end{filecontents*}
\csvreader[
longtable=lrrrrr,
table head=	\caption{Comparison of our variant using the \tilted relaxation \ref{C2**:tiltedgain} and the original LKH in minimal cost, average cost and average time per run on \textbf{large instances} with different candidate set types by computing the \textbf{median}. Green and bold indicate a better performance of our variant, red and italic a worse performance of our variant.
\label{fig:CR_median_large_instances}}\\
\hline
	Version & Candidates & \# Ex. & CostMin & CostAvg & TimeAvg% specify table head
\\ \hline \endfirsthead

\hline
\endfoot,
%table foot = \hline ,
late after line=\\,
head to column names, separator=semicolon]{Latex_Median_Summary_Large.csv}{}% use head of csv as column names
{\Version & \Candidates\ \CandidateModel & \Examples & \CostMin & \CostAvg & \TimeAvg \;s}
	
\begin{filecontents*}{Latex_Median_Summary_Medium.csv}
Version;CandidateModel;Candidates;Examples;Runs;Successes;TrialsMin;TrialsAvg;TrialsMax;Optimum;CostMin;CostAvg;CostMax;GapMin;GapAvg;GapMax;TimeMin;TimeAvg;TimeMax
LKH 3.0.8;ALPHA;5;140;10.00;3.000;3441;3441;3441;51826;60140;60148;60156;0.0;7.45e-05;2.21e-04;134.4;181.0;239.3
LKH with \ref{C2**:tiltedgain};ALPHA;5;140;10.00;2.500;3441;3441;3441;51826;\textcolor{green}{\textbf{60134}};\textcolor{red}{\textit{60150}};\textcolor{red}{\textit{60163}};0.0;\textcolor{red}{\textit{9.55e-05}};\textcolor{red}{\textit{3.03e-04}};\textcolor{red}{\textit{149.9}};\textcolor{red}{\textit{218.7}};\textcolor{red}{\textit{296.0}}
LKH 3.0.8;POPMUSIC;5;140;10.00;5.000;3441;3441;3441;51826;60134;60146;60160;0.0;4.60e-05;1.40e-04;150.5;194.8;258.2
LKH with \ref{C2**:tiltedgain};POPMUSIC;5;140;10.00;4.000;3441;3441;3441;51826;\textcolor{green}{\textbf{60132}};\textcolor{green}{\textbf{60141}};\textcolor{green}{\textbf{60150}};0.0;\textcolor{red}{\textit{4.85e-05}};\textcolor{red}{\textit{1.50e-04}};\textcolor{red}{\textit{162.5}};\textcolor{red}{\textit{229.1}};\textcolor{red}{\textit{310.7}}
LKH 3.0.8;POPMUSIC;100;4;100.0;0.0;3500;3500;3500;Error;65402158;65404645;65406064;0.0;2.75e-05;4.55e-05;5.390;9.735;17.93
LKH with \ref{C2**:tiltedgain};POPMUSIC;100;4;100.0;0.0;3500;3500;3500;Error;\textcolor{green}{\textbf{65398829}};\textcolor{green}{\textbf{65401888}};\textcolor{green}{\textbf{65402364}};0.0;\textcolor{green}{\textbf{2.70e-05}};\textcolor{green}{\textbf{3.30e-05}};\textcolor{green}{\textbf{5.370}};\textcolor{red}{\textit{12.62}};\textcolor{red}{\textit{24.75}}
\end{filecontents*}
\csvreader[
longtable=lrrrrr,
table head=	\caption{Comparison of our variant using the \tilted relaxation \ref{C2**:tiltedgain} and the original LKH in minimal gap, average gap and average time per run on \textbf{medium instances} with different candidate set types by computing the \textbf{median}.
	\label{fig:CR_median_medium_instances}}\\
\hline
Version & Candidates & \# Ex. & GapMin & GapAvg & TimeAvg  % specify table head
 \\ 
&&&  in \% &  in \% & in s\\ \hline \endfirsthead

\hline
\endfoot,
%table foot = \hline ,
late after line=\\,
head to column names, separator=semicolon]{Latex_Median_Summary_Medium.csv}{}% use head of csv as column names
{\Version & \Candidates\ \CandidateModel & \Examples & \GapMin & \GapAvg & \TimeAvg}
\newpage
\begin{filecontents*}{Latex_Median_Summary_Small.csv}
Version;CandidateModel;Candidates;Examples;Runs;Successes;TrialsMin;TrialsAvg;TrialsMax;Optimum;CostMin;CostAvg;CostMax;GapMin;GapAvg;GapMax;TimeMin;TimeAvg;TimeMax
LKH 3.0.8;ALPHA;5;142;100.0;100.0;227.5;227.5;227.5;18264;18264;18264;18264;0.0;0.0;0.0;0.190;0.415;0.600
LKH with \ref{C2**:tiltedgain};ALPHA;5;142;100.0;100.0;227.5;227.5;227.5;18264;18264;18264;18264;0.0;0.0;0.0;\textcolor{red}{\textit{0.195}};\textcolor{red}{\textit{0.470}};\textcolor{red}{\textit{0.675}}
LKH 3.0.8;POPMUSIC;5;123;100.0;100.0;280.0;280.0;280.0;26130;26130;26130;26130;0.0;0.0;0.0;0.260;0.480;0.750
LKH with \ref{C2**:tiltedgain};POPMUSIC;5;123;100.0;100.0;280.0;280.0;280.0;26130;26130;26130;26130;0.0;0.0;0.0;0.260;\textcolor{red}{\textit{0.650}};\textcolor{red}{\textit{0.980}}
LKH 3.0.8;POPMUSIC;100;123;100.0;0.0;235.0;235.0;235.0;Error;3803628;3803640;3803684;0.0;0.0;0.0;0.240;0.420;0.580
LKH with \ref{C2**:tiltedgain};POPMUSIC;100;123;100.0;0.0;235.0;235.0;235.0;Error;3803628;\textcolor{green}{\textbf{3803628}};\textcolor{green}{\textbf{3803628}};0.0;0.0;0.0;\textcolor{red}{\textit{0.450}};\textcolor{red}{\textit{0.640}};\textcolor{red}{\textit{1.020}}
\end{filecontents*}
\csvreader[
longtable=lrrrrr,
table head=	\caption{Comparison of our variant using the \tilted relaxation \ref{C2**:tiltedgain} and the original LKH in minimal gap, average gap and average time per run on \textbf{small instances} with different candidate set types by computing the \textbf{median}.
	\label{fig:CR_median_small_instances}}\\
\hline
 Version & Candidates & \# Ex. & GapMin & GapAvg & TimeAvg% specify table head
\\ 
&&&  in \% & in \% & in s\\ \hline \endfirsthead

\hline
\endfoot,
%table foot = \hline ,
late after line=\\,
head to column names, separator=semicolon]{Latex_Median_Summary_Small.csv}{}% use head of csv as column names
 {\Version & \Candidates\ \CandidateModel & \Examples & \GapMin & \TimeAvg & \TimeAvg}
 
\subsection{Percentile Summary}

\begin{filecontents*}{Latex_Quantil_Summary_Large.csv}
Version;CandidateModel;Candidates;Examples;Runs;Successes;TrialsMin;TrialsAvg;TrialsMax;Optimum;CostMin;CostAvg;CostMax;GapMin;GapAvg;GapMax;TimeMin;TimeAvg;TimeMax
LKH 3.0.8;ALPHA;5;20;10.00;0.050;100000;100000;100000;7640628;127297029;127307680;127321328;0.027;6.46e-04;2.22e-03;188223;245436;422117
LKH with \ref{C2**:tiltedgain};ALPHA;5;20;10.00;0.0;100000;100000;100000;7640628;\textcolor{red}{\textit{127300617}};\textcolor{red}{\textit{127307744}};\textcolor{green}{\textbf{127313466}};\textcolor{green}{\textbf{0.020}};\textcolor{red}{\textit{1.46e-03}};\textcolor{red}{\textit{4.41e-03}};\textcolor{green}{\textbf{153572}};\textcolor{green}{\textbf{177181}};\textcolor{green}{\textbf{267856}}
LKH 3.0.8;POPMUSIC;5;20;10.00;0.0;100000;100000;100000;7640628;127293502;127304503;127316702;0.022;3.24e-04;5.20e-04;205924;223968;245173
LKH with \ref{C2**:tiltedgain};POPMUSIC;5;20;10.00;0.0;100000;100000;100000;7640628;\textcolor{red}{\textit{127297125}};\textcolor{green}{\textbf{127304300}};\textcolor{green}{\textbf{127311539}};\textcolor{green}{\textbf{0.016}};\textcolor{red}{\textit{3.60e-04}};\textcolor{red}{\textit{1.00e-03}};\textcolor{red}{\textit{226074}};\textcolor{red}{\textit{238732}};\textcolor{red}{\textit{246300}}
\end{filecontents*}
\csvreader[
longtable=lrrrrr,
table head=	\caption{Comparison of our variant using the \tilted relaxation \ref{C2**:tiltedgain} and the original LKH in minimal cost, average cost and average time per run on \textbf{large instances} with different candidate set types by computing the \textbf{95th percentile}. Green and bold indicate a better performance of our variant, red and italic a worse performance of our variant.
	\label{fig:CR_quantil_large_instances}}\\
\hline
  Version & Candidates & \# Ex. & CostMin & CostAvg & TimeAvg% specify table head
\\ \hline \endfirsthead

\hline
\endfoot,
%table foot = \hline ,
late after line=\\,
head to column names, separator=semicolon]{Latex_Quantil_Summary_Large.csv}{}% use head of csv as column names
{\Version & \Candidates\ \CandidateModel & \Examples & \CostMin & \CostAvg & \TimeAvg \;s}%

\begin{filecontents*}{Latex_Quantil_Summary_Medium.csv}
Version;CandidateModel;Candidates;Examples;Runs;Successes;TrialsMin;TrialsAvg;TrialsMax;Optimum;CostMin;CostAvg;CostMax;GapMin;GapAvg;GapMax;TimeMin;TimeAvg;TimeMax
LKH 3.0.8;ALPHA;5;140;10.00;10.00;24148;24148;24148;19569316;40303989;40305342;40315395;0.017;6.12e-04;1.43e-03;12194;13688;16425
LKH with \ref{C2**:tiltedgain};ALPHA;5;140;10.00;10.00;24148;24148;24148;19569316;40303989;\textcolor{red}{\textit{40306223}};\textcolor{red}{\textit{40317643}};\textcolor{red}{\textit{0.023}};\textcolor{red}{\textit{6.90e-04}};\textcolor{red}{\textit{1.67e-03}};\textcolor{green}{\textbf{9874}};\textcolor{green}{\textbf{13352}};\textcolor{green}{\textbf{15431}}
LKH 3.0.8;POPMUSIC;5;140;10.00;10.00;24148;24148;24148;19569316;40303989;40304045;40304272;0.014;3.89e-04;8.50e-04;13719;15755;17272
LKH with \ref{C2**:tiltedgain};POPMUSIC;5;140;10.00;10.00;24148;24148;24148;19569316;\textcolor{red}{\textit{40304024}};\textcolor{green}{\textbf{40304024}};\textcolor{green}{\textbf{40304024}};\textcolor{red}{\textit{0.015}};\textcolor{red}{\textit{4.02e-04}};\textcolor{red}{\textit{9.06e-04}};\textcolor{green}{\textbf{12316}};\textcolor{green}{\textbf{14996}};\textcolor{red}{\textit{17739}}
LKH 3.0.8;POPMUSIC;100;4;100.0;0.0;4852;4852;4852;Error;90945752;90946250;90946348;0.0;3.98e-05;6.95e-05;10.13;15.77;22.63
LKH with \ref{C2**:tiltedgain};POPMUSIC;100;4;100.0;0.0;4852;4852;4852;Error;\textcolor{green}{\textbf{90945225}};\textcolor{red}{\textit{90948592}};\textcolor{red}{\textit{90949043}};0.0;\textcolor{red}{\textit{8.34e-05}};\textcolor{red}{\textit{9.00e-05}};\textcolor{red}{\textit{11.28}};\textcolor{red}{\textit{21.45}};\textcolor{red}{\textit{35.35}}
\end{filecontents*}
\csvreader[
longtable=lrrrrr,
table head=	\caption{Comparison of our variant using the \tilted relaxation \ref{C2**:tiltedgain} and the original LKH in minimal gap, average gap and average time per run on \textbf{medium instances} with different candidate set types by computing the \textbf{95th percentile}.
	\label{fig:CR_quantil_medium_instances}}\\
\hline
 Version & Candidates & \# Ex. & GapMin & GapAvg & TimeAvg % specify table head
\\ 
&&&  in \% & in \% & in s\\ \hline \endfirsthead

\hline
\endfoot,
%table foot = \hline ,
late after line=\\,
head to column names, separator=semicolon]{Latex_Quantil_Summary_Medium.csv}{}% use head of csv as column names
{\Version & \Candidates\ \CandidateModel & \Examples & \GapMin & \GapAvg & \TimeAvg}%

\begin{filecontents*}{Latex_Quantil_Summary_Small.csv}
Version;CandidateModel;Candidates;Examples;Runs;Successes;TrialsMin;TrialsAvg;TrialsMax;Optimum;CostMin;CostAvg;CostMax;GapMin;GapAvg;GapMax;TimeMin;TimeAvg;TimeMax
LKH 3.0.8;ALPHA;5;142;100.0;100.0;1000;1000;1000;22980365;22980365;22981320;22985197;0.0;5.25e-04;1.51e-03;7.199;15.09;29.50
LKH with \ref{C2**:tiltedgain};ALPHA;5;142;100.0;100.0;1000;1000;1000;22980365;22980365;\textcolor{red}{\textit{22982110}};\textcolor{red}{\textit{22985218}};0.0;\textcolor{red}{\textit{5.40e-04}};\textcolor{red}{\textit{2.09e-03}};\textcolor{red}{\textit{7.817}};\textcolor{green}{\textbf{14.79}};\textcolor{red}{\textit{33.44}}
LKH 3.0.8;POPMUSIC;5;123;100.0;100.0;1000;1000;1000;23019585;23019585;23019774;23022236;0.0;3.05e-04;1.73e-03;10.23;17.23;26.59
LKH with \ref{C2**:tiltedgain};POPMUSIC;5;123;100.0;100.0;1000;1000;1000;23019585;23019585;\textcolor{red}{\textit{23019902}};\textcolor{red}{\textit{23024568}};0.0;\textcolor{green}{\textbf{1.63e-04}};\textcolor{green}{\textbf{1.16e-03}};\textcolor{green}{\textbf{8.746}};\textcolor{green}{\textbf{14.91}};\textcolor{red}{\textit{46.30}}
LKH 3.0.8;POPMUSIC;100;123;100.0;100.0;399.7;399.7;399.7;2823926;6867910;6867910;6867910;2.934;0.029;0.029;0.835;1.368;1.913
LKH with \ref{C2**:tiltedgain};POPMUSIC;100;123;100.0;100.0;399.7;399.7;399.7;2823926;6867910;6867910;6867910;2.934;0.029;0.029;\textcolor{red}{\textit{1.564}};\textcolor{red}{\textit{2.863}};\textcolor{red}{\textit{4.379}}
\end{filecontents*}
\csvreader[
longtable=lrrrrr,
table head=	\caption{Comparison of our variant using the \tilted relaxation \ref{C2**:tiltedgain} and the original LKH in minimal gap, average gap and average time per run on \textbf{small instances} with different candidate set types by computing the \textbf{95th percentile}. 
	\label{fig:CR_quantil_small_instances}}\\
\hline
 Version & Candidates & \# Ex. & GapMin & GapAvg & TimeAvg % specify table head
\\ 
&&&  in \% &  in \% & in s\\ \hline \endfirsthead

\hline
\endfoot,
%table foot = \hline ,
late after line=\\,
head to column names, separator=semicolon]{Latex_Quantil_Summary_Small.csv}{}% use head of csv as column names
{\Version & \Candidates\ \CandidateModel & \Examples & \GapMin & \GapAvg & \TimeAvg}%

\subsection{Large instances}
\label{subsection:Add_Tables_large_instances}
\begin{filecontents*}{Latex_Comparison_Large_POPMUSIC.csv}
Problem;CandidateModel;Candidates;Optimum;RunsLKH;RunsOur;CostMinLKH;CostMinOur;CostAvgLKH;CostAvgOur;CostMaxLKH;CostMaxOur;GapMinLKH;GapMinOur;GapAvgLKH;GapAvgOur;GapMaxLKH;GapMaxOur;TimeMinLKH;TimeMinOur;TimeAvgLKH;TimeAvgOur;TimeMaxLKH;TimeMaxOur;CostMinChange;CostAvgChange;CostMaxChange;GapMinChange;GapAvgChange;GapMaxChange;TimeMinChange;TimeAvgChange;TimeMaxChange
DIMACS\_C31k.1;POPMUSIC;5;Error;10.00;10.00;59293266;59293266;59295163;\textcolor{red}{\textit{59304642}};59303429;\textcolor{red}{\textit{59349001}};0.0;0.0;3.20e-05;\textcolor{red}{\textit{1.92e-04}};1.71e-04;\textcolor{red}{\textit{9.40e-04}};14214;\textcolor{green}{\textbf{11629}};15658;\textcolor{green}{\textbf{13044}};17439;\textcolor{green}{\textbf{16144}};$\pm$ 0;+1.60e-02;+7.68e-02;$\pm$ 0;+500.00;+449.71;-18.19;-16.69;-7.43
DIMACS\_C31k.0;POPMUSIC;5;Error;10.00;10.00;59545390;59545390;59575202;\textcolor{red}{\textit{59604908}};59587149;\textcolor{red}{\textit{59673921}};0.0;0.0;5.01e-04;\textcolor{red}{\textit{1.00e-03}};7.01e-04;\textcolor{red}{\textit{2.16e-03}};14871;\textcolor{red}{\textit{15087}};18698;\textcolor{green}{\textbf{16701}};22081;\textcolor{green}{\textbf{20026}};$\pm$ 0;+4.99e-02;+0.146;$\pm$ 0;+99.60;+208.13;+1.45;-10.68;-9.31
vlsi\_bby34656;POPMUSIC;5;99159;10.00;10.00;99162;\textcolor{red}{\textit{99164}};99170;\textcolor{red}{\textit{99172}};99184;\textcolor{red}{\textit{99189}};3e-05;\textcolor{red}{\textit{5e-05}};1.12e-04;\textcolor{red}{\textit{1.31e-04}};2.52e-04;\textcolor{red}{\textit{3.03e-04}};26049;\textcolor{green}{\textbf{22380}};32498;\textcolor{green}{\textbf{29008}};38803;\textcolor{green}{\textbf{33051}};+2.02e-03;+2.02e-03;+5.04e-03;+66.67;+16.96;+20.24;-14.08;-10.74;-14.82
DIMACS\_E31k.0;POPMUSIC;5;Error;10.00;10.00;127284898;\textcolor{red}{\textit{127288632}};127295977;\textcolor{green}{\textbf{127295925}};127308208;\textcolor{green}{\textbf{127303176}};0.0;0.0;8.70e-05;\textcolor{green}{\textbf{5.70e-05}};1.83e-04;\textcolor{green}{\textbf{1.14e-04}};33934;\textcolor{green}{\textbf{22233}};38983;\textcolor{green}{\textbf{27460}};42552;\textcolor{green}{\textbf{33582}};+2.93e-03;-4.08e-05;-3.95e-03;$\pm$ 0;-34.48;-37.70;-34.48;-29.56;-21.08
DIMACS\_E31k.1;POPMUSIC;5;Error;10.00;10.00;127456986;\textcolor{red}{\textit{127458487}};127466483;\textcolor{green}{\textbf{127463432}};127478086;\textcolor{green}{\textbf{127470427}};0.0;0.0;7.50e-05;\textcolor{green}{\textbf{3.90e-05}};1.66e-04;\textcolor{green}{\textbf{9.40e-05}};35206;\textcolor{green}{\textbf{22571}};42157;\textcolor{green}{\textbf{29123}};49628;\textcolor{green}{\textbf{32684}};+1.18e-03;-2.39e-03;-6.01e-03;$\pm$ 0;-48.00;-43.37;-35.89;-30.92;-34.14
National\_TSPs\_bm33708;POPMUSIC;5;959289;10.00;10.00;959318;\textcolor{green}{\textbf{959308}};959356;\textcolor{green}{\textbf{959339}};959409;\textcolor{green}{\textbf{959398}};3e-05;\textcolor{green}{\textbf{2e-05}};6.90e-05;\textcolor{green}{\textbf{5.30e-05}};1.25e-04;\textcolor{green}{\textbf{1.14e-04}};37086;\textcolor{green}{\textbf{32965}};48001;\textcolor{green}{\textbf{36714}};54539;\textcolor{green}{\textbf{40882}};-1.04e-03;-1.77e-03;-1.15e-03;-33.33;-23.19;-8.80;-11.11;-23.51;-25.04
vlsi\_bna56769;POPMUSIC;5;158078;10.00;10.00;158112;\textcolor{green}{\textbf{158097}};158126;\textcolor{green}{\textbf{158112}};158144;\textcolor{green}{\textbf{158128}};2.15e-04;\textcolor{green}{\textbf{1.20e-04}};3.02e-04;\textcolor{green}{\textbf{2.15e-04}};4.18e-04;\textcolor{green}{\textbf{3.16e-04}};117332;\textcolor{green}{\textbf{97583}};130540;\textcolor{green}{\textbf{118318}};148043;\textcolor{green}{\textbf{136731}};-9.49e-03;-8.85e-03;-1.01e-02;-44.19;-28.81;-24.40;-16.83;-9.36;-7.64
vlsi\_dan59296;POPMUSIC;5;165371;10.00;10.00;165390;\textcolor{red}{\textit{165399}};165417;165417;165439;\textcolor{red}{\textit{165440}};1.15e-04;\textcolor{red}{\textit{1.69e-04}};2.78e-04;\textcolor{green}{\textbf{2.77e-04}};4.11e-04;\textcolor{red}{\textit{4.17e-04}};108444;\textcolor{green}{\textbf{100837}};140333;\textcolor{green}{\textbf{113975}};181983;\textcolor{green}{\textbf{140758}};+5.44e-03;$\pm$ 0;+6.04e-04;+46.96;-0.360;+1.46;-7.01;-18.78;-22.65
DIMACS\_C100k.0;POPMUSIC;5;Error;10.00;10.00;104633917;\textcolor{green}{\textbf{104617752}};104657832;\textcolor{green}{\textbf{104641010}};104672048;\textcolor{green}{\textbf{104657755}};0.0;0.0;2.29e-04;\textcolor{green}{\textbf{2.22e-04}};3.64e-04;\textcolor{red}{\textit{3.82e-04}};178845;\textcolor{green}{\textbf{145317}};200368;\textcolor{green}{\textbf{158532}};217955;\textcolor{green}{\textbf{174702}};-1.54e-02;-1.61e-02;-1.37e-02;$\pm$ 0;-3.06;+4.95;-18.75;-20.88;-19.84
TSP\_pla33810;POPMUSIC;5;66048945;10.00;10.00;66057180;\textcolor{green}{\textbf{66049097}};66067432;\textcolor{green}{\textbf{66061323}};66077718;\textcolor{green}{\textbf{66068696}};1.25e-04;\textcolor{green}{\textbf{2e-06}};2.80e-04;\textcolor{green}{\textbf{1.87e-04}};4.36e-04;\textcolor{green}{\textbf{2.99e-04}};197374;\textcolor{red}{\textit{225260}};207595;\textcolor{red}{\textit{236568}};224017;\textcolor{red}{\textit{243095}};-1.22e-02;-9.25e-03;-1.37e-02;-98.40;-33.21;-31.42;+14.13;+13.96;+8.52
DIMACS\_C100k.1;POPMUSIC;5;Error;10.00;10.00;105416845;\textcolor{green}{\textbf{105390777}};105446335;\textcolor{green}{\textbf{105425123}};105470562;\textcolor{green}{\textbf{105469042}};0.0;0.0;2.80e-04;\textcolor{red}{\textit{3.26e-04}};5.10e-04;\textcolor{red}{\textit{7.43e-04}};199253;\textcolor{green}{\textbf{171393}};216743;\textcolor{green}{\textbf{179302}};237226;\textcolor{green}{\textbf{191656}};-2.47e-02;-2.01e-02;-1.44e-03;$\pm$ 0;+16.43;+45.69;-13.98;-17.27;-19.21
National\_TSPs\_ch71009;POPMUSIC;5;4566506;10.00;10.00;4567023;\textcolor{green}{\textbf{4566684}};4567325;\textcolor{green}{\textbf{4566928}};4567697;\textcolor{green}{\textbf{4567207}};1.13e-04;\textcolor{green}{\textbf{3.90e-05}};1.79e-04;\textcolor{green}{\textbf{9.20e-05}};2.61e-04;\textcolor{green}{\textbf{1.54e-04}};332668;\textcolor{green}{\textbf{241528}};361241;\textcolor{green}{\textbf{279836}};396156;\textcolor{green}{\textbf{307192}};-7.42e-03;-8.69e-03;-1.07e-02;-65.49;-48.60;-41.00;-27.40;-22.53;-22.46
vlsi\_pbh30440;POPMUSIC;5;88313;10.00;10.00;88320;88320;88329;\textcolor{red}{\textit{88330}};88346;\textcolor{green}{\textbf{88342}};7.90e-05;7.90e-05;1.79e-04;\textcolor{red}{\textit{1.91e-04}};3.74e-04;\textcolor{green}{\textbf{3.28e-04}};15942;\textcolor{green}{\textbf{15214}};19752;\textcolor{green}{\textbf{18126}};24263;\textcolor{green}{\textbf{23369}};$\pm$ 0;+1.13e-03;-4.53e-03;$\pm$ 0;+6.70;-12.30;-4.57;-8.23;-3.68
vlsi\_fry33203;POPMUSIC;5;97240;10.00;10.00;97253;97253;97264;\textcolor{red}{\textit{97267}};97282;\textcolor{red}{\textit{97290}};1.34e-04;1.34e-04;2.51e-04;\textcolor{red}{\textit{2.79e-04}};4.32e-04;\textcolor{red}{\textit{5.14e-04}};22364;\textcolor{green}{\textbf{19295}};25535;\textcolor{green}{\textbf{23355}};29131;\textcolor{green}{\textbf{26976}};$\pm$ 0;+3.08e-03;+8.22e-03;$\pm$ 0;+11.16;+18.98;-13.72;-8.54;-7.40
vlsi\_xib32892;POPMUSIC;5;96757;10.00;10.00;96763;\textcolor{red}{\textit{96772}};96777;\textcolor{red}{\textit{96787}};96793;\textcolor{red}{\textit{96803}};6.20e-05;\textcolor{red}{\textit{1.55e-04}};2.04e-04;\textcolor{red}{\textit{3.10e-04}};3.72e-04;\textcolor{red}{\textit{4.75e-04}};26138;\textcolor{green}{\textbf{21434}};27994;\textcolor{green}{\textbf{26394}};34027;\textcolor{red}{\textit{34408}};+9.30e-03;+1.03e-02;+1.03e-02;+150.00;+51.96;+27.69;-18.00;-5.72;+1.12
vlsi\_pba38478;POPMUSIC;5;108318;10.00;10.00;108327;\textcolor{green}{\textbf{108326}};108340;\textcolor{green}{\textbf{108336}};108356;\textcolor{green}{\textbf{108344}};8.30e-05;\textcolor{green}{\textbf{7.40e-05}};2.05e-04;\textcolor{green}{\textbf{1.63e-04}};3.51e-04;\textcolor{green}{\textbf{2.40e-04}};28200;\textcolor{red}{\textit{29869}};36712;\textcolor{green}{\textbf{36368}};43687;\textcolor{green}{\textbf{42756}};-9.23e-04;-3.69e-03;-1.11e-02;-10.84;-20.49;-31.62;+5.92;-0.94;-2.13
vlsi\_ics39603;POPMUSIC;5;106819;10.00;10.00;106825;\textcolor{green}{\textbf{106824}};106834;\textcolor{green}{\textbf{106831}};106854;\textcolor{green}{\textbf{106844}};5.60e-05;\textcolor{green}{\textbf{4.70e-05}};1.45e-04;\textcolor{green}{\textbf{1.15e-04}};3.28e-04;\textcolor{green}{\textbf{2.34e-04}};46288;\textcolor{green}{\textbf{42195}};55097;\textcolor{green}{\textbf{51369}};66486;\textcolor{green}{\textbf{61853}};-9.36e-04;-2.81e-03;-9.36e-03;-16.07;-20.69;-28.66;-8.84;-6.77;-6.97
vlsi\_rbz43748;POPMUSIC;5;125183;10.00;10.00;125190;\textcolor{red}{\textit{125200}};125206;125206;125222;\textcolor{green}{\textbf{125217}};5.60e-05;\textcolor{red}{\textit{1.36e-04}};1.87e-04;\textcolor{green}{\textbf{1.85e-04}};3.12e-04;\textcolor{green}{\textbf{2.72e-04}};56174;\textcolor{green}{\textbf{49261}};64024;\textcolor{green}{\textbf{57064}};75148;\textcolor{green}{\textbf{68673}};+7.99e-03;$\pm$ 0;-3.99e-03;+142.86;-1.07;-12.82;-12.31;-10.87;-8.62
vlsi\_fna52057;POPMUSIC;5;147789;10.00;10.00;147812;\textcolor{green}{\textbf{147806}};147817;147817;147832;\textcolor{green}{\textbf{147825}};1.56e-04;\textcolor{green}{\textbf{1.15e-04}};1.91e-04;\textcolor{green}{\textbf{1.90e-04}};2.91e-04;\textcolor{green}{\textbf{2.44e-04}};62476;\textcolor{green}{\textbf{55101}};70470;\textcolor{green}{\textbf{63181}};87790;\textcolor{green}{\textbf{67682}};-4.06e-03;$\pm$ 0;-4.74e-03;-26.28;-0.524;-16.15;-11.80;-10.34;-22.90
vlsi\_fht47608;POPMUSIC;5;125104;10.00;10.00;125134;\textcolor{green}{\textbf{125121}};125143;\textcolor{green}{\textbf{125136}};125162;\textcolor{green}{\textbf{125149}};2.40e-04;\textcolor{green}{\textbf{1.36e-04}};3.15e-04;\textcolor{green}{\textbf{2.57e-04}};4.64e-04;\textcolor{green}{\textbf{3.60e-04}};63003;\textcolor{green}{\textbf{51880}};73115;\textcolor{green}{\textbf{62841}};82211;\textcolor{green}{\textbf{72675}};-1.04e-02;-5.59e-03;-1.04e-02;-43.33;-18.41;-22.41;-17.65;-14.05;-11.60
\hline Average;POPMUSIC;5.000;Error;10.00;10.00;32826656;\textcolor{green}{\textbf{32824384}};32832477;\textcolor{green}{\textbf{32832052}};32837146;\textcolor{red}{\textit{32841860}};7.47e-05;\textcolor{green}{\textbf{6.38e-05}};2.05e-04;\textcolor{red}{\textit{2.24e-04}};3.46e-04;\textcolor{red}{\textit{4.35e-04}};80793;\textcolor{green}{\textbf{69651}};91276;\textcolor{green}{\textbf{78864}};103658;\textcolor{green}{\textbf{88445}};-2.89e-03;1.09e-04;7.09e-03;nan;21.10;23.30;-13.16;-13.12;-12.86
Median;POPMUSIC;5.000;102989;10.00;10.00;161751;\textcolor{green}{\textbf{161748}};161771;\textcolor{green}{\textbf{161764}};161792;\textcolor{green}{\textbf{161784}};5.90e-05;\textcolor{green}{\textbf{4.85e-05}};1.98e-04;\textcolor{green}{\textbf{1.90e-04}};3.57e-04;\textcolor{green}{\textbf{3.09e-04}};41687;\textcolor{green}{\textbf{37580}};51549;\textcolor{green}{\textbf{44042}};60512;\textcolor{green}{\textbf{52305}};-4.62e-04;-9.64e-04;-4.26e-03;nan;-2.06;-12.56;-13.85;-10.81;-10.45
95th Quantil;POPMUSIC;5.000;7640628;10.00;10.00;127293502;\textcolor{red}{\textit{127297125}};127304503;\textcolor{green}{\textbf{127304300}};127316702;\textcolor{green}{\textbf{127311539}};2.16e-04;\textcolor{green}{\textbf{1.56e-04}};3.24e-04;\textcolor{red}{\textit{3.60e-04}};5.20e-04;\textcolor{red}{\textit{1.00e-03}};205924;\textcolor{red}{\textit{226074}};223968;\textcolor{red}{\textit{238732}};245173;\textcolor{red}{\textit{246300}};8.05e-03;1.80e-02;8.00e-02;nan;119.60;220.10;6.33;-0.192;1.49
\end{filecontents*}
\csvreader[
  longtable=lrrrrrr,
  table head=\caption{Direct comparison of our variant using the \tilted relaxation \ref{C2**:tiltedgain} and the original LKH with 5 POPMUSIC candidates in minimal cost and average time per run for all \textbf{large finished instances}. Green and bold indicate a better performance of our variant, red and italic a worse performance of our variant.
\label{fig:CR_large_instances_POPMUSIC_direct}}\\
    \hline
 \multicolumn{1}{c}{Problem} & \multicolumn{2}{c}{CostMin} & %\multicolumn{2}{c}{CostAvg} & 
\multicolumn{2}{c}{TimeAvg in s} & \multicolumn{1}{c}{TimeAvg} \\ 
& LKH 3.0.8 & \ref{C2**:tiltedgain} & LKH 3.0.8 & \ref{C2**:tiltedgain} & \multicolumn{1}{c}{ratio in \%}\\ \hline \endfirsthead
\hline \multicolumn{1}{c}{Problem} & \multicolumn{2}{c}{CostMin} & %\multicolumn{2}{c}{CostAvg} & 
\multicolumn{2}{c}{TimeAvg in s} & \multicolumn{1}{c}{TimeAvg} \\   & LKH 3.0.8 & \ref{C2**:tiltedgain} & LKH 3.0.8 & \ref{C2**:tiltedgain} & \multicolumn{1}{c}{ratio in \%} \\ \hline \endhead
\hline
    \endfoot,
    %table foot = \hline ,
  late after line=\\,
head to column names, separator=semicolon]{Latex_Comparison_Large_POPMUSIC.csv}{}% use head of csv as column names
    {\Problem & \CostMinLKH & \CostMinOur & \TimeAvgLKH & \TimeAvgOur & \TimeAvgChange}
    
    \begin{filecontents*}{Latex_Comparison_Large_ALPHA.csv}
Problem;CandidateModel;Candidates;Optimum;RunsLKH;RunsOur;CostMinLKH;CostMinOur;CostAvgLKH;CostAvgOur;CostMaxLKH;CostMaxOur;GapMinLKH;GapMinOur;GapAvgLKH;GapAvgOur;GapMaxLKH;GapMaxOur;TimeMinLKH;TimeMinOur;TimeAvgLKH;TimeAvgOur;TimeMaxLKH;TimeMaxOur;CostMinChange;CostAvgChange;CostMaxChange;GapMinChange;GapAvgChange;GapMaxChange;TimeMinChange;TimeAvgChange;TimeMaxChange
DIMACS\_C31k.0;ALPHA;5;Error;10.00;10.00;59556652;\textcolor{green}{\textbf{59545390}};59576964;\textcolor{green}{\textbf{59548465}};59591993;\textcolor{green}{\textbf{59557881}};0.0;0.0;3.41e-04;\textcolor{green}{\textbf{5.20e-05}};5.93e-04;\textcolor{green}{\textbf{2.10e-04}};12879;\textcolor{green}{\textbf{10096}};13819;\textcolor{green}{\textbf{10936}};16729;\textcolor{green}{\textbf{14142}};-1.89e-02;-4.78e-02;-5.72e-02;$\pm$ 0;-84.75;-64.59;-21.61;-20.86;-15.46
DIMACS\_C31k.1;ALPHA;5;Error;10.00;10.00;59407594;\textcolor{red}{\textit{59655194}};59530130;\textcolor{red}{\textit{59836059}};59786763;\textcolor{red}{\textit{60281404}};0.0;0.0;2.06e-03;\textcolor{red}{\textit{3.03e-03}};6.38e-03;\textcolor{red}{\textit{0.010}};17356;\textcolor{green}{\textbf{10631}};21859;\textcolor{green}{\textbf{12458}};39064;\textcolor{green}{\textbf{17528}};+0.417;+0.514;+0.827;$\pm$ 0;+47.09;+56.74;-38.75;-43.01;-55.13
vlsi\_pbh30440;ALPHA;5;88313;10.00;10.00;88317;\textcolor{red}{\textit{88324}};88331;\textcolor{red}{\textit{88332}};88349;\textcolor{green}{\textbf{88340}};4.50e-05;\textcolor{red}{\textit{1.25e-04}};1.99e-04;\textcolor{red}{\textit{2.14e-04}};4.08e-04;\textcolor{green}{\textbf{3.06e-04}};13729;\textcolor{green}{\textbf{13346}};17770;\textcolor{green}{\textbf{16367}};25943;\textcolor{green}{\textbf{22287}};+7.93e-03;+1.13e-03;-1.02e-02;+177.78;+7.54;-25.00;-2.79;-7.90;-14.09
vlsi\_fry33203;ALPHA;5;97240;10.00;10.00;97259;\textcolor{green}{\textbf{97245}};97266;\textcolor{green}{\textbf{97261}};97272;\textcolor{red}{\textit{97280}};1.95e-04;\textcolor{green}{\textbf{5.10e-05}};2.68e-04;\textcolor{green}{\textbf{2.17e-04}};3.29e-04;\textcolor{red}{\textit{4.11e-04}};19629;\textcolor{green}{\textbf{17705}};22033;\textcolor{green}{\textbf{20035}};25391;\textcolor{green}{\textbf{24342}};-1.44e-02;-5.14e-03;+8.22e-03;-73.85;-19.03;+24.92;-9.80;-9.07;-4.13
vlsi\_xib32892;ALPHA;5;96757;10.00;10.00;96779;\textcolor{red}{\textit{96782}};96797;\textcolor{red}{\textit{96798}};96823;\textcolor{green}{\textbf{96813}};2.27e-04;\textcolor{red}{\textit{2.58e-04}};4.17e-04;\textcolor{red}{\textit{4.26e-04}};6.82e-04;\textcolor{green}{\textbf{5.79e-04}};19194;\textcolor{green}{\textbf{17477}};23416;\textcolor{green}{\textbf{21503}};26286;\textcolor{red}{\textit{27140}};+3.10e-03;+1.03e-03;-1.03e-02;+13.66;+2.16;-15.10;-8.95;-8.17;+3.25
vlsi\_bby34656;ALPHA;5;99159;10.00;10.00;99166;99166;99177;99177;99185;\textcolor{red}{\textit{99189}};7.10e-05;7.10e-05;1.81e-04;\textcolor{red}{\textit{1.84e-04}};2.62e-04;\textcolor{red}{\textit{3.03e-04}};21777;\textcolor{green}{\textbf{16945}};24384;\textcolor{green}{\textbf{22177}};31299;\textcolor{green}{\textbf{28180}};$\pm$ 0;$\pm$ 0;+4.03e-03;$\pm$ 0;+1.66;+15.65;-22.19;-9.05;-9.97
vlsi\_pba38478;ALPHA;5;108318;10.00;10.00;108318;\textcolor{red}{\textit{108329}};108338;\textcolor{red}{\textit{108343}};108362;\textcolor{green}{\textbf{108357}};0.0;\textcolor{red}{\textit{1.02e-04}};1.85e-04;\textcolor{red}{\textit{2.33e-04}};4.06e-04;\textcolor{green}{\textbf{3.60e-04}};24431;\textcolor{green}{\textbf{22471}};27484;\textcolor{red}{\textit{28693}};34843;\textcolor{red}{\textit{35045}};+1.02e-02;+4.62e-03;-4.61e-03;+0.00e+00;+25.95;-11.33;-8.02;+4.40;+0.58
DIMACS\_E31k.0;ALPHA;5;Error;10.00;10.00;127288512;\textcolor{red}{\textit{127292616}};127299383;\textcolor{red}{\textit{127299452}};127313390;\textcolor{green}{\textbf{127305035}};0.0;0.0;8.50e-05;\textcolor{green}{\textbf{5.40e-05}};1.95e-04;\textcolor{green}{\textbf{9.80e-05}};31503;\textcolor{green}{\textbf{19540}};36829;\textcolor{green}{\textbf{22063}};41713;\textcolor{green}{\textbf{24634}};+3.22e-03;+5.42e-05;-6.56e-03;$\pm$ 0;-36.47;-49.74;-37.97;-40.09;-40.94
DIMACS\_E31k.1;ALPHA;5;Error;10.00;10.00;127458856;\textcolor{green}{\textbf{127452642}};127465322;\textcolor{green}{\textbf{127465293}};127472143;\textcolor{red}{\textit{127473658}};0.0;0.0;5.10e-05;\textcolor{red}{\textit{9.90e-05}};1.04e-04;\textcolor{red}{\textit{1.65e-04}};31865;\textcolor{green}{\textbf{19139}};37729;\textcolor{green}{\textbf{24552}};44994;\textcolor{green}{\textbf{31036}};-4.88e-03;-2.28e-05;+1.19e-03;$\pm$ 0;+94.12;+58.65;-39.94;-34.93;-31.02
National\_TSPs\_bm33708;ALPHA;5;959289;10.00;10.00;959328;\textcolor{green}{\textbf{959300}};959376;\textcolor{green}{\textbf{959344}};959464;\textcolor{green}{\textbf{959406}};4.10e-05;\textcolor{green}{\textbf{1.10e-05}};9.00e-05;\textcolor{green}{\textbf{5.70e-05}};1.82e-04;\textcolor{green}{\textbf{1.22e-04}};34680;\textcolor{green}{\textbf{30049}};40564;\textcolor{green}{\textbf{33538}};44806;\textcolor{green}{\textbf{41072}};-2.92e-03;-3.34e-03;-6.05e-03;-73.17;-36.67;-32.97;-13.35;-17.32;-8.33
vlsi\_ics39603;ALPHA;5;106819;10.00;10.00;106832;\textcolor{green}{\textbf{106828}};106844;\textcolor{green}{\textbf{106837}};106861;\textcolor{green}{\textbf{106849}};1.22e-04;\textcolor{green}{\textbf{8.40e-05}};2.31e-04;\textcolor{green}{\textbf{1.70e-04}};3.93e-04;\textcolor{green}{\textbf{2.81e-04}};38871;\textcolor{green}{\textbf{35231}};43284;\textcolor{green}{\textbf{39507}};47791;\textcolor{green}{\textbf{44438}};-3.74e-03;-6.55e-03;-1.12e-02;-31.15;-26.41;-28.50;-9.36;-8.73;-7.02
vlsi\_rbz43748;ALPHA;5;125183;10.00;10.00;125193;\textcolor{red}{\textit{125202}};125214;\textcolor{red}{\textit{125217}};125227;\textcolor{red}{\textit{125233}};8e-05;\textcolor{red}{\textit{1.52e-04}};2.49e-04;\textcolor{red}{\textit{2.69e-04}};3.51e-04;\textcolor{red}{\textit{3.99e-04}};38900;\textcolor{green}{\textbf{29419}};46904;\textcolor{green}{\textbf{39891}};56693;\textcolor{green}{\textbf{48208}};+7.19e-03;+2.40e-03;+4.79e-03;+90.00;+8.03;+13.68;-24.37;-14.95;-14.97
vlsi\_fna52057;ALPHA;5;147789;10.00;10.00;147802;\textcolor{red}{\textit{147805}};147813;\textcolor{red}{\textit{147820}};147831;\textcolor{red}{\textit{147836}};8.80e-05;\textcolor{red}{\textit{1.08e-04}};1.64e-04;\textcolor{red}{\textit{2.07e-04}};2.84e-04;\textcolor{red}{\textit{3.18e-04}};49032;\textcolor{green}{\textbf{43539}};57145;\textcolor{green}{\textbf{49625}};62956;\textcolor{green}{\textbf{55690}};+2.03e-03;+4.74e-03;+3.38e-03;+22.73;+26.22;+11.97;-11.20;-13.16;-11.54
vlsi\_fht47608;ALPHA;5;125104;10.00;10.00;125141;\textcolor{green}{\textbf{125128}};125154;\textcolor{green}{\textbf{125147}};125167;\textcolor{red}{\textit{125173}};2.96e-04;\textcolor{green}{\textbf{1.92e-04}};4.00e-04;\textcolor{green}{\textbf{3.41e-04}};5.04e-04;\textcolor{red}{\textit{5.52e-04}};65102;\textcolor{green}{\textbf{44416}};68980;\textcolor{green}{\textbf{51986}};76751;\textcolor{green}{\textbf{58712}};-1.04e-02;-5.59e-03;+4.79e-03;-35.14;-14.75;+9.52;-31.77;-24.64;-23.50
vlsi\_bna56769;ALPHA;5;158078;10.00;10.00;158110;\textcolor{green}{\textbf{158100}};158122;\textcolor{green}{\textbf{158115}};158130;\textcolor{red}{\textit{158132}};2.02e-04;\textcolor{green}{\textbf{1.39e-04}};2.78e-04;\textcolor{green}{\textbf{2.34e-04}};3.29e-04;\textcolor{red}{\textit{3.42e-04}};80157;\textcolor{green}{\textbf{67788}};92607;\textcolor{green}{\textbf{81798}};107314;\textcolor{green}{\textbf{97641}};-6.32e-03;-4.43e-03;+1.26e-03;-31.19;-15.83;+3.95;-15.43;-11.67;-9.01
vlsi\_dan59296;ALPHA;5;165371;10.00;10.00;165416;\textcolor{green}{\textbf{165402}};165430;\textcolor{green}{\textbf{165426}};165442;\textcolor{red}{\textit{165468}};2.72e-04;\textcolor{green}{\textbf{1.87e-04}};3.54e-04;\textcolor{green}{\textbf{3.33e-04}};4.29e-04;\textcolor{red}{\textit{5.87e-04}};91154;\textcolor{green}{\textbf{72738}};100765;\textcolor{green}{\textbf{83096}};119332;\textcolor{green}{\textbf{94770}};-8.46e-03;-2.42e-03;+1.57e-02;-31.25;-5.93;+36.83;-20.20;-17.53;-20.58
TSP\_pla33810;ALPHA;5;66048945;10.00;10.00;66054712;\textcolor{green}{\textbf{66053728}};66061156;\textcolor{red}{\textit{66061673}};66075351;\textcolor{green}{\textbf{66067857}};8.70e-05;\textcolor{green}{\textbf{7.20e-05}};1.85e-04;\textcolor{red}{\textit{1.93e-04}};4.00e-04;\textcolor{green}{\textbf{2.86e-04}};100823;\textcolor{red}{\textit{101673}};119333;\textcolor{red}{\textit{120625}};135925;\textcolor{red}{\textit{143467}};-1.49e-03;+7.83e-04;-1.13e-02;-17.24;+4.32;-28.50;+0.84;+1.08;+5.55
DIMACS\_C100k.1;ALPHA;5;Error;10.00;10.00;105512671;\textcolor{green}{\textbf{105465175}};105572893;\textcolor{red}{\textit{105610958}};105635409;\textcolor{red}{\textit{105832131}};0.0;0.0;5.71e-04;\textcolor{red}{\textit{1.38e-03}};1.16e-03;\textcolor{red}{\textit{3.48e-03}};184124;\textcolor{green}{\textbf{151755}};220936;\textcolor{green}{\textbf{175830}};414418;\textcolor{green}{\textbf{267055}};-4.50e-02;+3.61e-02;+0.186;$\pm$ 0;+141.68;+200.00;-17.58;-20.42;-35.56
National\_TSPs\_ch71009;ALPHA;5;4566506;10.00;10.00;4566747;\textcolor{green}{\textbf{4566746}};4567078;\textcolor{green}{\textbf{4566935}};4567415;\textcolor{green}{\textbf{4567152}};5.30e-05;5.30e-05;1.25e-04;\textcolor{green}{\textbf{9.40e-05}};1.99e-04;\textcolor{green}{\textbf{1.41e-04}};266111;\textcolor{green}{\textbf{188100}};284184;\textcolor{green}{\textbf{202845}};308198;\textcolor{green}{\textbf{212445}};-2.19e-05;-3.13e-03;-5.76e-03;$\pm$ 0;-24.80;-29.15;-29.32;-28.62;-31.07
DIMACS\_C100k.0;ALPHA;5;Error;10.00;10.00;104716990;\textcolor{red}{\textit{104814645}};104769222;\textcolor{red}{\textit{104940696}};104926802;\textcolor{red}{\textit{105243529}};0.0;0.0;4.99e-04;\textcolor{red}{\textit{1.20e-03}};2.00e-03;\textcolor{red}{\textit{4.09e-03}};177392;\textcolor{green}{\textbf{122397}};243397;\textcolor{green}{\textbf{165001}};568401;\textcolor{green}{\textbf{283092}};+9.33e-02;+0.164;+0.302;$\pm$ 0;+140.48;+104.50;-31.00;-32.21;-50.20
\hline Average;ALPHA;5.000;Error;10.00;10.00;32842020;\textcolor{red}{\textit{32856187}};32856000;\textcolor{red}{\textit{32880367}};32882369;\textcolor{red}{\textit{32930336}};8.89e-05;\textcolor{green}{\textbf{8.03e-05}};3.47e-04;\textcolor{red}{\textit{4.50e-04}};7.80e-04;\textcolor{red}{\textit{1.18e-03}};65935;\textcolor{green}{\textbf{51723}};77171;\textcolor{green}{\textbf{61126}};111442;\textcolor{green}{\textbf{78546}};2.10e-02;3.30e-02;6.20e-02;nan;11.77;12.90;-19.64;-17.84;-18.66
Median;ALPHA;5.000;102989;10.00;10.00;161763;\textcolor{green}{\textbf{161751}};161776;\textcolor{green}{\textbf{161771}};161786;\textcolor{red}{\textit{161800}};6.20e-05;\textcolor{red}{\textit{7.15e-05}};2.40e-04;\textcolor{green}{\textbf{2.16e-04}};3.97e-04;\textcolor{green}{\textbf{3.30e-04}};36776;\textcolor{green}{\textbf{29734}};41924;\textcolor{green}{\textbf{36522}};46393;\textcolor{green}{\textbf{42755}};-7.56e-04;1.79e-04;1.23e-03;nan;1.91;6.74;-18.89;-16.14;-14.53
95th Quantil;ALPHA;5.000;7640628;10.00;10.00;127297029;\textcolor{red}{\textit{127300617}};127307680;\textcolor{red}{\textit{127307744}};127321328;\textcolor{green}{\textbf{127313466}};2.73e-04;\textcolor{green}{\textbf{1.95e-04}};6.46e-04;\textcolor{red}{\textit{1.46e-03}};2.22e-03;\textcolor{red}{\textit{4.41e-03}};188223;\textcolor{green}{\textbf{153572}};245436;\textcolor{green}{\textbf{177181}};422117;\textcolor{green}{\textbf{267856}};0.109;0.181;0.328;nan;141.10;108.90;-2.61;1.25;3.37
\end{filecontents*}
\csvreader[
  longtable=lrrrrrr,
  table head=\caption{Direct comparison of our variant using the \tilted relaxation \ref{C2**:tiltedgain} and the original LKH with 5 ALPHA candidates in minimal cost and average time per run for all \textbf{large finished instances}.
\label{fig:CR_large_instances_ALPHA_direct}}\\
    \hline
 \multicolumn{1}{c}{Problem} & \multicolumn{2}{c}{CostMin} & %\multicolumn{2}{c}{CostAvg} & 
\multicolumn{2}{c}{TimeAvg in s} & \multicolumn{1}{c}{TimeAvg} \\ 
& LKH 3.0.8 & \ref{C2**:tiltedgain} & LKH 3.0.8 & \ref{C2**:tiltedgain} & \multicolumn{1}{c}{ratio in \%}\\ \hline \endfirsthead
\hline \multicolumn{1}{c}{Problem} & \multicolumn{2}{c}{CostMin} & %\multicolumn{2}{c}{CostAvg} & 
\multicolumn{2}{c}{TimeAvg in s} & \multicolumn{1}{c}{TimeAvg} \\   & LKH 3.0.8 & \ref{C2**:tiltedgain} & LKH 3.0.8 & \ref{C2**:tiltedgain} & \multicolumn{1}{c}{ratio in \%} \\ \hline \endhead
\hline
    \endfoot,
  late after line=\\,
head to column names, separator=semicolon]{Latex_Comparison_Large_ALPHA.csv}{}% use head of csv as column names
    {\Problem & \CostMinLKH & \CostMinOur & \TimeAvgLKH & \TimeAvgOur & \TimeAvgChange}
\newpage
\subsection{Medium instances}

\csvreader[
  longtable=lrrrrrr,
  table head=\caption{Direct comparison of our variant using the \tilted relaxation \ref{C2**:tiltedgain} and the original LKH with 5 POPMUSIC candidates in minimal cost and average time per run for all \textbf{medium instances}. Green and bold indicate a better performance of our variant, red and italic a worse performance of our variant.
\label{fig:CR_medium_instances_POPMUSIC_direct}}\\
    \hline
 \multicolumn{1}{c}{Problem} & \multicolumn{2}{c}{GapMin in \%} & %\multicolumn{2}{c}{CostAvg} & 
\multicolumn{2}{c}{TimeAvg in s} & \multicolumn{1}{c}{TimeAvg} \\ 
& LKH 3.0.8 & \ref{C2**:tiltedgain} & LKH 3.0.8 & \ref{C2**:tiltedgain} & \multicolumn{1}{c}{ratio in \%}\\ \hline \endfirsthead
\hline \multicolumn{1}{c}{Problem} & \multicolumn{2}{c}{GapMin in \%} & %\multicolumn{2}{c}{CostAvg} & 
\multicolumn{2}{c}{TimeAvg in s} & \multicolumn{1}{c}{TimeAvg} \\   & LKH 3.0.8 & \ref{C2**:tiltedgain} & LKH 3.0.8 & \ref{C2**:tiltedgain} & \multicolumn{1}{c}{ratio in \%} \\ \hline \endhead
\hline
    \endfoot,
    %table foot = \hline ,
  late after line=\\,
head to column names, separator=semicolon]{Latex_Comparison_Medium_POPMUSIC.csv}{}% use head of csv as column names
    {\Problem & \GapMinLKH & \GapMinOur & \TimeAvgLKH & \TimeAvgOur & \TimeAvgChange}
    
\begin{filecontents*}{Latex_Comparison_Medium_POPMUSIC_100.csv}
Problem;CandidateModel;Candidates;Optimum;RunsLKH;RunsOur;CostMinLKH;CostMinOur;CostAvgLKH;CostAvgOur;CostMaxLKH;CostMaxOur;GapMinLKH;GapMinOur;GapAvgLKH;GapAvgOur;GapMaxLKH;GapMaxOur;TimeMinLKH;TimeMinOur;TimeAvgLKH;TimeAvgOur;TimeMaxLKH;TimeMaxOur;CostMinChange;CostAvgChange;CostMaxChange;GapMinChange;GapAvgChange;GapMaxChange;TimeMinChange;TimeAvgChange;TimeMaxChange
Tnm\_instances\_Tnm2002;POPMUSIC;100;Error;100.0;100.0;37031143;\textcolor{green}{\textbf{37029600}};37031967;\textcolor{green}{\textbf{37030011}};37032674;\textcolor{green}{\textbf{37030111}};0.0;0.0;2.20e-05;\textcolor{green}{\textbf{1.10e-05}};4.10e-05;\textcolor{green}{\textbf{1.40e-05}};1.900;\textcolor{red}{\textit{3.400}};2.800;\textcolor{red}{\textit{7.450}};4.210;\textcolor{red}{\textit{16.54}};-4.17e-03;-5.28e-03;-6.92e-03;$\pm$ 0;-50.00;-65.85;+78.95;+166.07;+292.87
Tnm\_instances\_Tnm3001;POPMUSIC;100;Error;100.0;100.0;55942427;\textcolor{green}{\textbf{55939425}};55944300;\textcolor{red}{\textit{55944493}};55946519;\textcolor{green}{\textbf{55944885}};0.0;0.0;3.30e-05;\textcolor{red}{\textit{9.10e-05}};7.30e-05;\textcolor{red}{\textit{9.80e-05}};4.670;\textcolor{green}{\textbf{3.720}};8.850;\textcolor{red}{\textit{9.860}};23.10;\textcolor{green}{\textbf{15.97}};-5.37e-03;+3.45e-04;-2.92e-03;$\pm$ 0;+175.76;+34.25;-20.34;+11.41;-30.87
Tnm\_instances\_Tnm4000;POPMUSIC;100;Error;100.0;100.0;74861889;\textcolor{green}{\textbf{74858233}};74864990;\textcolor{green}{\textbf{74859282}};74865609;\textcolor{green}{\textbf{74859842}};0.0;0.0;4.10e-05;\textcolor{green}{\textbf{1.40e-05}};5e-05;\textcolor{green}{\textbf{2.10e-05}};6.110;\textcolor{red}{\textit{7.020}};10.62;\textcolor{red}{\textit{15.39}};15.88;\textcolor{red}{\textit{32.96}};-4.88e-03;-7.62e-03;-7.70e-03;$\pm$ 0;-65.85;-58.00;+14.89;+44.92;+107.56
Tnm\_instances\_Tnm5002;POPMUSIC;100;Error;100.0;100.0;93784081;\textcolor{red}{\textit{93784106}};93784119;\textcolor{red}{\textit{93787882}};93784125;\textcolor{red}{\textit{93788314}};0.0;0.0;0.0;\textcolor{red}{\textit{4e-05}};0.0;\textcolor{red}{\textit{4.50e-05}};10.84;\textcolor{red}{\textit{12.03}};16.68;\textcolor{red}{\textit{22.52}};19.98;\textcolor{red}{\textit{35.77}};+2.67e-05;+4.01e-03;+4.47e-03;$\pm$ 0;+0.00e+00;+0.00e+00;+10.98;+35.01;+79.03
\hline Average;POPMUSIC;100.0;Error;100.0;100.0;65404885;\textcolor{green}{\textbf{65402841}};65406344;\textcolor{green}{\textbf{65405417}};65407232;\textcolor{green}{\textbf{65405788}};0.0;0.0;2.40e-05;\textcolor{red}{\textit{3.90e-05}};4.10e-05;\textcolor{red}{\textit{4.45e-05}};5.880;\textcolor{red}{\textit{6.543}};9.737;\textcolor{red}{\textit{13.80}};15.79;\textcolor{red}{\textit{25.31}};-3.60e-03;-2.14e-03;-3.27e-03;nan;inf;inf;21.12;64.35;112.10
Median;POPMUSIC;100.0;Error;100.0;100.0;65402158;\textcolor{green}{\textbf{65398829}};65404645;\textcolor{green}{\textbf{65401888}};65406064;\textcolor{green}{\textbf{65402364}};0.0;0.0;2.75e-05;\textcolor{green}{\textbf{2.70e-05}};4.55e-05;\textcolor{green}{\textbf{3.30e-05}};5.390;\textcolor{green}{\textbf{5.370}};9.735;\textcolor{red}{\textit{12.62}};17.93;\textcolor{red}{\textit{24.75}};-4.53e-03;-2.47e-03;-4.92e-03;nan;62.88;-11.88;12.94;39.96;93.29
95th Quantil;POPMUSIC;100.0;Error;100.0;100.0;90945752;\textcolor{green}{\textbf{90945225}};90946250;\textcolor{red}{\textit{90948592}};90946348;\textcolor{red}{\textit{90949043}};0.0;0.0;3.98e-05;\textcolor{red}{\textit{8.33e-05}};6.95e-05;\textcolor{red}{\textit{9.00e-05}};10.13;\textcolor{red}{\textit{11.28}};15.77;\textcolor{red}{\textit{21.45}};22.63;\textcolor{red}{\textit{35.35}};-6.02e-04;3.46e-03;3.36e-03;nan;inf;inf;69.34;147.90;265.10
\end{filecontents*}
\csvreader[
  longtable=lrrrrrr,
  table head=\caption{Direct comparison of our variant using the \tilted relaxation \ref{C2**:tiltedgain} and the original LKH with 100 POPMUSIC candidates in minimal cost and average time per run for all \textbf{medium instances}.
\label{fig:CR_medium_instances_POPMUSIC_direct_100}}\\
    \hline
 \multicolumn{1}{c}{Problem} & \multicolumn{2}{c}{GapMin \%} & %\multicolumn{2}{c}{CostAvg} & 
\multicolumn{2}{c}{TimeAvg in s} & \multicolumn{1}{c}{TimeAvg} \\ 
& LKH 3.0.8 & \ref{C2**:tiltedgain} & LKH 3.0.8 & \ref{C2**:tiltedgain} & \multicolumn{1}{c}{ratio in \%}\\ \hline \endfirsthead
\hline \multicolumn{1}{c}{Problem} & \multicolumn{2}{c}{GapMin in \%} & %\multicolumn{2}{c}{CostAvg} & 
\multicolumn{2}{c}{TimeAvg in s} & \multicolumn{1}{c}{TimeAvg} \\   & LKH 3.0.8 & \ref{C2**:tiltedgain} & LKH 3.0.8 & \ref{C2**:tiltedgain} & \multicolumn{1}{c}{ratio in \%} \\ \hline \endhead
\hline
    \endfoot,
    %table foot = \hline ,
  late after line=\\,
head to column names, separator=semicolon]{Latex_Comparison_Medium_POPMUSIC_100.csv}{}% use head of csv as column names
    {\Problem & \GapMinLKH & \GapMinOur & \TimeAvgLKH & \TimeAvgOur & \TimeAvgChange}

\csvreader[
  longtable=lrrrrrr,
  table head=\caption{Direct comparison of our variant using the \tilted relaxation \ref{C2**:tiltedgain} and the original LKH with 5 ALPHA candidates in minimal cost and average time per run for all \textbf{medium instances}.
\label{fig:CR_medium_instances_ALPHA_direct}}\\
    \hline
 \multicolumn{1}{c}{Problem} & \multicolumn{2}{c}{GapMin in \%} & %\multicolumn{2}{c}{CostAvg} & 
\multicolumn{2}{c}{TimeAvg in s} & \multicolumn{1}{c}{TimeAvg} \\ 
& LKH 3.0.8 & \ref{C2**:tiltedgain} & LKH 3.0.8 & \ref{C2**:tiltedgain} & \multicolumn{1}{c}{ratio in \%}\\ \hline \endfirsthead
\hline \multicolumn{1}{c}{Problem} & \multicolumn{2}{c}{GapMin in \%} & %\multicolumn{2}{c}{CostAvg} & 
\multicolumn{2}{c}{TimeAvg in s} & \multicolumn{1}{c}{TimeAvg} \\   & LKH 3.0.8 & \ref{C2**:tiltedgain} & LKH 3.0.8 & \ref{C2**:tiltedgain} & \multicolumn{1}{c}{ratio in \%} \\ \hline \endhead
\hline
    \endfoot,
    %table foot = \hline ,
  late after line=\\,
head to column names, separator=semicolon]{Latex_Comparison_Medium_ALPHA.csv}{}% use head of csv as column names
    {\Problem & \GapMinLKH & \GapMinOur & \TimeAvgLKH & \TimeAvgOur & \TimeAvgChange}
\subsection{Small instances}

\csvreader[
  longtable=lrrrrrr,
  table head=\caption{Direct comparison of our variant using the \tilted relaxation \ref{C2**:tiltedgain} and the original LKH with 5 POPMUSIC candidates in minimal cost and average time per run for all \textbf{small instances}. Green and bold indicate a better performance of our variant, red and italic a worse performance of our variant.
\label{fig:CR_small_instances_POPMUSIC_direct}}\\
    \hline
 \multicolumn{1}{c}{Problem} & \multicolumn{2}{c}{GapMin in \%} & %\multicolumn{2}{c}{CostAvg} & 
\multicolumn{2}{c}{TimeAvg in s} & \multicolumn{1}{c}{TimeAvg} \\ 
& LKH 3.0.8 & \ref{C2**:tiltedgain} & LKH 3.0.8 & \ref{C2**:tiltedgain} & \multicolumn{1}{c}{ratio in \%}\\ \hline \endfirsthead
\hline \multicolumn{1}{c}{Problem} & \multicolumn{2}{c}{GapMin in \%} & %\multicolumn{2}{c}{CostAvg} & 
\multicolumn{2}{c}{TimeAvg in s} & \multicolumn{1}{c}{TimeAvg} \\   & LKH 3.0.8 & \ref{C2**:tiltedgain} & LKH 3.0.8 & \ref{C2**:tiltedgain} & \multicolumn{1}{c}{ratio in \%} \\ \hline \endhead
\hline
    \endfoot,
    %table foot = \hline ,
  late after line=\\,
head to column names, separator=semicolon]{Latex_Comparison_Small_POPMUSIC.csv}{}% use head of csv as column names
    {\Problem & \GapMinLKH & \GapMinOur & \TimeAvgLKH & \TimeAvgOur & \TimeAvgChange}

%\stepcounter{footnote}\footnotetext{As seen in \cref{fig:CR_small_instances_POPMUSIC_direct,%
%fig:CR_small_instances_POPMUSIC_direct_100,%
%fig:CR_small_instances_ALPHA_direct} the relative time difference between the original LKH and our modification is not measurable for all small instances. Therefore we do not specify the average, median or 95-th quantile for these instances.}


\csvreader[
  longtable=lrrrrrr,
  table head=\caption{Direct comparison of our variant using the \tilted relaxation \ref{C2**:tiltedgain} and the original LKH with 100 POPMUSIC candidates in minimal cost and average time per run for all \textbf{small instances}.
\label{fig:CR_small_instances_POPMUSIC_direct_100}}\\
    \hline
 \multicolumn{1}{c}{Problem} & \multicolumn{2}{c}{GapMin in \%} & %\multicolumn{2}{c}{CostAvg} & 
\multicolumn{2}{c}{TimeAvg in s} & \multicolumn{1}{c}{TimeAvg} \\ 
& LKH 3.0.8 & \ref{C2**:tiltedgain} & LKH 3.0.8 & \ref{C2**:tiltedgain} & \multicolumn{1}{c}{ratio in \%}\\ \hline \endfirsthead
\hline \multicolumn{1}{c}{Problem} & \multicolumn{2}{c}{GapMin in \%} & %\multicolumn{2}{c}{CostAvg} & 
\multicolumn{2}{c}{TimeAvg in s} & \multicolumn{1}{c}{TimeAvg} \\   & LKH 3.0.8 & \ref{C2**:tiltedgain} & LKH 3.0.8 & \ref{C2**:tiltedgain} & \multicolumn{1}{c}{ratio in \%} \\ \hline \endhead
\hline
    \endfoot,
    %table foot = \hline ,
  late after line=\\,
head to column names, separator=semicolon]{Latex_Comparison_Small_POPMUSIC_100.csv}{}% use head of csv as column names
    {\Problem & \GapMinLKH & \GapMinOur & \TimeAvgLKH & \TimeAvgOur & \TimeAvgChange}
    

\csvreader[
  longtable=lrrrrrr,
  table head=\caption{Direct comparison of our variant using the \tilted relaxation \ref{C2**:tiltedgain} and the original LKH with 5 ALPHA candidates in minimal cost and average time per run for all \textbf{small instances}.
\label{fig:CR_small_instances_ALPHA_direct}}\\
    \hline
 \multicolumn{1}{c}{Problem} & \multicolumn{2}{c}{GapMin in \%} & %\multicolumn{2}{c}{CostAvg} & 
\multicolumn{2}{c}{TimeAvg in s} & \multicolumn{1}{c}{TimeAvg} \\ 
& LKH 3.0.8 & \ref{C2**:tiltedgain} & LKH 3.0.8 & \ref{C2**:tiltedgain} & \multicolumn{1}{c}{ratio in \%}\\ \hline \endfirsthead
\hline \multicolumn{1}{c}{Problem} & \multicolumn{2}{c}{GapMin in \%} & %\multicolumn{2}{c}{CostAvg} & 
\multicolumn{2}{c}{TimeAvg in s} & \multicolumn{1}{c}{TimeAvg} \\   & LKH 3.0.8 & \ref{C2**:tiltedgain} & LKH 3.0.8 & \ref{C2**:tiltedgain} & \multicolumn{1}{c}{ratio in \%} \\ \hline \endhead
\hline
    \endfoot,
    %table foot = \hline ,
  late after line=\\,
head to column names, separator=semicolon]{Latex_Comparison_Small_ALPHA.csv}{}% use head of csv as column names
    {\Problem & \GapMinLKH & \GapMinOur & \TimeAvgLKH & \TimeAvgOur & \TimeAvgChange}

\section{Code}
\begin{figure}[h] 
\lstset{
frame=single,
language=C,
keywordstyle=\color{blue},
commentstyle=\color{green},
numbers=left,
basicstyle=\ttfamily\small%Added this to make it fit in the frame
}
\begin{lstlisting}
for (Nt2 = t2->CandidateSet; (t3 = Nt2->To); Nt2++) {
    if (t3 == t2->Pred || t3 == t2->Suc) {
        continue;
    } else if ((G1 = *G0 - Nt2->Cost) <= 0 && GainCriterionUsed && 
        ProblemType != HCP && ProblemType != HPP) {
        if ((G1 = *G0 - Nt2->Cost) <= 0 && GainCriterionPartial) {
            if (GainCriterionViolated  == 1) {
                continue;
            } else
                GainCriterionViolated = 1;
            }
        } else if ((G1 = *G0 - Nt2->Cost) <= 0) {
            continue; 
        } else {
            continue;
        }
    /*Reset of GainCriterionViolated*/
    } else if ((G1 = *G0 - Nt2->Cost) > 0 && GainCriterionUsed &&
        GainCriterionPartial && 
        ProblemType != HCP && ProblemType != HPP) {
        if (GainCriterionViolated) {
            GainCriterionViolated = 0;
        }
    }
    <more code>
}
\end{lstlisting}
\caption{Excerpt of function \texttt{Node *Best5OptMove(Node * t1, Node * t2, GainType * G0, GainType * Gain)} in Best5OptMove.c using the \homogeneous relaxation \ref{C2*:newgain}.}
\label{lst:CodeBest5opt}
\end{figure}
In \cref{lst:CodeBest5opt} we display an excerpt the adjusted function \texttt{Node *Best5OptMove(Node * t1, Node * t2, GainType * G0, GainType * Gain)} in Best5OptMove.c using the \homogeneous relaxation \ref{C2*:newgain}. 
The excerpt gives an idea of our crucial changes in the implementation. 
In contrast to the pseudocode in \cref{alg:LKH_0.1_Gain_Criterion}, the original implementation by Helsgaun \cite{LKH3} does not provide access to the value of the previous gain $G_{i-1}$ in every step.
Therefore, we introduce a boolean variable \struc{GainCriterionViolated} which marks, whether $G_{i-1} < 0$ is true or not.
To that end, we initialize this variable in the files \enquote{LKH.h} and \enquote{ReadParameters.c}. 
If the positive gain criterion \ref{C2:gain} was not violated in the last step but is not satisfied in the current step, we update the variable GainCriterionViolated to \enquote{true}. 
Otherwise, if the positive gain criterion was violated in the last step but is satisfied in the current step we update GainCriterionViolated to \enquote{false}. 
We implement this update in \enquote{BestKOptMove.c} and \enquote{Best5KOptMove.c}. 

\end{document}